\documentclass[12pt,reqno]{amsart}

\usepackage{geometry}
\usepackage{fullpage}
\usepackage[all]{xy}

\usepackage[OT2,T1]{fontenc}
\usepackage{textcomp}

\usepackage[scaled=0.92]{helvet}
\usepackage{newtxmath}
\usepackage[english]{babel}
\usepackage{booktabs}

\usepackage{etoolbox}
\patchcmd{\section}{\scshape}{\bfseries}{}{}
\makeatletter
\renewcommand{\@secnumfont}{\bfseries}
\makeatother

\newtheorem{introtheorem}{Theorem}

\newtheorem{introconjecture}[introtheorem]{Conjecture}
\theoremstyle{definition}

\newtheorem{introdef}[introtheorem]{Definition}

\theoremstyle{plain}
\newtheorem{theorem}{Theorem}[subsection]
\newtheorem{proposition}[theorem]{Proposition}
\newtheorem{lemma}[theorem]{Lemma}

\newtheorem{conjecture}[theorem]{Conjecture}
\theoremstyle{definition}
\newtheorem{definition}[theorem]{Definition}
\newtheorem{notation}[theorem]{Notation}

\newtheorem*{examplesn}{Examples}
\newtheorem{remark}[theorem]{Remark}

\newtheorem{example}[theorem]{Example}

\counterwithin{equation}{section}

\usepackage[usenames, dvipsnames, table]{xcolor}
\usepackage[hidelinks, hypertexnames=false]{hyperref}
\usepackage{mathtools}
\usepackage{bigints} 
\usepackage{calrsfs}

\usepackage{subcaption} 

\usepackage{listings}
\definecolor{codegreen}{rgb}{0,0.6,0}
\definecolor{codegray}{rgb}{0.5,0.5,0.5}
\definecolor{codepurple}{rgb}{0.58,0,0.82}
\definecolor{backcolour}{RGB}{240,240,240}
\lstdefinestyle{mystyle}{
    commentstyle=\color{codegreen},
    keywordstyle=\color{magenta},
    numberstyle=\tiny\color{codegray},
    stringstyle=\color{codepurple},
    basicstyle=\ttfamily\footnotesize,
    breakatwhitespace=false,         
    breaklines=true,                 
    captionpos=b,                    
    keepspaces=true,                                 
    numbersep=5pt,                  
    showspaces=false,                
    showstringspaces=false,
    showtabs=false,                  
    tabsize=2
}
\newcommand{\mf}{\mathfrak}
\newcommand{\mc}{\mathcal}
\newcommand{\Z}{\mathbf{Z}}
\newcommand{\R}{\mathbf{R}}
\newcommand{\Q}{\mathbf{Q}}
\newcommand{\FF}{\mathbf{F}}
\newcommand{\Cc}{\mathbf{C}}
\newcommand{\T}{\mathbf{T}}
\newcommand{\PP}{\mathbf{P}}
\newcommand\h{\mathbf{h}}
\newcommand\x{\mathbf{x}}

\newcommand\m{\mathrm{m}}
\newcommand{\mfq}{\mathfrak{q}}
\newcommand{\MM}{\mathrm{M}}
\newcommand\NN{{\mathcal N}}

\newcommand\vtheta{\boldsymbol{\theta}}
\newcommand{\Gm}{\mathbf{G}_m}
\newcommand{\sym}{\mathfrak{S}}
\newcommand{\err}{\mathcal{E}}
\DeclareMathOperator{\tr}{tr}
\renewcommand{\geq}{\geqslant}
\renewcommand{\leq}{\leqslant}
\DeclareMathOperator{\Gal}{Gal}
\renewcommand{\Re}{\mathrm{Re}}
\renewcommand{\d}{\mathrm{d}}
\newcommand{\hypgeo}[2]{
  {\vphantom{F}}_{#1}\kern-\scriptspace F_{#2}
}
\DeclarePairedDelimiter\abs{\lvert}{\rvert}

\usepackage{tikz}
\def\nudge{.5}
\tikzset{axis/.style={thick, Black!75!black, -latex, shorten <=-\nudge cm, shorten >=-2*\nudge cm}}
\tikzset{line/.style={ultra thick,Black}}

\begin{document}

\date{\today}
\title{The asymptotic Mahler measure of Gaussian periods}
\author[G.~Cornelissen]{Gunther Cornelissen}
\address{\parbox{\linewidth}{\normalfont Mathematisch Instituut, Universiteit Utrecht, Postbus 80.010, 3508 TA Utrecht, Nederland \\
Max-Planck-Institut f\"ur Mathematik, Postfach 7280, 53072 Bonn, Deutschland \\ \vspace*{-3mm} \mbox{ } }}
\email{g.cornelissen@uu.nl}
\author[D.~Hokken]{David Hokken}
\address{\normalfont Mathematisch Instituut, Universiteit Utrecht, Postbus 80.010, 3508 TA Utrecht, Nederland}
\email{d.p.t.hokken@uu.nl}
\author[B.~Ringeling]{Berend Ringeling}
\address{\normalfont CRM / Universit\'e de Montr\'eal, P.O. Box 6128, Centre-ville Station
Montr\'eal (Qu\'ebec) H3C 3J7, Canada}
\email{b.j.ringeling@gmail.com}

\thanks{We thank Francesco Amoroso, Frits Beukers, Peter Koymans, Matilde Lal\'in, Junxian Li, Riccardo Pengo, Subham Roy, Jan Stienstra, Seth Sullivant, Th\'eo Untrau, Don Zagier and Wadim Zudilin for valuable input. We are also grateful to the anonymous referee for various suggestions. G.~C. thanks the Max-Planck-Institut f\"ur Mathematik (Bonn) for excellent conditions during work on this paper. D.~H. is supported by the Dutch Research Council (NWO) through the grant OCENW.M20.233.}  

\begin{abstract}
We construct a sequence of cyclotomic integers (Gaussian periods) of particularly small Mahler measure/height. We study the asymptotics of their Mahler measure as a function of their conductor, to find that the growth rate is the (multivariate) Mahler measure of a family of log Calabi--Yau varieties of increasing dimension. In turn, we study the asymptotics of some of these Mahler measures as the dimension increases, as well as properties of the associated algebraic dynamical system. We describe computational experiments that suggest that these cyclotomic integers realize the smallest non-zero logarithmic Mahler measure in the set of algebraic integers with cyclic Galois group of a given odd order. Finally, we discuss some precise conjectures that imply double logarithmic growth for those Mahler measures as a function of that order. 

The proofs use ideas from the theory of quantitative equidistribution, reflexive polytopes and toric varieties, the theory of random walks, Bessel functions, class field theory, and Linnik's constant.
\end{abstract}

\maketitle

\section{Introduction}

How `small' can a number be, given some algebraic constraints on that number? We investigate a question of the said type, namely: how small can the logarithmic Mahler measure $\m(\alpha)$ of an algebraic integer $\alpha$ ---not itself a root of unity--- be, if we assume that $\alpha$ generates a field with cyclic Galois group of order $n>1$? The current best lower bound for $\m(\alpha)$ is linear in $n$ (\cite{Schinzel,AmorosoDvornicich}). After some numerical experiments, we identified an interesting sequence of algebraic integers (in fact, Gaussian periods) whose Mahler measure might be the true minimum for odd $n$; hence we study their Mahler measure asymptotically in $n$, taking us into the world of reflexive polytopes, Calabi--Yau varieties and random walks. Using this detour through geometry, we discuss the growth of the Mahler measures of these Gaussian periods, which we believe to be typically of order $n \log{\log{n}}$. By the relation between Mahler measure and height, soon to be spelled out in detail, all our results translate into statements on minimal heights of such algebraic integers. 
  
\subsection*{Definition of a sequence of cyclotomic integers} 
For a positive integer $n$, denote by $\kappa(n)$ the smallest positive integer such that $p = p(n) \coloneqq \kappa(n) n + 1$ is prime (sequence \cite[A034693]{oeis} in the OEIS). 
Note that $\kappa(n)$ is even when $n$ is odd. See Table \ref{table-k} for some values. 

\begin{table}[h]
\begin{center}
\begin{tabular}{ccccccccccccc}
\toprule
$n$  & 7 & 13 & 17 & 19 & 24 & 25 & 27 & 31 & 32 & 34 & 37 & 38 \\ \midrule
$\kappa(n)$ & 4 & 4 & 6 & 10 & 3 & 4 & 4 & 10 & 3 & 3 & 4 & 5 \\
$p$ & 29 & 53 & 103 & 191 & 73 & 101 & 109 & 311 &  97 & 103 & 149 & 191  \\ \bottomrule
\end{tabular}
\end{center} 
\caption{Values of $n\leq 40$ with $\kappa(n) \geq 3$ and $p=kn+1$.}   
\label{table-k}  
\end{table}

\begin{figure}
 $$ \xymatrix@C=2mm@R=4mm{  \Q(\zeta_p) \ar@{-}[dd]_{G=(\Z/p)^*} \ar@{-}[dr]^H  \\ &  \Q(\alpha_n) \ar@{-}[dl]^{\Z/n\Z}    \\ \Q  &  } $$
\caption{A diagram of field extensions}  \label{fields} 
\end{figure}
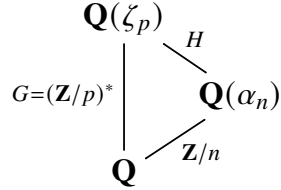

For fixed $k \in \Z_{\geq 1}$, define
 \begin{equation} 
 	\label{defnnk} 
 	\NN_k \coloneqq  \{ n \colon \kappa(n)=k \}. 
 \end{equation} 
 We prove that $\NN_k$ is infinite for any $k$, cf.\ \S \ref{asNNk}. For now, fix $k$, assume that $n \in \NN_k$, set $p=kn+1$, and let $\zeta_p$ denote a primitive $p$-th root of unity. The Galois group of the cyclotomic field $\Q(\zeta_p)$ is cyclic of order $p-1 = k n$; let $H$ denote its unique subgroup of order $k$ and define the \emph{Gaussian period} 
$$ \alpha_n \coloneqq \tr_H(\zeta_p) \coloneqq  \sum\limits_{\sigma \in H} \zeta_p^{\sigma}. $$
Standard theory of cyclotomic fields (see, e.g., \cite[II]{Marcus}) implies that  $\Q(\alpha_n)/\Q$ is a cyclic Galois extension of order $n$, with $\alpha_n$ as a normal generator and $p$ as conductor, see Figure \ref{fields}. Also notice that for $k>1$ (i.e., $n+1$ composite), $\alpha_n$ is itself not a root of unity. These Gaussian periods $\alpha_n$ are the algebraic integers alluded to above.  

\begin{examplesn} \mbox{ } 
\begin{enumerate}
\item For $k = 2$ (so $2n+1$ is prime but $n+1$ is not) we have $\alpha_n=2 \cos(2 \pi/p)$. 
\item For $k=4$, we have $\alpha_n = 2 \cos(2\pi/p) + 2 \cos(2\pi\nu/p)$ for any $\nu \in (\Z/p)^*$ of order $4$. For example, if $n=7$, then $p=29$ and we may take $\nu=12$.
\end{enumerate}
\end{examplesn}

\subsection*{Mahler measure and height} 

The zero locus of a Laurent polynomial $f \in \Cc[x_1^{\pm 1},\dots,x_N^{\pm 1}]$ defines a (Laurent) hypersurface $V$ in $\Gm^N$, and the (logarithmic, multivariate) \emph{Mahler measure} of $f$ or $V$ is defined by 
\begin{align*} \m(f) = \m(V) & \coloneqq  \int_0^1 \dots \int_0^1 \log |f(e^{2 \pi i \theta_1},\dots,e^{2 \pi i \theta_N}) | \, \d\theta_1 \dots \d\theta_N \\ &= \int_{\T^N} \log|f(x_1,\dots,x_N)| \, \d\mu, \end{align*} 
where $\T^N =\{ (x_1,\dots,x_N) \in \Cc^N \colon |x_i|=1, \forall i\}$ is the complex unit torus with normalized Haar measure $\d\mu = \frac{1}{(2\pi i)^N} \frac{\d X_1}{X_1} \cdots \frac{\d X_N}{X_N}$. The \emph{exponential Mahler measure} is defined by $\MM(f)=\MM(V)\coloneqq \exp(\m(f))$.

We denote the ring of algebraic integers by $\overline \Z$. 
By definition, the Mahler measure $\m(\alpha)$ of an algebraic integer $\alpha \in \overline \Z$ with minimal polynomial $f_\alpha$ over $\Q$ is $\m(\alpha)\coloneqq \m(f_\alpha)$. Write $$\log^+{x} = \max \{0, \log{x}\}.$$ By Jensen's formula \cite[Lemma 1.9]{EverestWard},
\begin{equation} \label{mahlerheight} \m(\alpha) = \sum \log^+ |\alpha'|  = [\Q(\alpha):\Q]\, h(\alpha), \end{equation} 
where the sum is over all conjugates $\alpha'$ of $\alpha$ (the roots of $f_\alpha$ in $\overline \Q$), and $h(\alpha)$ is the \emph{absolute logarithmic Weil height} of $\alpha$. Thus, $\MM(\alpha)\coloneqq \MM(f_\alpha)$ is the product of the absolute value of all conjugates of $\alpha$ outside the complex unit circle. (The multivariable Mahler measure is also a height, related to a line bundle on the projective toric variety defined by the Newton polygon of $f$, see \cite[\S 7.2--7.3]{Maillot}.) 

\subsection*{The cyclopolytope and the cyclovariety} 

For an integer $N$, let $\mu_N$ denote the set of (not necessarily primitive) $N$-th roots of unity, and denote by $\mu_\infty \coloneqq  \bigcup \mu_N$ the set of all roots of unity. 

\begin{introdef}
\label{def:A}
Fix an integer $k>0$. Any element of $\mu_k$ can be expressed as an integer linear combination of $\{1,\zeta_k,\dots,\zeta_k^{d-1}\}$ with $d=\phi(k)$; let $A_k$ denote the corresponding set of vectors in $\R^d$. 
\begin{enumerate}
\item The \emph{cyclopolytope} $N_k$ is the convex hull of $A_k$ in $\R^d$. \item The \emph{cyclovariety} $C_k$ is the Laurent hypersurface  in $\Gm^{d+1}$ of dimension $d$ defined by the equation 
\begin{equation} \label{defpk} P_k(x_0,\mathbf{x}) \coloneqq  x_0 + F_k(\mathbf{x})=0  \mbox{ with } F_k(\mathbf{x}) \coloneqq  \sum_{\mathbf a \in A_k}  \mathbf{x}^{\mathbf{a}} \end{equation} 
where $\mathbf{x} = (x_1,\dots,x_d)$ is a multi-index. Note that $N_k$ is the Newton polytope of $F_k$. 
\end{enumerate}

\end{introdef}  

\begin{examplesn}[cyclopolytopes]
\mbox{ } 
\begin{enumerate}
\item $\mu_4 = \{ \pm 1, \pm i\}$ gives the four vectors $A_4=\{\pm(1,0),\pm(0,1)\}$ in the basis $\{ 1,i\}$. 
\item Writing $\rho$ for a primitive third root of unity, $\mu_6 = \{ \pm 1, \pm \rho, \pm \rho^2 = \mp (\rho + 1)\}$ gives the six vectors $A_6=\{\pm(1,0),\pm(0,1),\pm(1,-1)\}$ in the basis $\{1,-\rho\}$. 
\end{enumerate} 
Figure \ref{np23} shows the two corresponding cyclopolytopes. 

\begin{figure} 
\scalebox{0.8}{\begin{tikzpicture}
\draw[axis] (-1.1,0) -- (1.1,0) node[right=2* \nudge cm] {\(x_1\)};
\draw[axis] (0,-1.1) -- (0,1.1) node[above=2*\nudge cm] {\(x_2\)};
\draw[line] (-1,0) -- (0,1);
\draw[line] (0,1) -- (1,0);
\draw[line] (1,0) -- (0,-1);
\draw[line] (0,-1) -- (-1,0);
\end{tikzpicture}}
\, \, \, \, \, \, \, 
\scalebox{0.8}{\begin{tikzpicture}
\draw[axis] (-1.1,0) -- (1.1,0) node[right=2* \nudge cm] {\(x_1\)};
\draw[axis] (0,-1.1) -- (0,1.1) node[above=2*\nudge cm] {\(x_2\)};
\draw[line] (-1,0) -- (-1,1);
\draw[line] (-1,1) -- (0,1);
\draw[line] (0,1) -- (1,0);
\draw[line] (1,0) -- (1,-1);
\draw[line] (1,-1) -- (0,-1);
\draw[line] (0,-1) -- (-1,0);
\end{tikzpicture}}\\[2mm]
\caption{The cyclopolytopes (polygons) $N_4$ (left) and $N_6$ (right)} \label{np23}
\end{figure}
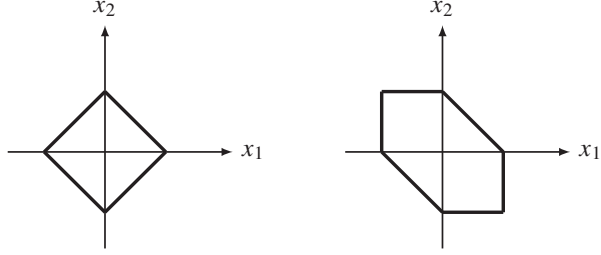 

\end{examplesn} 

\begin{examplesn}[cyclovarieties, cf.\ \S \ref{cvarcpol}]
\label{excl} 
\mbox{ }
\begin{enumerate}
\item $C_1$ is defined by $x_0+x_1=0$; notice that $\m(C_1)=0$. 
\item For $k=2^r$ with $r \geq 1$, $C_{k}$ is defined by $x_0+x_1+1/x_1+x_2+1/x_2+\dots+x_{k/2}+1/x_{k/2}=0$. 
\item For $k=q$ with $q$ an odd prime, $C_k$ is defined by $x_0+x_1+\dots+x_{k-1}+1/(x_1 \cdots x_{k-1})=0$.
\item $C_6$ is defined by $x_0+x_1+1/x_1+x_2 + 1/x_2 + x_1/x_2+x_2/x_1=0$. 
\end{enumerate} 
In Table \ref{cycth}, we list the Mahler measures of some of these cyclovarieties. 
 \begin{table}[ht!]
\centering
\begin{tabular}{ll}
\toprule
Value & Reference \\ 
\midrule
$\m(C_2)=L'(\chi_{-3},-1)$ & $\chi_{-3}$ odd character of conductor $3$; Smyth \cite{Smyth-C2} \\
$\m(C_3)=- 14 \zeta'(-2)$	 & Smyth \cite{Smyth-C2}\\
$\m(C_4)=-2L'(E,-1)$  & $E$ elliptic curve of conductor $15$ with Cremona label 15a8;  \\
&  conjecture of Boyd and Rodriguez-Villegas \cite{Boyd-RV-conj};\\
& proven by Brunault \cite{Brunault} \\[-0.7mm]
$\m(C_5)\overset{?}{=}-8 L'(f,-1)$ & $f(\tau) \coloneqq  (\eta(\tau) \eta(2\tau) \eta(3\tau) \eta(6\tau))^2 \in M_4(\Gamma_0(6))$; \\
& conjecture of Rodriguez-Villegas (2003); see, e.g., \cite[p.~75]{BrunaultZudilin}   \\
$\m(C_6)=0.6439432099\dots$ &  more digits computed in \S \ref{MC6}; no conjectured closed expression\\
\bottomrule
\end{tabular}
\caption{Known and conjectured values $\m(C_k)$ for $k \leq 6$}\label{cycth}
\end{table}
To switch between $C_4$ and the variety of Brunault, $z+(x+1)(y+1)=0$, change variables by $x=x_1 x_2, y=x_1/x_2, z=x_0x_1$ and divide by $x_1$; to switch between $C_5$ and the variety of Rodriguez-Villegas, $1+y_1+\dots+y_5=0$, change variables by $y_i=x_i/x_0$ for $i=1,\dots,4$, $y_5 = 1/(x_0x_1\cdots x_4)$ and multiply by $x_0$.

\end{examplesn} 

\subsection*{Asymptotic Mahler measure: fixing $k$} 

Our first main result is that the asymptotic growth of $\m(\alpha_n)$ in $n$ is well-defined by the above setup of \emph{first} fixing $k$ and then considering $n$ for which $\kappa(n)=k$. Since $n$ and $p$ are linearly related for fixed $k$, the result is not just asymptotic in $n \in \NN_k$ for fixed $k$, but also asymptotic in the conductor $p \in k \NN_k+1$ for fixed $k$. 

\begin{introtheorem} \label{mainas} For $n \in \NN_1$, $\m(\alpha_n)=0$, which equals $\m(C_1)=0$. 
For fixed $k \in \Z_{\geq 2}$, 
\begin{equation} \label{nowas} \m(\alpha_n) \sim  \m(C_k) \cdot n \mbox{ as } n \rightarrow + \infty, n \in \NN_k.\end{equation}
 In fact, for all $k \in \Z_{\geq 1}$ and $n \in \NN_k$,  
\begin{equation} \label{was} \Big| \frac{\m(\alpha_n)}{n} - \m(C_k) \Big| \leq\frac{4\sqrt{3} k(k+1)}{p^{1/d}} + \frac{2\log k}{p},    
\end{equation} 
where, as before, $p=kn+1$ and $d=\phi(k)$. 
 \end{introtheorem}

The proof of Theorem \ref{mainas} is given in Section \ref{proof1}. The basic idea is that (up to some manipulation), $\m(\alpha_n)/n$ is a sum of a function evaluated at sample points or integrated against a delta-measure, and $\m(C_k)$ is a corresponding integral against an absolutely continuous measure. This allows us to interpret the result as a comparison statement, either between Riemann sums and Riemann integrals, or between integrals over different measures. To prove \eqref{was}, we use the second approach applied to a one-variable function, relying on a recent method of Kowalski and Untrau \cite{KU} involving the Wasserstein distance. We also discuss a more standard approach using quasi-Monte-Carlo integration techniques for a multivariable function, where Weyl equidistribution (in this case proven by Myerson \cite{Myerson}) immediately proves \eqref{nowas}, and  the Hlawka--Koksma inequality gives an upper bound that captures the $k$-dependence in the Hardy--Krause variation (that we show to be finite but exponential in $k$), and the $p$-dependence by bounding the discrepancy of the set of sample points.  For $k=2$, we determine some lower order terms by Euler--Maclaurin summation, see Proposition \ref{LOT}. 

\subsection*{Geometry of the cyclovariety} 

In Section \ref{sec:geometry of the cyclovariety}, we study the cyclopolytope and the geometry of the cyclovariety in greater detail, focussing on the case where $k$ is even. For $k=4$ and $k=6$, considering $x_0$ as a parameter, $C_k$ compactifies to a family of elliptic curves over $\PP^1_{x_0}$ with finitely many singular fibers. This generalizes to the following statement about other cyclovarieties. 

\begin{introtheorem} \label{mainCY}
For $k=2 q^r$ with $q$ an odd prime, the cyclopolytope $N_k$ \textup{(}which is the Newton polytope of $C_k \cap \{x_0=0\}$\textup{)} is centrally symmetric and reflexive, and defines a smooth projective toric Fano variety $\PP_k$ which is a direct product of toric del Pezzo varieties. 

The Zariski closure of $C_k(x_0)$ in $\PP_k$, parametrized by $x_0 \in \PP^1$, is a pencil in the anticanonical linear system $|-K_{\PP_k}|$, and after resolving the base locus \textup{(}that is contained in the toric boundary of $\PP_k$\textup{)}, $C_k(x_0)$ admits a compactification to a family, all but finitely many of whose fibers are smooth varieties with canonical bundle the restriction of the exceptional divisor of the base point resolution. 
\end{introtheorem}

The members $C_k(x_0)$ of the original family for which the compactification is smooth are `log Calabi--Yau varieties' in the sense of Green, Hacking and Keel \cite[Def.\ 1.1]{GHK}, but since there are so many versions of the notion `Calabi--Yau', we have chosen to spell out the exact properties in Theorem \ref{mainCY}. In Proposition \ref{cohCk}, we further compute the cohomology groups of the structure sheaf on the smooth fibers of the compactified family and find that they are vanishing in all middle degrees, the Euler characteristic is zero, and the varieties are simply connected for $k \neq 6$, like for a `strict' Calabi--Yau variety (except that their canonical bundle is not trivial, due to resolving the base locus).

In Section \ref{NP}, we also elaborate on the more general case, and study in greater detail the polytope $N_k$ in relation to Batyrev's work on Mirror Symmetry \cite{Batyrev}.  The cyclopolytope has been studied before in \cite{BeckHosten}, and the associated toric varieties were first studied by Klyachko and Voskresenkii \cite{KlyachkoVoskresenskii}. In  condensed matter physics, the varieties with $x_0$ as variable (playing the role of energy level) are known as `Bloch varieties' and their fibers for fixed $x_0$ are called `Fermi varieties'. 

\subsection*{The associated algebraic dynamical system} 

By work of Schmidt, Lind and Ward \cite{SchmidtLindWard}, the cyclovariety corresponds to an algebraic dynamical system where $\Z^{d+1}$ acts on a topological group with (topological) entropy given by $\m(\alpha_n)$: the space is the Pontryagin dual of the additive group of the coordinate ring of $C_k$ and the action is dual to multiplication by monomials, cf.\ 
\cite{Schmidt}. In Section \ref{DS}, we prove the following.

\begin{introtheorem} \label{mainDyn} 
The algebraic dynamical system corresponding to the cyclovariety $C_k$ for $k \geq 2$ has the following properties: 
\begin{enumerate} 
\item It has positive entropy $\m(C_k$), is ergodic, mixing of all orders and Bernoulli. 
\item For $k\in \{2^r, q^r, 2q^r\}$ with $q$ an odd prime and $r$ arbitrary, it is not expansive. 
\end{enumerate} 
\end{introtheorem}

The proof uses the work of Schmidt and others to relate dynamical properties to arithmetic-geometric properties of the variety that are also of independent interest. The theorem follows from studying, for example, generalized cyclotomic points on the cyclovariety, or points on the intersection of cyclovariety and the complex torus; see Proposition \ref{toricpointsDS} where we prove that \begin{equation} \label{introinter} C_k(\Cc) \cap \T^{d+1} \neq \emptyset \end{equation} for $k$ as in the second part of the theorem. We believe that non-expansiveness always holds, but our proof breaks down if $\phi(k)/k$ becomes small, see Remark \ref{remprimorial}.  

\subsection*{Asymptotics of $\m(C_k)$ in $k$} 
For $k 
\geq 2$, $\m(C_k)$ is bounded away from zero, see Proposition \ref{infmck}; this follows from Smyth's lower bound for Mahler measures of non-reciprocal one-variable polynomials together with Lawton's theorem relating single- and multivalued Mahler measures through a limiting process. Concerning upper bounds, in Proposition~\ref{prop:upper-bound}, we use $L^2$-norms to show that  $\m(C_k) \leq {\log(k+1)}/{2}$ for general $k$. 

Next, we study the asymptotics of $\m(C_{k})$ as a function of $k$; we were unable to do this for general $k$, but  we have the following result, where $\gamma$ is the Euler--Mascheroni constant. 

\begin{introtheorem} 
\label{mainRW} 
For $k$ of the form $k=2q^r$ with fixed $r \geq 1$ and a prime number $q$ tending to infinity, or $q=2$ and $r$ tending to infinity, we have 
\begin{equation}
\m(C_k) = \frac{\log k}{2} - \frac{\gamma}{2} - \frac{\log 2 }{2} + \sqrt{\frac{2}{\pi k }} + \frac{1}{4 k} + O\left(\frac{\log k }{k^{3/2}} \right).
\end{equation}
For $k=q$ a prime number tending to infinity, we have 
\begin{equation}
\m(C_k) = \frac{\log k}{2} - \frac{\gamma}{2} + \frac{5}{8 k} + O\left(\frac{\log k }{k^{3/2}} \right).
\end{equation}
\end{introtheorem}
The result is surprisingly close to the previously mentioned general upper bound. 
The proof allows us to determine higher order approximations of $\m(C_k)$; see \eqref{moreterms}.  
The method of proof is a mix of probability theory and the theory of special functions. The Mahler measure is the integral of the logarithm against the probability density function $\rho_k$ of the random variable $|F_k(X_1,\dots,X_d)|$. Such an approach has been used for  `linear' Mahler measures, i.e., of the form $\m(x_1+\dots+x_d)$, which can be understood in terms of random planar walks, as studied by Rayleigh, Pearson, and others, see, e.g., \cite{ProfWalk}; the asymptotics for linear Mahler measures as $d$ tends to infinity can then be obtained via the central limit theorem, see, e.g., \cite{MyersonSmyth}. It seems we cannot use versions of the central limit theorem directly in our case, due to the dependency of the random variables being summed; dependencies that are rather intricate and depend on the structure of the cyclopolytope. Rather, after studying the existence of $\rho_k$, we adopt an old method of Kluyver \cite{kluyver1906local} that expresses the corresponding probability density function $\rho_k$ in terms of Bessel functions (Proposition 
\ref{prop:gtrepr}), and then use asymptotic results for those functions, thus avoiding any statistics with dependent variables. We work this out in detail for the most complicated case, $k=2q^r$ with $q$ tending to infinity and $r$ fixed, where an infinite sum of Bessel functions appears. The other cases can be more easily reduced to the known case of linear Mahler measures, see \S \ref{subsecas}. 

\subsection*{An algorithm to compute Mahler measure} This leads us to discuss, in Section \ref{MC6}, an algorithm to compute the probability density functions $\rho_k$ efficiently, by using a Picard--Fuchs equation satisfied by the constant term sequence of powers of $F_k$, and applying Hilbert transforms and analytic continuation. The approach also provides a fast algorithm to compute the Mahler measure accurately (up to hundreds of digits); to do this directly from the definition is numerically very unstable, since, in general, the integrand has singularities along the intersection of the cyclovariety and the torus, cf.\ Formula \eqref{introinter}.   The practicality of the algorithm hinges on the computation of the differential equation, which can be done in theory by using the creative telescoping algorithm or the Griffiths--Dwork method, and in practice for values of $k$ roughly up to $10$. 

\subsection*{A Lehmer-style question for cyclotomic integers, and a computer search using class field theory} 
Lehmer \cite{Lehmer} asked in 1933 what is the value of $$ \inf \{ \m(\alpha) \colon \alpha \in \overline \Z - \mu_\infty \}.$$ The current most common belief is that the infimum equals a minimum and is bounded away from zero; more precisely, that the infimum equals $\m(x^{10}+x^{9}-x^{7}-x^{6}-x^{5}-x^{4}-x^{3}+x+1)$. 

Fixing $n \in \Z_{\geq 2}$, we can consider the best lower bound for the Mahler measure of cyclotomic integers (themselves not roots of unity), with cyclic Galois group of order $n$: 
\begin{equation}
\label{eq:m_n-t_n-definition}
\begin{aligned} 
m_n & \coloneqq  \min \{ \m(\alpha) \colon \alpha \in \overline \Z - \mu_\infty, \Gal(\Q(\alpha)/\Q) \cong \Z/n\Z \}.
\end{aligned}
\end{equation}
This is a minimum, not just an infimum, since there are only finitely many monic polynomials of a fixed degree $n$ with an upper bound on their Mahler measure; see Lemma~\ref{lem:Mahler coefficient bound}. 
Applying Schinzel's result \cite[Addendum, Corollary 1']{Schinzel} to the polynomial $z-\alpha$ for any algebraic integer $\alpha$ that is not a root of unity and has Galois group $\Z/n\Z$ yields
\begin{equation}
\label{eq:Schinzel}
m_n \geq \frac{1}{2} \log \Big(\frac{1+\sqrt{5}}{2}\Big)n;
\end{equation}
see also \cite{AmorosoDvornicich}, where a lower bound is derived in the more general setting of algebraic numbers, not just integers. We show the following.

\begin{introtheorem} \label{mainMM} 
We have $m_q = \m(\alpha_q)$ for $q=3,5,7$. 
\end{introtheorem}

Note that $\kappa(3)=\kappa(5)=2$ whereas $\kappa(7)=4$. Theorem \ref{mainMM} is proven by computer search. A naive search for $n=7$ would involve $\approx 8 \cdot 10^{18}$ cases, but we use some algebraic number theory and class field theory to cut down the size of the search space by restricting the factors and absolute value of conductors. 
We also compute $m_n$ for other small values of $n$. Certain even values of $n$ behave differently; for example, $m_4 = \m(\zeta_5 + \zeta_5^3)$.  But we dare to ask whether it is always true that 
\begin{equation} 
\label{con} 
m_n \stackrel{?}{=}\m(\alpha_n) \ \ (n \mbox{ odd}).
\end{equation} 

\subsection*{Towards a conjecture on the growth rate of $m_n$ and $\m(\alpha_n)$; bounds for $\kappa(n)$} 
In Lemma~\ref{lem:m(an) unconditional bounds} we show that $n \ll m(\alpha_n) \ll n\log{n}$ for all $n$. We pose the following more precise conjecture for odd $n$.
\begin{introconjecture} 
\label{mainconj} 
We have $\m(\alpha_n) \asymp n \log{\log{n}}$ for all but $o(X)$ odd values of $n \leq X$ as $X \rightarrow + \infty$.
\end{introconjecture} 
{Note that in this conjecture, the values of the odd integers $n$ are no longer restricted by $n \in \NN_k$ for fixed $k$, as was the case in Theorem \ref{mainas}. }
Under \eqref{con}, Conjecture~\ref{mainconj} would yield $m_n \asymp n \log{\log{n}}$ for all but $o(X)$ of the odd integers $n \leq X$, as $X \to +\infty$. This indicates that the linear lower bound \eqref{eq:Schinzel} would, for typical odd $n$, be tight up to a $\log{\log{n}}$ factor. On the other hand, Theorem~\ref{mainas} yields
\begin{equation*}
\liminf_{n \to \infty} \frac{m_n}{n} < \infty,
\end{equation*}
thus \eqref{eq:Schinzel} cannot be improved upon much in full generality.

Conjecture~\ref{mainconj} fits into a general framework of conjectures (mostly not of a quantitative nature) on height growth in Galois extensions and for certain families of Laurent polynomials, see Remarks \ref{RemGalGrow} and \ref{RemLauGrow}. 

We show that Conjecture~\ref{mainconj} is implied by a particular form of the prime tuplets conjecture (Conjecture \ref{con:PT_with_error_term}) in combination with the two further conjectures (\ref{con:mCk logk} and \ref{con:malphan mCk}) that $\m(C_{\kappa(n)})  \overset{?}{\gg}  \log \kappa(n)$ and $\m(\alpha_n)/n  \overset{?}{\asymp}  \m(C_{\kappa(n)})$, both for all but $o(X)$ odd values of $n \leq X$, as $X$ tends to infinity.  These two conjectures are cousins of Theorems \ref{mainas} and \ref{mainRW}, but now with $k=\kappa(n)$ a function of $n$, rather than being fixed. To show that the three conjectures imply our main conjecture we need to involve almost-always bounds for $\kappa(n)$, entering the circle of ideas around Linnik's constant. Using the Siegel--Walfisz theorem, we prove an unconditional lower bound, stating that for all but $O(X/\log{\log{X}})$ values of $n \leq X$ we have $\kappa(n) > \log{n}/\log{\log{n}}$ (Proposition \ref{lem:k(n) lower bound}); and we invoke an almost-always upper bound $\log \kappa(n) \ll \log \log n$ proven by Granville assuming the prime tuplets conjecture.  

To summarize the situation, the overarching conjecture about the true growth behaviour of the minimal height for algebraic integers with cyclic Galois group of fixed odd order has now been decoupled into open problems on asymptotics of random walks and quantitative equidistribution results. 

\medskip

In Section \ref{open}, we have collected open problems and generalisations. See Section \ref{sec:code} for a description of \texttt{SageMath}  \cite{SageMath} routines used in this paper and posted online. 
 
\section{On the structure of the set $\mathcal N_k$} 
  
\subsection{Density of the set $\NN_k$} \label{asNNk}
 
Theorem~\ref{mainas} only makes sense if $\NN_k$ is of infinite size. A more precise result is the following proposition. 

\begin{notation} For $X>0$ and a set $A$ of positive integers, let $A(X) \coloneqq  A \cap [1,X]$ denote the set of elements of $A$ smaller than or equal to $X$. 
 \end{notation} 
 
\begin{proposition} 
\label{nnkinf} 
Fix $k \geq 1$. For sufficiently large $X$, we have
$$
\# \mathcal{N}_k(X) = \frac{kX}{\phi(k) \log{X}} + O_k \Big( \frac{X}{(\log{X})^2} \Big).
$$
\end{proposition} 
\begin{proof}
Define $\mathcal{N}'_k(X) = \{n \in \Z_{>0} : n \leq X, \, kn + 1 \text{ prime}\}$ and, for each $j=1,2,\ldots, k-1$, the set
\begin{equation*}
E_j(X) \coloneqq  \{n \in \Z_{>0} : n \leq X, \, kn + 1 \text{ prime}, \, jn + 1 \text{ prime}\}.
\end{equation*}
Then
$$\mathcal{N}_k(X) = \mathcal{N}'_k(X) \setminus \bigcup_{j=1}^{k-1} E_j(X).$$
Taking cardinalities on both sides and using the union bound, we obtain the inequalities
\begin{equation}
\label{eq:inequalities_NkX}
\# \mathcal{N}'_k(X) - \sum_{j=1}^{k-1} \# E_j(X) \leq \# \mathcal{N}_k(X) \leq \mathcal{N}_k'(X).
\end{equation}
The set $\mathcal{N}_k'(X)$ simply lists all prime numbers $\leq kX+1$ congruent to $1 \bmod{k}$; hence, by the Siegel--Walfisz theorem \cite[Corollary 5.29]{IK},
\begin{equation}
\label{eq:primes_in_AP}
\# \mathcal{N}_k'(X) = \frac{kX}{\phi(k) \log{kX}} + O \Big( \frac{kX}{(\log{kX})^2} \Big)
\end{equation}
where the implied constant is independent of $k$ (in particular, \eqref{eq:primes_in_AP} holds even when $k$ depends on $X$, although we do not need that here). Furthermore, by \cite[Lemma~5.1]{Ford-et-al}, each $E_j(X)$ is of cardinality $ \# E_j(X) \ll_k X/(\log{X})^2$. Combining this with \eqref{eq:primes_in_AP} and \eqref{eq:inequalities_NkX} proves the proposition.
\end{proof}
  
\subsection{Explicit construction of elements in $\NN_k$} 
 
The following is an easy construction of infinitely many elements in $\NN_k$ for fixed $k$. For this, choose $k-1$ prime numbers $p_1,\dots,p_{k-1}$ such that 
$$ k < p_1 < p_2 < \dots < p_{k-1}.$$ Consider the system of congruences  
\begin{equation} \label{CRT1} in+1 = 0 \mbox{ mod } p_i, \mbox{ for all } i=1,\dots,k-1, \end{equation} 
and notice that $i<k<p_i$, so that $i$ is invertible modulo $p_i$, so \eqref{CRT1} is in fact $n = -1/i \mbox{ mod } p_i$. 
By the Chinese Remainder Theorem, solutions to \eqref{CRT1} are the same as solutions to a single congruence 
\begin{equation} \label{CRT2} n= r \mbox{ mod } p_1\cdots p_{k-1}\end{equation} 
for a suitable $r$ depending on $k$. 
Now consider the arithmetic progression 
\begin{equation} \label{AP} kr+1 + \lambda kp_1\cdots p_{k-1} \mbox{ with } \lambda \in \Z. \end{equation} 
By Dirichlet's theorem, this sequence contains infinitely many primes, as soon as $kr+1$ and $kp_1\cdots p_{k-1}$ are coprime, i.e., as soon as $p_i$ does not divide $kr+1$ for all $i$. Now $$kr+1 \mbox{ mod } p_i = kn+1 \mbox{ mod } p_i = -k/i+1 \mbox{ mod } p_i \neq 0 \mbox{ mod } p_i$$ since $i \neq k \mbox{ mod } p_i$. 
Any of the infinitely many primes $p$ in \eqref{AP} corresponds with an integer $n$ for which the equations \eqref{CRT2}, and hence \eqref{CRT1} are satisfied. Choosing $p$ so large that  for the corresponding $n$, we have $in+1 \neq p_i$ for all $i=1,\dots,k-1$, then also $in+1$ is composite for all such $i$. This means that, for all but finitely many $p$, we have $n \in \NN_k$.

\section{Asymptotics of $\m(\alpha_n)$ for fixed $k$; proofs of Theorem~\ref{mainas}} \label{proof1} 
\label{sec:asymp_m_alpha_n_fixed_k}

We first give a short proof of Theorem~\ref{mainas} using a new approach of Kowalski and Untrau via Wasserstein distances (here, the `multidimensionality' is encoded in the correct choice of measures, that are integrated against a suitable one-variable function). After that, we present a more standard approach using convergence of Riemann sums to Riemann integrals via Weyl equidistribution (here, the `multidimensionality' is encoded by integrating a multivariable function), and the quantitative version using the Hlawka--Koksma inequality. 

\subsection{The case $k=1$} Concerning the `degenerate' case $k=1$, we have the following. 

 \begin{lemma} \label{pos}
$\m(C_k)=0 \iff k=1$. 
\end{lemma} 
\begin{proof} We have $\m(C_1)=0$.  On the other hand, for $k>1$, by Jensen's formula  \cite[Lemma 1.9]{EverestWard}, $\m(C_k) = \int_{\T^d} \log^+|F_k| $, which is strictly positive by noticing that, with $\mathbf 1=(1,\dots,1)$ the all-one vector, $F_k(\mathbf{1}) = k>1$ for $k \geq 2$, so that $\log^+|F_k|$ is positive on an open neighbourhood of a point in $\T^d$, hence $\m(C_k)>0$. 
\end{proof} 

If $n \in \NN_1$, then also $\m(\alpha_n)=0$, so \eqref{was} holds. (Formally, \eqref{nowas} makes no sense for $k=1$.) 
From now on, we assume $k>1$. 

\subsection{Using Wasserstein distance} 

\begin{proof}[Proof of Theorem \ref{mainas}] To prove Theorem \ref{mainas} for $k>1$, it suffices to prove \eqref{was}. For this, we will use the \emph{Wasserstein distance} $W(\nu_1,\nu_2) = \inf \mathbf{E}(|X-Y|)$ between two probability measures $\nu_1, \nu_2$ on $\Cc$, where $X$ and $Y$ run over all random variables on $\Cc$ such that $X$ has law $\nu_1$ and $Y$ has law $\nu_2$. (To be more precise, this is the $1$-Wasserstein distance on the complex plane equipped with the standard metric.) Denote by $\mathrm{Lip_1}(\Cc,\R)$ the set of $1$-Lipschitz functions $\Cc \to \R$, i.e.,  those functions satisfying $|u(z)-u(w)| \leq |z-w|$ for all $z,w \in \Cc$.
The Wasserstein distance satisfies Kantorovich--Rubinstein duality 
\begin{equation} 
\label{w1} 
W(\nu_1,\nu_2) = \sup_{u \in \mathrm{Lip_1}(\Cc,\R)} \Big| \int_\Cc u \, \d \nu_1 - \int_\Cc u \, \d \nu_2 \Big|,
\end{equation}
see \cite[Theorem~1.2(6)]{KU}.
Let $\delta_x$ denote the Dirac measure at $x$. Denote by $F_{k\ast} \mu_{\T^d}$ the pushforward of the Haar measure $\mu_{\T^d}$ on the torus $\T^d$ by the Laurent polynomial $F_k$, seen as a map $F_k \colon \T^d \rightarrow \Cc$. Define the probability measures
\begin{equation*}
\nu_1 = \frac{1}{p} \sum_{a \in \FF_p} \delta_{ \tr_H(\zeta_p^a)} \ , \quad \nu_2 = F_{k\ast} \mu_{\T^d}
\end{equation*}
and the function $h(z)=\log^+|z|$. It suffices to prove
\begin{equation}
\label{eq:Wasserstein}
\Big|\frac{k}{p} \m(\alpha_n) + \frac{\log{k}}{p} - \m(C_k)\Big| = \Big| \int_\Cc h \, \d \nu_1 - \int_\Cc h \, \d \nu_2 \Big| \leq W(\nu_1, \nu_2) \leq \frac{4\sqrt{3} k(k+1)}{p^{1/\phi(k)}}.
\end{equation}
Indeed, assuming \eqref{eq:Wasserstein}, the triangle inequality yields
\begin{equation} 
\label{changebnd} 
 \Big| \frac{\m(\alpha_n)}{n} - \m(C_k) \Big|  \leq \frac{4\sqrt{3} k(k+1)}{p^{1/\phi(k)}} +\Big(\frac1n-\frac{k}{p}\Big)\, \m(\alpha_n) + \frac{\log{k}}{p} 
\end{equation} 
From \eqref{mahlerheight}, $\m(\alpha_n) \leq n  \max \log |\alpha_n'|$, where the maximum extends over all Galois conjugates $\alpha_n'$ of $\alpha_n$. Since each conjugate $\alpha_n'$ is a sum of $k$ roots of unity, we obtain that $\m(\alpha_n) \leq n \log{k}$. Combining this with \eqref{changebnd}, noticing $1/n-k/p=1/(np)$, yields the theorem.

It remains to show \eqref{eq:Wasserstein}. For the first equality, observe that 
\begin{equation} 
\label{w2}  
\int_\Cc h \, \d \nu_1 = \frac{1}{p} \Big( \log{k} + \sum_{a \in \FF_p^{\times}} \log^+ |\tr_H(\zeta_p^a)| \Big) = \frac{k}{p} \m(\alpha_n) + \frac{\log{k}}{p};
\end{equation} 
Furthermore, 
\begin{align} 
\label{w3} 
\int_\Cc h \, \d \nu_2 &= \int_\Cc \log^+|t| \, \d(F_{k\ast} \mu_{\T^d}) = \int_{\T^d} \log^+|F_k(\mathbf x)|\, \d \mu_{\T^d} \nonumber \\ 
&\overset{(\dagger)}{=} \int_{\T^{d+1}} \log | x_0 + F_k(\mathbf x)|\, \d \mu_{\T^{d+1}} = \m(C_k), 
\end{align} 
where $(\dagger)$ is an application of Jensen's formula.

For the inequality in the middle of \eqref{eq:Wasserstein}, we apply \eqref{w1} to the function $u=h$; hence it suffices to show that $h \in \mathrm{Lip_1}(\Cc,\R)$. This is trivial if both arguments $z,w$ are inside the complex unit circle; if one of them is inside (say, $|w| \leq 1$) and one is outside (say, $|z| \geq 1$), then the claim is $\log|z| \leq  |z-w|$, which follows from $\log x \leq x-1$ for all $x \in \R_{\geq 1}$ and the reverse triangle inequality. Finally, if both $z$ and $w$ are outside the complex unit circle, we can use that the real function $\log x$ is everywhere differentiable on $\R_{\geq 1}$ with derivative $1/x$ upper bounded by $1$, so is $1$-Lipschitz in that region; hence $|\log |z| - \log |w|| \leq ||z|-|w|| \leq |z-w|$ by the reverse triangle inequality. 

For the last inequality in \eqref{eq:Wasserstein}, {we apply \cite[Corollary~3.6 and Remark~3.7(2)]{KU} with the set  $\mathrm{Z}$ in that reference consisting of the complex zeros of our polynomial $g(X) = X^k-1$.} The prime denoted $q$ in \cite[Remark~3.7(2)]{KU} is our $p$, which is $1 \bmod{k}$ and therefore totally split in the splitting field of $g$ over $\Q$. In our notation, their measure $\nu_q$ is then exactly our $\nu_1$. Furthermore, to see that their $\mu_g$ can be identified with our $\nu_2$ one may consult \cite[Section~3]{KU2}. To conclude, the proofs of \cite[Prop.~3.4 and Cor.~3.6]{KU} imply the last inequality in \eqref{eq:Wasserstein}, where we also used that the constant denoted $\mathrm{C}_{\mathrm{Z}}$ in the proof of \cite[Prop.~3.4]{KU} is $1$ by \cite[Lem.~3.3]{KU}.
\end{proof} 

\begin{remark}
Let $\mathcal{P}_k$ be the set of primes $p = kn+1$ with $n$ ranging over $\mathcal{N}_k$. Then $\mathcal{P}_k$ is a subset of the set $\mathcal{A}_k$ of \emph{$k$-admissible integers} in \cite{Untrau}. Taking $d=k$ and $F = \log^+|\cdot|$ in \cite[Theorem~1.3]{Untrau} then also yields the qualitative convergence \eqref{nowas}.
\end{remark}

\subsection{Using Hlawka--Koksma quantitative equidistribution} 
Another approach to bounding the left hand side in \eqref{was} is by using the Hlawka--Koksma inequality. For this, we first interpret that expression as the `error' in approximating a $d$-dimensional Riemann integral by a Riemann sum. We now explain which function and which sample points this applies to. 

\subsubsection*{Mahler measure and `integration error'} Fix a choice of a primitive root $\lambda$ modulo $p$; then $\nu\coloneqq  \lambda^{n}$ is an element of exact multiplicative order $k$ modulo $p$, so a generator for the unique cyclic subgroup of $(\Z/p)^*$ of order $k$. 
We can rewrite the Gaussian period $\alpha_n$ as  
\begin{equation} 
\label{alphan} 
\alpha_n \coloneqq  \sum_{r=0}^{k- 1} e^{\frac{2 \pi i}{p} \nu^r}. 
\end{equation} 

\begin{definition} 
\label{def-f} 
Recall that $d=\phi(k)$. With $F_k(\x)$ as in  Definition~\ref{def:A}, define the function
\begin{equation*} 
\label{deff} 
f(\vtheta)\coloneqq  f_k(\theta_1,\dots,\theta_{d})\coloneqq  \log^+{|F_k(e^{2 \pi i \theta_1},\dots,e^{2 \pi i \theta_d})|}. \end{equation*}
\end{definition} 

\begin{definition} 
\label{def-sp} 
Define the set $S_p = \{ \vtheta_1, \vtheta_2, \ldots, \vtheta_{p-1}\}$ of sample points in $[0,1)^d$ by 
\begin{equation*}
\vtheta_j\coloneqq (\{\lambda^j/p\},\{\nu \lambda^j/p\},\dots,\{\nu^{d-1} \lambda^j/p\}) \in [0,1)^{d} \quad \text{for } j = 1, 2, \ldots, p-1.
\end{equation*}
\end{definition}

\begin{lemma}
For fixed $k \in \Z_{\geq 2}$, the `integration error'
\begin{equation*} 
\err_p \coloneqq \frac{1}{p-1} \sum_{j=1}^{p-1} f(\vtheta_j) - \int_{[0,1)^d} f(\vtheta) \, \d\vtheta
\end{equation*}
satisfies $\err_p = {\m(\alpha_n)}/{n} - \m(C_k)$.
\end{lemma}
\begin{proof}
Because $\lambda^j$ for $j=1,\dots,n$ is a set of representatives for the cosets of $\langle \nu=\lambda^n \rangle \cong \Z/k$ in $(\Z/p)^* = \langle \lambda \rangle \cong \Z/(kn)$, the conjugates of $\alpha_n$ are precisely given by   
$$\alpha_n^{(j)} \coloneqq  \sum_{r=0}^{k- 1}  e^{\frac{2 \pi i}{p} \lambda^j \nu^r}$$ where $j=1,\dots, n$. The vectors $\mathbf{a}^{(r)}=(a^{(r)}_0,\dots,a^{(r)}_{d-1})$ for $r=0,\dots,k-1$ in $A_k$ are defined so that 
$$ \nu^r = \sum_{s=0}^{d-1} a_s^{(r)} \nu^s \mbox{ mod } p, $$
so the Mahler measure of $\alpha_n$ is 
$$ \m(\alpha_n) =  \sum_{j=1}^{n} \log^+{\left| \sum_{r=0}^{k- 1}  \prod_{s=0}^{d-1} \left(e^{\frac{2 \pi i}{p} \lambda^j}\right)^{a^{(r)}_s \nu^s} \right|} = \sum_{j=1}^{n} f(\vtheta_j). $$
The left hand side of this equation does not change when replacing $\alpha_n$ by a conjugate element --- in particular, when changing $\zeta_p$ to $\zeta_p^{\lambda^t}$ for some $t$. Hence the right hand side also does not change when, instead of summing over $j=1,\dots,n$, we sum over $j=bn+1,\dots,(b+1)n$ for some positive integer $b$. 
Thus, we have an equality $$\m(\alpha_n)/n = \frac{1}{p-1} \sum_{j=1}^{p-1} f(\vtheta_j).$$

The Mahler measure of $C_k$ is (using Jensen's formula) 
\begin{equation*} \m(C_k) = \int_{\T^{d+1}} \log |x_0+F_k(\mathbf{x})| \, \d\mu = \int_{\T^d} \log^+{|F_k(\x)|} \, \d\mu =  \int_{[0,1)^d} f(\vtheta) \, \d\vtheta, \end{equation*}
which implies the result. 
\end{proof}

\begin{remark}[Weyl equidistribution] \label{weyl} Notice that the function $f$ is bounded and continuous on $[0,1)^d$. The qualitative statement \eqref{nowas} (for $k>1$) is equivalent to
\begin{equation} 
\label{eq:lim} 
\lim_{\substack{p \rightarrow + \infty \\ \frac{p-1}{k} \in \NN_{k}}} \left|\err_p \right| = 0,
\end{equation} 
which, by Weyl's equidistribution theorems (see, e.g.\ \cite[Ch.\ 1, \S 6]{KuipersNiederreiter}), is equivalent to 
\begin{equation} 
\label{limweyl} 
\lim_{p \rightarrow +\infty} \frac{1}{p} \sum_{j=1}^{p-1} e^{2 \pi i\, \mathbf{h} \cdot \vtheta_j}  = 0 
\end{equation} 
for all $\mathbf{h} \in \Z^d \setminus \{\mathbf{0}\}$. This was proven by Myerson in \cite[Theorem 12]{Myerson}; for the sake of completeness, we recall the argument. 
Write $P(t)=h_0 + h_1 t + \dots + h_{d-1} t^{d-1} \in \Z[t]$. Then
\begin{equation*}
 \sum_{j=1}^{p-1}  e^{2\pi i \h \cdot \vtheta_j} = \sum_{j=1}^{p-1} e^{\frac{2 \pi i}{p} \lambda^j P(\nu)} = \sum_{s=1}^{p-1}  e^{\frac{2\pi i}{p} s P(\nu)} = 
\begin{cases}
p-1 & \text{if } P(\nu) \equiv 0 \bmod{p} \\ 
-1 & \text{else,}
\end{cases}
\end{equation*}
where the penultimate step follows 
since $\lambda$ is a primitive root modulo $p$, and the last equality holds because the last sum is a complete exponential sum barring the term corresponding to $s=0$.
{Write $\Phi_k$ for the $k$-th cyclotomic polynomial.} With $d=\phi(k)$ and $\deg P<d$, the polynomials $P$ and $\Phi_k$ are coprime, so there exists $n \in \Z$ and $a(t), b(t) \in \Z[t]$ such that  $a(t)P(t)+b(t)\Phi_k(t) = n$ in $\Z[t]$. Since $\Phi_k(\nu)=0 \mbox{ mod } p$, we have that $P(\nu)=0 \mbox{ mod } p$ only if $p$ divides $n$. Hence there are only finitely many $p$ with $P(\nu)=0 \mbox{ mod } p$, and this proves that the sum equals $-1$ for any sufficiently large $p$, so the limit in Equation \eqref{limweyl} is indeed zero. 
 \end{remark} 

\subsubsection*{Discrepancy and the Hlawka--Koksma bound}  The argument in Remark \ref{weyl} says precisely that $S_p$, as in Definition \ref{def-sp}, is equidistributed. 
A method to study the convergence rate of the integration error $\err_p$ is by using the Hlawka--Koksma inequality, which says that 
\begin{equation} 
\label{eq:Hlawka-Koksma}
\abs{\err_p} \leq V(f) \cdot D(S_p),
\end{equation}
where $V(f)$ is the variation, in the sense of Hardy and Krause (an epithet that we will leave out from now on), of the function $f$, see, e.g., \cite[Thm. 3.11]{Leobacher-et-al}, and $D(S_p)$ is the discrepancy of the set $S_p$, in the sense of the following definition. 
\begin{definition} We define $\mathcal B_d$ to be the set of axis-parallel boxes in $\R^d$ with $\mathbf{0}$ as vertex and contained in the unit cube, so $\mathcal B_d$ consists of sets of the form 
 $ B= \prod_{i=1}^d [0,b_i) \subset \R^d$ for $0<b_i<1$. These have Lebesgue measure $\lambda(B) = \prod_{i=1}^d b_i$. 
 A finite subset $S$ of $[0,1)^d$ has \emph{discrepancy}
\begin{equation*}
 D(S)\coloneqq  \sup_{B \in \mathcal B_d} \left| \frac{\# (S \cap B)}{\#S} - \lambda(B)    \right|;
 \end{equation*}
this type of discrepancy is sometimes called \emph{star discrepancy} in order to distinguish it from the \emph{extreme discrepancy} (which is sometimes also called simply discrepancy) \cite[Chapter~2.2]{Leobacher-et-al}. A sequence $\{S_n\}_{n=1}^\infty$ of finite subsets of $[0,1)^d$ is \emph{equidistributed} if $\lim\limits_{n \rightarrow + \infty} D(S_n) = 0$. 
\end{definition} 

\subsubsection*{Finite variation by smooth approximation} We first prove that $V(f)$ is finite, which we do by approximating the (non-smooth) function $f$ by smooth ones. We will not define the variation $V(f)$ in general, but use the following result. 
\begin{lemma} \label{ulim} 
Uniform limits of functions on $[0,1)^d$ in the sup-norm with continuous mixed partial derivatives up to order one in each variable are of bounded variation. 
\end{lemma} 
\begin{proof} 
As in \cite{Aistleitner-et-al}, consider the set $\mathcal V_\infty$ of measurable functions on $[0,1)^d$ that are limits, in the sup-norm, of finite linear combinations of characteristic functions of sets in $\mathcal B_d$, and denote by $\mathcal V$ the set of functions in $\mathcal V_\infty$ of bounded variation.

If $\phi \colon [0,1)^d \rightarrow \R$ is a function for which all mixed partial derivatives up to order one in all variables are continuous, then its variation is given by 
\begin{equation} \label{zaremba} 
V(\phi) = \sum_{\emptyset \neq u \subseteq \{1,\dots,d\}} \int_{[0,1)^{|u|}} \left| \partial_{\mathbf \vtheta_u} \phi(\vtheta_u,\mathbf 1) \right| \d\vtheta_u, 
\end{equation} 
where $(\vtheta_{u},\mathbf 1)$ is the vector where all entries not in $u$ are replaced by $1$ (cf.\ \cite[p.\ 69]{Leobacher-et-al}).

Assume that we have a sequence $\{f_i\}$ of functions converging to some function $f$, as in the lemma. Formula \eqref{zaremba} applied to $\phi=f_i$ implies that $f_i$ are of bounded variation, and hence in particular, are in $\mathcal V_\infty$ \cite[Theorem 4.2]{Aistleitner-et-al}. Since $f$ is a uniform limit of $f_i$, by the triangle inequality it can also be approximated by finite linear combinations of characteristic functions of sets in $\mathcal{B}_d$, so $f \in \mathcal V_\infty$, too. 
We can hence use the following estimate \cite[Prop.\ 2.3(ii)]{Aistleitner-et-al}: 
\begin{equation}  \mbox{ if\ }f_i,f \in \mathcal V_\infty \mbox{ with\ } ||f-f_i||_\infty \rightarrow 0,  \mbox{ then\ }V(f) \leq \liminf V(f_i), \end{equation} 
which implies in particular that $f$ is also of bounded variation. 
\end{proof} 

\begin{lemma} 
\label{lem:f is of bounded variation}
The function $f$ in Definition \ref{deff} is of bounded variation. 
\end{lemma} 

\begin{proof} 
By Lemma \ref{ulim}, it suffices to approximate $f$ uniformly by suitably smooth functions. For this, we use the shorthand $F\coloneqq F_k$ and define, for $i > 0$, 
$$ f_i(\vtheta) = 
\left\{ \begin{array}{ll} \log |F(\vtheta)| & \mbox{ if } |F(\vtheta)|> 1+1/i; \\
g_i(\vtheta)  & \mbox{ if } |F(\vtheta)| \leq 1+1/i; 
\end{array} \right. 
$$
for suitable functions $g_i$ to be specified.

We first note that on $|F|>1+1/i$,  $f_i$ is smooth, and for an non-empty index $\alpha \subseteq \{1,\dots,d\}$ we have $|\partial_\alpha f_i| = |P|/|F|^{|\alpha|} \leq |P| (i/(i+1))^{|\alpha|}$ for some polynomial of degree $\leq d$ in $\{\partial_\beta F\}$ for $|\beta|\leq |\alpha|$.  Now since $F(\vtheta) = \sum_{\mathbf a \in A_k} e^{2 \pi i \mathbf a \cdot \vtheta}, $ the partial derivatives satisfy $$|\partial_\beta F| \leq k \max\{ (2 \pi ||\mathbf a||)^{|\beta|} \colon \mathbf a \in A_k \}$$ (with $||(a_i)|| \coloneqq \sum |a_i|$), which is uniformly bounded in terms of $k$ only. 
By choosing smooth functions $g_i$ that are zero on $|F(\vtheta)|< 1$ and all of whose first partial derivatives equal those of $\log |F(\vtheta)|$ on  $|F(\vtheta)| = 1+1/i$, we can apply Lemma \ref{ulim} to $f_i$. 
\end{proof} 

The result of Myerson in Remark \ref{weyl} says precisely that $D(S_p) \rightarrow 0$ as $p \rightarrow +\infty$ with $(p-1)/k \in \mathcal N_k$. Hence \eqref{eq:Hlawka-Koksma} also implies \eqref{nowas}, and splits questions of quantitative control of the integration error into those about variation and discrepancy. 

\subsubsection*{Explicit upper bounds?} Numerics suggest a decay rate of $\err_p$ on the order of $1/p$ for $k\in \{2,3,4\}$ and a more irregular decay rate $\ll 1/p$ for higher $k$.
What is the quantitative nature of the bound arising from the Hlawka--Koksma inequality? 
\begin{enumerate} 
\item Concerning the variation, it follows that $V(f)$ is a finite constant for fixed $k$. For varying $k$, we see from the proof of 
\ref{lem:f is of bounded variation} that 
\begin{equation} \label{upperVbad} V(f) \leq k^d (2 \pi \max_{a \in A_k} ||a||)^{d^2}, \end{equation} 
an upper bound that grows faster than exponential in $k$ (with $d \approx \phi(k) \approx k$). 
\item Concerning the discrepancy, Roth's general lower bound on the discrepancy \cite[Theorem~2.24]{Leobacher-et-al} shows that $D(S_p)\gg_k \log(p)^{(d-1)/2}/p$, and Schmidt showed that $D(S_p) \gg \log(p)/p$ for $k=4$ (i.e., $d=2$; this is essentially best possible for that dimension \cite[Theorem~2.25]{Leobacher-et-al}). Based on numerical experiments for $k=4$, we speculate that this lower bound is sharp for almost all $n \in \mathcal{N}_4$. Nevertheless, if $\mathcal{N}_4$ contains infinitely many (almost-)squares, then the corresponding subsequence of values seems to have an asymptotics at least on the order of $1/\sqrt{p}$. 
\item The Hlawka--Koksma inequality arises in general from a reproducing Hilbert kernel method in which a lot of sums with alternating signs are upper bounded by the sum of their absolute values;  in reality, there might be more cancellation in the `integration error', see for example Zaremba's formula \cite[Thm. 3.10]{Leobacher-et-al}.
\end{enumerate}

\begin{remark} 
In \cite{DimitrovHabegger}, Dimitrov and Habegger derive discrepancy estimates for Galois orbits of torsion points on $\Gm^d$ (they use yet another integration error bound, the `modulus of continuity'). We are instead dealing with cosets of a subgroup of the Galois group. 
\end{remark}

\subsection{Lower order terms} With some extra work, it is possible to determine subdominant terms in an expansion of $\m(\alpha_n)$; for example, we have the following for $k=2$. 

\begin{proposition} 
\label{LOT} 
As $n \rightarrow + \infty$ with $n \in \NN_2$, 
$$ \m(\alpha_n) \sim \m(C_2) \cdot n + \frac12(\m(C_2) - \log(2)) + O\left(\frac{1}{n}\right).$$
\end{proposition} 
  
  \begin{proof} This is an application of Euler--Maclaurin summation. More precisely, $$\m(\alpha_n) = \sum_{j=1}^n \log \max \{ 1, 2 \cos(2 \pi j/p) \},$$ and in this summation, the maximum is larger than $1$ precisely if $1 \leq j \leq p/6$  or $p/3 \leq j \leq n$. If we set $f(t)\coloneqq \log (2 \cos(2 \pi t/p))$, then $\m(\alpha_n) = \Sigma_1 + \Sigma_2$ where 
  $$ \Sigma_1
 = \sum_{i=1}^{\lfloor p/6 \rfloor} f(i) \mbox{ and }   \Sigma_2
 = \sum_{i=\lceil p/3 \rceil}^n f(i). $$
 Euler--Maclaurin summation gives 
  \begin{equation} \label{EM1}  \Sigma_1 = \int_0^{\lfloor p/6 \rfloor} f(t) \, \d t  + \frac{1}{2} (f( \lfloor p/6 \rfloor) - f(0)) + \frac{B_2}{2}  (f'( \lfloor p/6 \rfloor) - f'(0)) + R_2 \end{equation} 
 with $|R_2| \leq  \int_0^{\lfloor p/6 \rfloor} |f''(t)| \, \d t$.
 \begin{enumerate} 
 \item Since $f''(t) = -4 \pi^2/p^2 \cdot \frac{1}{\cos(2 \pi t/p)^2}$ and  $\cos(2 \pi t/p)>1/2$ in this region, we find $$|R_2| \leq 16 \lfloor p/6 \rfloor \pi^2/p^2 = O(1/p).$$
\item If we substitute $u=t/p$ in the first term of \eqref{EM1} and use $\lfloor p/6 \rfloor = p/6-\{p/6\}$, we find that it equals 
 \begin{equation} \label{EM2} p \int_0^{1/6} \log(2 \cos 2 \pi u) \, \d u - p \int_{1/6-\{p/6\}/p}^{1/6}  \log(2 \cos 2 \pi u) \, \d u; \end{equation} 
 Note that 
 \begin{align*} 
 \m(C_2) & = \int_0^1 \log^+ |2 \cos 2\pi u| \, \d u = 2( \int_0^{1/6} + \int_{1/3}^{1/2}) \log (2 \cos 2\pi u) \, \d u \\ & = 4\int_0^{1/6} \log (2 \cos 2\pi u) \, \d u,
 \end{align*} 
 so the first term in \eqref{EM2} equals $p/4 \cdot \m(C_2)$, whereas the second term in \eqref{EM2}, using the expansion $ \log(2 \cos 2 \pi u) = - 2 \sqrt{3} \pi (u-1/6) +O((u-1/6)^2)$, is equal to $-\sqrt{3} \pi / p + O(1/p^2) = O(1/p)$. 
\item The second term in \eqref{EM1} is \begin{align*} f(\lfloor p/6 \rfloor) -f(0) & = f(p/6-\{p/6\}) - \log 2= \log ( 2 \cos ( 2 \pi/6 - \{p/6\}/p)) - \log 2 \\ &= - \log 2 + O(1/p), \end{align*} again using the Taylor series. 
\item For the third term in \eqref{EM1}, we use $f'(t) = - 2 \pi/p \tan (2 \pi t/p)$, hence $f'(0)=0$ and $f'(\lfloor p/6 \rfloor)= -2 \pi/p \tan(\pi/3-2 \pi \{p/6\}/p) = O(1/p)$.
\end{enumerate} 
In the end, we find that 
$$\Sigma_1 = p/4 \cdot \m(C_2) - (\log 2)/2 + O(1/p).$$ 
In a similar way, 
$$ \Sigma_2 = p/4 \cdot \m(C_2) + O(1/p),$$
so that in the end, using $p=2n+1$, 
\begin{align*} 
\m(\alpha_n) &=  \Sigma_1 + \Sigma_2 = p/2 \cdot \m(C_2) - (\log 2)/2 + O(1/p) \\ & = (n+1/2) \m(C_2) -( \log 2)/2 + O(1/n)\\ &= \m(C_2) \cdot n + (\m(C_2)-\log 2)/2 + O(1/n), 
\end{align*}
as was to be proven. 
\end{proof} 
  
\begin{remark}   
The asymptotics of the Mahler measure of $2 \cos(2\pi/N)$ as $N \rightarrow + \infty$ is also studied in \cite[\S 6]{SmythTR}, where it is shown that only finitely many of these attain the minimal Mahler measure amongst \emph{totally real} algebraic integers of given degree. For odd $n$, $\alpha_n$ is totally real, but our Galois-theoretic condition is much more restrictive than that. 
\end{remark}

\section{The geometry of the cyclopolytope and the cyclovariety;  proof of Theorem \ref{mainCY}} 
\label{sec:geometry of the cyclovariety} \label{NP}

Considering $x_0$ as a parameter, the curves $C_4$ and $C_6$ (in the remaining two variables $x_1$ and $x_2$) can be compactified into elliptic curves. In this section, we study how this result generalized to other $C_k$.

\subsection{The cyclopolytopes and the cyclovarieties} \label{cvarcpol}

\subsubsection*{Terminology} We recall some terminology concerning polytopes (the bounded convex hull of a finite set of ``spanning'' vertices), referring to \cite{polysurvey} for details. A subset $S$ of $\R^d$ is called \emph{centrally symmetric} if whenever $v \in S$, then also $-v \in S$. The \emph{(polar) dual} of a subset $S \subseteq \R^d$ is defined by $S^*\coloneqq  \{ w \in \R^d \colon \langle v,w \rangle \leq 1, \forall v \in S\}$. Note that the dual is sometimes taken to be the (isomorphic) point reflection of $S^*$ in the origin, $-S^* =  \{ w \in \R^d \colon \langle v,w \rangle \geq -1, \forall v \in S\}$; for centrally symmetric $S$, we have $-S^*=S^*$.

We consider only polytopes $P \subseteq \R^d$ that are full-dimensional and contain the origin. A polytope $P$ is called \emph{integral} if the spanning vertices are in $\Z^d$. A polytope $P$ is called \emph{reflexive} (a concept due to Batyrev \cite{Batyrev}) if the dual polytope is also integral.  A polytope $P$ is reflexive precisely if $P$ can be described as $P = \{ v \in \R^d \colon A v \leq \mathbf{1} \}$ for some \emph{integral} matrix $A$ (where $\mathbf 1$ is the all-one vector). If $P$ is the convex hull of finitely many vectors $v_1,\dots,v_r$, then the dual is the intersection $\bigcap_{i=1}^r v_i^*$. Reflexive polytopes have precisely one interior integral point, and a result of Lagarias and Ziegler \cite{LagariasZiegler} implies that in fixed dimension $d$, there are only finitely many reflexive polytopes up to translations and unimodular linear transformations (i.e., elements from $\mathrm{GL}(d,\Z)$). 

\subsubsection*{Computing the cyclopolytopes and cyclovarieties} 

The cyclopolytope $N_k$ is the convex hull of the set $$A_k = \{ \mathbf{a}^{(r)}=(a^{(r)}_0,\dots,a^{(r)}_{d-1}) \colon r=0,\dots,k-1 \} \mbox { where } d=\phi(k) \mbox{ and } \zeta_k^r = \sum_{s=0}^{d-1} a_s^{(r)} \zeta_k^s. $$
In coordinate-free notation, $N_k = \mathrm{conv}(\mu_k) \subseteq \Z[\zeta_k] \otimes_{\Z} \R$. 

\begin{remark} Notice that this representation depends on the choice of a primitive $k$-th root of unity, but a different choice will only change the resulting polytope by an invertible integral linear transformation of $\R^d$ (and the resulting variety by an isomorphism). 
\end{remark} 

Let $\Phi_k(t)$ denote the $k$-th cyclotomic polynomial. The vectors in $A_k$ can be computed as the coefficients in the remainder of division of $t^r$ by $\Phi_k(t)$.

\begin{remark} \label{algocycsage} See Section~\ref{sec:code} for a quick routine to compute cyclopolytopes exploiting the built-in cyclotomic field functionality in \texttt{SageMath}. 
In \texttt{SageMath}, one can test reflexivity of $N_k$ by using the PALP-library, see \cite{PALP}, but that algorithm will terminate in reasonable time only in small dimensions. 
\end{remark} 

 An equivalent procedure that directly computes the cyclovariety (from which the cyclopolytope can be read off as exponent vectors) is the following. Uniquely define the transformation $T \colon \Z[t] \rightarrow \{0,1\}[X_1^{\pm 1}, X_2^{\pm 1}, \dots]$ from the polynomials over the integers to the Laurent polynomials in countably many variables with all non-zero coefficients equal to $1$ by 
the properties that  it is a morphism of the additive monoid to the multiplicative monoid (i.e.\ $T(a+b)=T(a)T(b)$ for all $a,b \in \Z[t]$) and that  
$T(t^n)=x_{n+1}$. 
The cyclovariety is the subscheme of $\Gm^{k+1}$ defined by the ideal $$ \langle \sum_{i=0}^{k} x_i, T(\Phi_k(t))-1,T(t\Phi_k(t))-1,\dots, T(t^{k-d-1}\Phi_k(t))-1 \rangle.$$ Eliminating all variables $x_i$ with $i>d$ from these equations gives the hypersurface equation of the cyclovariety in $\Gm^{d+1}$.

\begin{remark} 
The cyclopolytopes and cyclovarieties occur in many guises in the literature, e.g., in \cite{Myerson} (as in this paper), \cite{BeckHosten} (with a different definition, in connection with word lengths and transportation polytopes), \cite{Reiner} (as matroid), \cite{KlyachkoVoskresenskii} (in connection with algebraic groups), \cite{Lala} (in connection with planar walks). 
\end{remark} 

We now present some examples. First, notice the following. 

\begin{lemma} 
If $k$ is even, then $P_k(x_0, \mathbf{x})$ is symmetric, i.e., $$P_k(\mathbf{x}) = x_0 + G_k(x_1,\dots,x_d) + G_k(x_1^{-1},\dots,x_d^{-1})$$ for some polynomial $G_k \in \Z[\mathbf x]$. 
\end{lemma}

\begin{proof} 
If $k$ is even, $H$ contains complex conjugation ($-1 \in (\Z/p)^*$), so $$\alpha_n = \sum_{\sigma \in H/\langle -1 \rangle} \left( \zeta_p^\sigma + \zeta_p^{-\sigma} \right),$$ from which the result follows. 
\end{proof}

Therefore, it is useful to introduce the following notation. 

\begin{notation} 
For a Laurent polynomial $P$ of the form $$P=x_0 + F(x_1,\dots,x_d) \in \Z[x_0^{\pm 1},\dots,x_d^{\pm 1}]$$ we denote by $P+\sym$ the following `symmetrized' Laurent polynomial
$$ P+\sym\coloneqq  x_0 + F(x_1,\dots,x_d) + F(x_1^{-1},\dots,x_d^{-1})  \in  \Z[x_0^{\pm 1},\dots,x_d^{\pm 1}].$$
\end{notation}

Notice that the result of symmetrisation is invariant under $x_i \rightarrow 1/x_i$ for all $i \geq 1$, but not for $i=0$. 

\begin{example} 
\label{N2k} 
If $k=2^r$ for some $n \geq 1$, $\Phi_k(t) = t^{2^{r-1}}+1$, and the cyclovariety is defined by
$$ x_0+x_1+\dots+x_{2^{r-1}} + \sym. $$
In fact, a cyclovariety $C_k$ has equation $x_0+x_1+\dots+x_m +\sym$ for some $m$ precisely if $m=\phi(k) = k/2$, which happens if and only if $k=2^r$ for some $r \geq 1$. 
\end{example}

\begin{example} 
\label{Nq} 
If $k=q$ for an odd prime $q$, $\Phi_k(t) = (t^q-1)/(t-1)$, and the cyclovariety is defined by the (non-symmetrized) equation
$$ x_0+x_1+\dots+x_{q-1} + \frac{1}{x_1 \cdots x_{q-1}}. $$
\end{example}

 \begin{example} 
 \label{N2qn} 
 If $k=2 q^r$ for an odd prime $q$, $\Phi_{k}(t)=\Phi_q(-t^{q^{r-1}}) = t^{q^r-q^{r-1}} - t^{q^r-2q^{r-1}} + \dots + 1$ with $d=q^r-q^{r-1}$, and the cyclovariety $C_k$ is defined by the Laurent polynomial
$$x_0 + x_1 + \dots + x_{q^r-q^{r-1}} + \sum_{j=1}^{q^{r-1}} \prod_{i=2}^q x_{q^r-iq^{r-1}+j}^{(-1)^i} +  \sym. $$
A particularly important case occurs when $k=2q$ for an odd prime $q$, and then $C_k$ is defined by 
\begin{align} \label{N2ell} & x_0 + x_1 + \dots + x_{q-1} + \frac{x_1 \cdot x_3 \cdots x_{q-2}}{x_2 \cdot x_5 \cdots x_{q-1}}+  \sym \nonumber \\ 
& = x_0 + x_1 + \frac{1}{x_1} + \dots + x_{q-1}+\frac{1}{x_{q-1}} + \frac{x_1 \cdot x_3 \cdots x_{q-2}}{x_2 \cdot x_4 \cdots x_{q-1}} + \frac{x_2 \cdot x_4 \cdots x_{q-1}}{x_1 \cdot x_3 \cdots x_{q-2}}. \end{align} 
\end{example}

\subsubsection*{Results} 

\begin{proposition} \label{polNk} For all even $k$, $N_k$ is a centrally symmetric polytope. 
\begin{enumerate} 
\item If $k=2^r$ is a power of $2$, the Newton polytope $N_k$ is the \emph{cross-polytope} in dimension $d=2^{r-1}$, i.e., the convex hull of $\{\pm e_i\}_ {i=1}^{d}$, where $e_i$ are the standard base vectors in $\R^d$. This is also the direct sum of $d$ copies of the interval $[-1,1]$. It is reflexive with dual polytope the $d$-dimensional hypercube. 

\item \label{DPpol} If $k=2q$ where $q$ is an odd prime, the Newton polytope $N_k$ is a \emph{del Pezzo polytope}, i.e., the convex hull of $\{\pm e_i\}_ {i=1}^{q-1} \cup \{\pm v\}$, where $e_i$ are the standard base vectors in $\R^{q-1}$ and $v=(1,-1,1,\dots,-1) \in \R^{q-1}.$  It is reflexive. 

\item \label{polNksum} If $k=2q^r$ where $q$ is an odd prime, the Newton polytope $N_k$ is the direct sum of $q^{r-1}$ copies of the del Pezzo polytope $N_{2q}$.  It is reflexive. 
\end{enumerate} 
\end{proposition} 

\begin{remark} 
The del Pezzo polytope is usually defined by using the vectors $\pm(1,\dots,1)$, leading to the same polytope as in the above formulation up to a unimodular transformation. 
\end{remark}

\begin{proof} If $k$ is even, the set of all $k$-th roots of unity $\mu_k$ is invariant under $x \mapsto -x$, hence so is $N_k$. 

The shape of the Newton polytopes follows immediately from the explicit expressions of the cyclovarieties given in Examples \ref{N2k} and \ref{N2qn}. 

The reflexivity of these polytopes is well-known, but we include the following brief discussion on how to find the actual dual polytopes. 

If $k=2^r$, then $N_k$ is the convex hull of $\{\pm e_i\}_ {i=1}^{d}$ with $d=2^{r-1}$.  The duals of the standard basis vectors are $e_i^* = \{ x \in \R^{d} \colon x_i \leq 1\}$ and $(-e_i)^* = \{ x \in \R^{d} \colon x_i \geq -1\}$. Thus, $\bigcap (\pm e_i)^*$ is precisely the convex hull of all vectors of the form $(\pm 1, \pm 1, \dots, \pm 1)$, a side-length-two \emph{hypercube}. Since it is integral, $N_k$ is reflexive.  

If $k=2q$ where $q$ is an odd prime, then $N_{2q}$ is the convex hull of $\{\pm e_i\}_ {i=1}^{q-1} \cup \{\pm v\}$, where $e_i$ are the standard base vectors in $\R^{q-1}$ and $v=(1,-1,1,\dots,-1) \in \R^{q-1}$. We compute a set of integral points whose convex hull is the polar dual $N_\ell^*$. As before, $\bigcap (\pm e_i)^*$ is the side-length-two hypercube spanned by $(\pm 1, \pm 1, \dots, \pm 1)$. The dual of the vector $v$ is $v^*=\{ x \in \R^{\ell-1} \colon x_1-x_2+x_3-x_4+\dots + x_{q-2} - x_{q-1} \leq 1\}$. Hence the dual $N^*_{2q}$ is the part of the hypercube in between the boundary hyperplanes of $v^*$ and $(-v)^*$, i.e., \begin{equation} \label{ineqv}  -1 \leq x_1-x_2+x_3-x_4+\dots + x_{q-2} - x_{q-1} \leq 1. \end{equation} 
The hypercube intersects $v^*$ is the hyperplane spanned by $\{e_1,-e_2,\dots,e_{q-2}, -e_{q-1}\}$ and it intersects $(-v)^*$ in the hyperplane spanned by $\{-e_1,e_2,\dots,-e_{q-2},e_{q-1}\}.$
Thus, $N^*_{2q}$ is the convex hull of the set $\{\pm e_i\}$, together with all hypercube corners $(\pm 1, \pm 1, \dots, \pm 1)$ that satisfy the inequalities in \eqref{ineqv}, and since all these spanning points have integral coordinates, the dual lattice is integral, too. We conclude that $N_{2q}$ is reflexive. 

Finally, for $k=2q^r$ with $q$ an odd prime, $N_{2q^r}$ is spanned as the convex hull of the set $\{\pm e_i\}_{i=1}^{q^r-q^{r-1}}$, together with points of the form 
 $$ \pm (\underbrace{0, \dots, 0}_{j},1, \underbrace{0, \dots, 0}_{q^{r-1}}, -1, \underbrace{0, \dots, 0}_{q^{r-1}}, 1,\dots)$$
 for $j=0,\dots,q^{r-1}-1$, and this is indeed precisely $N_{2q}^{q^{r-1}}$. Since the direct sum of reflexive polytopes is reflexive (the dual of a direct sum is the product of duals), we conclude that $N_{2q^r}$ is also reflexive. 
 \end{proof} 
 
\begin{remark}  
The Newton polytope of $C_k$ itself (including the variable $x_0$)  is a pyramid over the Newton polytope $N_k$. More precisely, under the standard embedding  $\iota \colon \R^{d} \subseteq \R^{d+1}$ mapping $v$ to $(v,0)$, it is the convex hull of $\iota(N_k) \cup \{e_{d+1}\}$, where $e_{d+1}$ is the new basis vector. Thus, the Newton polytope of $C_k$ is $\iota(N_k)^* \cap e_{d+1}^*$, which is the convex hull of the image under $\iota$ of the spanning points of $N_k^*$, together with the point $e_{d+1}$. Since these are all integral if the points in $N_k^*$ are, the Newton polytope of $C_k$ is itself reflexive when $N_k$ is reflexive. 
\end{remark}  

\begin{remark} \label{guessreflexive} 
One may wonder whether $N_k$ is reflexive for all $k$. Reflexivity, introduced by Batyrev in his study of mirror symmetry, is a rare polytope property. We checked in PALP that this it holds for $N_k$ with $k=2\ell$ and $\ell=6,10,12,15$. An interesting challenge occurs for $k=3\cdot 5 \cdot 7 = 105$ and the related centrally symmetric case $k=2 \cdot 105 = 210$. Beck and Ho\c{s}ten have shown that, for distinct primes $p_1,\dots,p_r$, $N_{p_1\cdots p_r}$ is \emph{combinatorially} dual (i.e., with the posets given the face inclusions being dual) to a generalized transportation polytope $T(p_1,\dots,p_r)$ \cite[Theorem 18]{BeckHosten}, and they remark that Seth Sullivant has computed enough faces of $T(3,5,7)$ to show that it is not integral (see the arguments following \cite[Conjecture 4]{BeckHosten}). We have checked in \texttt{SageMath} (see Section~\ref{sec:code}) that $T(3,5)$ is not unimodularily equivalent to the \emph{polar}  dual of $N_{3\cdot 5}$, so the question of reflexivity remains open even in this case. Recall that combinatorial duality is not a metric property, but polar duality is. 
\end{remark} 

\begin{remark} The cyclovariety $C_k$ equals the \emph{Bloch variety} (named after the physicist Felix Bloch) of the lattice Laplacian on the $\R^d$ lattice with fundamental cell $N_k$, in which $-x_0$ is the energy of the system. The fibers of this variety for fixed $x_0$ are called \emph{Fermi varieties} for a fixed energy. Our example $C_6$ is closely related to a honeycomb crystal structure, and the tight-binding model for electrons placed on this lattice is the simplest theory of graphite/graphene, see Section \ref{MC6}. 

These varieties have been studied by methods of algebraic geometry in the one-dimensional case in the monograph  \cite{Fermibook}, and questions of irreducibility of fibers have been studied in great generality  in \cite{FillmanBloch} and \cite{FillmanFermi}. In our examples, the generic connectedness and irreducibility will follow from their description via specific linear systems. 
\end{remark}

\begin{remark} The equations have apparent symmetries, given by interchanging or inverting the variables. For $k=4,6$, these produce automorphisms of the corresponding cubics, some of which are given by translations over torsion points, once an elliptic curve structure has been chosen. The elliptic curves do not admit complex multiplication in general. It would be interesting to understand the automorphism groups of the cyclovarieties in general. 
\end{remark}
 
\subsection{Toric geometry related to the cyclovariety} 

Assume that $P$ is a reflexive polytope in dimension $d$. By considering the face fan of $P$, to it one associates a projective (possibly singular) toric Fano variety $\PP_P$, compactifying a given $\Gm^d$. If $P$ is centrally symmetric, then $\PP_P$ admits an involution $\sigma$ such that $\sigma(tx)=t^{-1}\sigma(x)$ for all $t \in \Gm^d$. 

\begin{lemma} \label{dp} Let $\PP_k=\PP_{N_k}$ denote the projective toric variety associated to the polytope $N_k$. 
\begin{enumerate}
\item \label{dp2} If $k=2q$ where $q$ is an odd prime, $\PP_{k}$ is the $(q-1)$-dimensional \emph{toric del Pezzo variety}, a smooth Fano variety; i.e., the anticanonical divisor is ample. 
\item \label{dp3} If $k=2q^r$ with $q$ an odd prime, $\PP_{q^n} = \PP_{2q}^{q^{n-1}}$ is a smooth Fano variety. 
\end{enumerate}
\end{lemma} 

\begin{remark} 
The toric del Pezzo surfaces $\PP_{2\ell}$ were introduced by Klyachko and Voskresenskii \cite{KlyachkoVoskresenskii} (denoted by $V^{2\ell}$ in that reference) as basic building block in their classification of centrally symmetric smooth Fano polytopes and are related to root systems. They can also be defined by blowing up $(\PP^1)^{2 \ell}$ in two torus-invariant points; compare \cite[VII.8.6]{Ewaldbook}. 
\end{remark}

\begin{proof} 
For \eqref{dp2}, $\PP_k$ is constructed from the del Pezzo polytope in Proposition \ref{polNk}\eqref{DPpol}, see \cite[p.\ 234]{KlyachkoVoskresenskii}, including smoothness.  For the correspondence between the smooth Fano property for polytopes and varieties, see, e.g., \cite[V.8.3-8.4; VII.8.4-8.6; pp.\ 329--330]{Ewaldbook}.

For \eqref{dp3}, note that direct sum of polytopes corresponds to direct product of the corresponding toric varieties \cite[pp.\ 329--330]{Ewaldbook}, so the result follows from Proposition \ref{polNk}\eqref{polNksum}, and that a product of smooth Fano varieties is a smooth Fano variety. 
\end{proof}

\subsection{The log Calabi--Yau property for certain $k$}

In this subsection, we first prove Theorem \ref{mainCY}, and then turn to compute the cohomological properties of the compactified resolved pencil. Let $k=2q^r$ for an odd prime $q$. The proof is given by a sequence of constructions, and we will illustrate this by explicitly dealing with the case $k=6$. To distinguish better between the parameter of the family $x_0$ and the coordinates $x_1,\dots,x_d$, in this subsection we will use the notation $\lambda \coloneqq  x_0$.

\begin{proof}[Proof of Theorem \ref{mainCY}] The (affine) cyclovariety $C_k(\lambda) \subseteq \Gm^d$ is given by 
$$ C_k(\lambda) \colon \lambda + F_k(\mathbf x) = 0. $$
We define the projective cyclovariety $\overline C_k(\lambda)$ as the closure of 
$ C_k(\lambda)$ in $\PP_k$. Since $N_k$ is reflexive, $\PP_k$ is Fano, so $-K_{\PP_k}=\mathcal O_{\PP_k}(1)$ is ample \cite[Thm.\ 4.1.9]{Batyrev}. Since $\PP_k$ is projective smooth toric  by Lemma \ref{dp}\eqref{dp2}, $-K_{\PP_k}$ is then automatically very ample (Demazure, see \cite[Cor.\ 2.15]{Oda}).

Note that since $F_k$ is a Laurent polynomial, in the toric compactification $\PP_k$ we should consider this pencil after `clearing denominators': if we set $$D(\mathbf x) \coloneqq  \prod_{i=1}^d x_i^{\max \{|a_i| \colon a_i \leq 0\}}, $$ it is the pencil in $\PP_k$ defined by the (generally incomplete) linear subsystem of the complete anticanonical linear system 
$$ \mathfrak d  \coloneqq  \langle D(\mathbf x), D(\mathbf x) F_k(\mathbf x) \rangle \subseteq |-K_{\PP_k}|. $$
 The incomplete linear system $\mathfrak d$ has base points defined by the simultaneous vanishing of $D$ and $DF_k$, which defines a subscheme $B$ of the toric boundary $\PP_k \setminus \Gm^d$. Infinitely many fibers might have singularities at these boundary components, so we need to resolve the base point locus in order to be able to apply the theorem of Bertini to the entire fiber (by that theorem, the fibers are generically smooth outside the base point locus).

We now take a smooth resolution of the base point locus $B$ of $\mathfrak d$, i.e., a composition of blow-ups such that the inverse image of the ideal sheaf of the base locus becomes invertible; existence follows from the principalisation theorem, see, e.g., \cite[Thm.\ 2.5]{BEV}. This produces a smooth variety $\widetilde \PP_k \overset{\pi}{\rightarrow} \PP_k$, and since $B$ is of codimension two in $\PP_k$, we find that 
\begin{equation} \label{blowupK} K_{\widetilde \PP_k} = \pi^* K_{\PP_k} + E, 
\end{equation} 
where $E = \pi^{-1}(B)$ (see, e.g.,  \cite[Ex.\ II.8.5]{Hartshorne} or \cite[I\S 4, p.\ 187]{GH}). We define the resolved family $\widetilde C_k \rightarrow \PP^1_\lambda$ by pullback \begin{equation} \label{deftildec} \widetilde C_k(\lambda) \coloneqq  \overline C_k(\lambda) \times_{\PP_k} \widetilde \PP_k. \end{equation} This gives a base-point free pencil $\widetilde{\mathfrak d}$ defined by the pull-back of the original linear system $\mathfrak d$ along $\pi$, i.e., defined by the divisor $E$ corresponding to the invertible ideal sheaf defining the blow-up $\widetilde \PP_k$ (compare \cite[Ex.\ II.7.17.3]{Hartshorne}), with $$ \widetilde{\mathfrak d} \subseteq |- \pi^* K_{\PP_k} | = |-K_{\widetilde {\PP}_k}+E|. $$ 
We apply the theorem of Bertini \cite[I \S1, p.\ 157]{GH} to this, now base-point free, pencil: we find that, for all but finitely many $\lambda$, $\widetilde C_k(\lambda)$ is smooth and irreducible. 

Since for the resolution of the base locus, we iteratively blow-up only in the toric boundary, $C_k(\lambda)$ is isomorphic to its strict transform $\pi^{-1}(C_k(\lambda))$ \cite[II.7.13(b)]{Hartshorne}. By the adjunction formula \cite[I \S 1, p. 147]{GH} and using \eqref{blowupK}, we find that $$K_{\widetilde C_k(\lambda)} = (K_{\widetilde \PP_k} - \pi^* K_{\PP_k})|_{\widetilde C_k(\lambda)} = E|_{\widetilde C_k(\lambda)}$$ hold for all but finitely many $\lambda$. Furthermore, in the principalisation theorem, it is guaranteed that $E|_{\widetilde C_k(\lambda)}$ is a normal crossing divisor. 
\end{proof} 

\begin{remark} 
	\label{rem:complete-anticanonical}
	If we consider instead the complete anticanonical linear system (so the analogue of $F_k$ but with varying coefficients in front of all non-trivial monomials, leading to a higher-dimensional parameter space), the generic element will be smooth; this is the set-up of `regularity' in \cite{Batyrev}, which is rather different from our one-dimensional parameter space. 
\end{remark}

\begin{example} \label{AGex} 
For $k=6$, $\PP_6$ is the del Pezzo variety given by blowing up $\PP^2$ in three general points. The anticanonical embedding in $\PP^6$ is given by 
$$ \PP_6 = \{ (b:a_1:\dots:a_6) \colon ba_i = a_{i-1} a_{i+1}, a_i a_{i+3}=b^2 \} $$
(see, e.g., \cite[Thm.\ 8.4.1]{Dol}), 
and our pencil is given by intersecting with the variable hyperplane $$H_\lambda \coloneqq  \{ \lambda b + \sum a_i = 0 \},$$ 
which has base point locus $B \coloneqq  \{ (0:a_1:\dots:a_6) \colon \sum a_i = 0 \}\cap \PP_6 \subseteq \PP^6$ of codimension $2$. The fibers of the family $\overline C_6(\lambda)$ can contain singular points in the boundary locus $b=0$. 

The resolution, given by blowing up 6 smooth points in the toric boundary of $\PP_6$, is worked out in \cite[\S 5.7]{Grassi}, including a description of the remaining singular fibers. In this case, the family defines a rational elliptic surface with Mordell--Weil lattice $ A_1 \oplus A_2 \oplus A_5$ of index $6$ in $E_8$ (number 66 in \cite[Table 8.2]{SS}) so the surface has Mordell--Weil group $\Z/6\Z$. It turns out there are four singular fibers, of type $I_1,I_2,I_3$ and $I_6$. Note that the original equations for $C_6(\lambda)$ give a family of affine models of elliptic curves whose discriminant, as a function of $\lambda$, vanishes at $\lambda=2,3,-6$. 
\end{example}

\begin{remark} Another formulation of Theorem \ref{mainCY} is that the original family $C_k(\lambda)$ is a \emph{pencil of log Calabi-Yau varieties}, i.e., the fibers with smooth compactification $\widetilde C_k(\lambda)$ satisfy the conditions in \cite{GHK} (in the surface case, this is also known as a `Looijenga pair' \cite{Loo}). 
\end{remark} 

\subsubsection*{Cohomology and topology of the compactified resolved pencil} Finally, we compute the fundamental group and cohomology of the structure sheaf on the smooth compactified fibers. We see that all middle cohomology vanishes, a requirement that for usual Calabi--Yau varieties is called `strictness'. 

\begin{proposition} \label{cohCk} The smooth fibers of the compactification $\widetilde C_k(\lambda)$ of $C_k(\lambda)$, of dimension $d-1$, have the following cohomology groups: 
$$ \mathrm H^i(\widetilde C_k(\lambda), \mathcal O_{\widetilde C_k(\lambda)}) =  \left\{ \begin{array}{ll} \Cc  & \mbox{ if } i=0 \mbox{ or } i=d-1, \\ 0 & \mbox{ if } i=1,\dots,d-2 \mbox{ or } i\geq d; \end{array} \right.$$
in particular, their Euler characteristic $\chi(\mathcal O_{\widetilde C_k(\lambda)})=0$ is zero, and their  Hodge numbers satisfy $$h^{p,0}(\widetilde C_k(\lambda))=h^{0,p}(\widetilde C_k(\lambda))=0 \mbox{ for } 0<p<d-1.$$ Furthermore, for $k \neq 6$, they are simply connected \textup{(}i.e., $\pi_1(\widetilde C_k(\lambda))=0$\textup{)}, whereas $\pi_1(\widetilde C_6(\lambda)) \cong \Z$. 
\end{proposition} 

\begin{proof} For the sake of simplifying notation, for the duration of this proof we set 
$$ P \coloneqq  \PP_k;\ \ \widetilde P \coloneqq  \widetilde \PP_k; \ \ C \coloneqq  \widetilde C_k(\lambda)$$
(note that the letter `$C$' refers to the cyclovariety, not to a `curve').  
The result is obvious for $i=0$ (by connectedness of the fibers) and $i\geq d$ (since $C$ is smooth projective of dimension $d-1$), so we focus on the remaining values of $i$. 
If we denote by $j \colon C \hookrightarrow \widetilde P$ the embedding defined by the linear system $|-\pi^*K_P|$, then we have the following short exact sequence of sheaves on $\widetilde P$:  
$$ 0 \rightarrow \mathcal O_{\widetilde P}(\pi^*K_P) 
\rightarrow \mathcal O_{\widetilde P} \rightarrow j_* \mathcal O_C \rightarrow 0, $$
which, using \eqref{blowupK}, leads to a long exact sequence in cohomology
\begin{equation} \label{themiddle} \dots  \rightarrow \mathrm H^i(\widetilde P,\mathcal O_P) \rightarrow \mathrm H^i(\widetilde P,j_* \mathcal O_C) \rightarrow \mathrm H^{i+1}(\widetilde P,\mathcal O_{\widetilde P}(\pi^*K_{P})) \rightarrow \dots,  \end{equation} 
where the middle cohomology group equals $\mathrm H^i(C, \mathcal O_C)$ \cite[III Lem.~2.1]{Hartshorne}, the group we wish to compute. To compute the outside cohomology groups, notice that the resolution $\pi \colon \widetilde P\rightarrow P$ is done via the principalisation theorem, i.e., through iterated blow-ups in smooth centers (though not necessarily in a torus-invariant center, so $\widetilde P$ might not be a toric variety). 

We now interject the following observation: since $P$ is a smooth toric projective variety, it is birational to a projective space, and since $\tilde P$ is given by iterated blowups, it is birational to $P$. Hence $\widetilde P$ is a projective smooth rational variety, and thus, it is simply connected (both algebraically and topologically) and it has vanishing Hodge numbers $h^{p,0}(\widetilde P) = \dim \mathrm H^0(\widetilde P, \Omega_{\widetilde P}^p)$ for all $p>0$ (see, e.g.,  \cite[Cor.~4.18]{Debarre}). By symmetry of the Hodge diamond \cite[pp.\ 116--117]{GH}, we also find $$ \dim \mathrm H^q(\widetilde P, \mathcal O_{\widetilde P})= h^{0,q}(\widetilde P) = h^{q,0}(\widetilde P) = 0 $$ for $q>0$. This shows that the leftmost cohomology groups in \eqref{themiddle} vanish for $i>0$. If one prefers a proof avoiding Hodge theory, one may reason as follows: since we are blowing up a smooth variety, a result of Hironaka \cite[Cor. 2, p. 153]{Hironaka} implies the vanishing of higher direct images of the structure sheaf: $\mathrm R^i \pi_* \mathcal O_{\widetilde P} = 0$, for $i>0$, with $\pi_* \mathcal O_{\widetilde P} = \mathcal O_P$. Thus (a degenerate case of) the Leray spectral sequence \cite[III.~Ex.\ 8.1]{Hartshorne} gives 
$ \mathrm H^i(\widetilde P, \mathcal O_{ \widetilde P}) = \mathrm H^i(P, \pi_* \mathcal O_{ \widetilde P}) = \mathrm H^i(P, \mathcal O_P). $ 
 Since $P$ is a smooth toric variety, Demazure vanishing \cite[Theorem 9.2.3]{Cox} says $\mathrm H^i(P, \mathcal O_P) = 0$ for $i>0$.
 
For the rightmost cohomology groups, we  mimic the above algebraic proof for another sheaf, and for this, we use the Grauert--Riemenschneider vanishing theorem \cite[Satz 2.3]{Grauert} in combination with \eqref{blowupK} to conclude that  $$\mathrm R^i \pi_* \pi^* K_P= \mathrm R^i \pi_* (K_{\widetilde P}-E) =0 \ \ (i>0)$$ (since $\mathcal O_{\widetilde P}(-E)$ is ample).  Again, from the Leray spectral sequence, we find that $$\mathrm H^{i+1}(\widetilde P, \pi^* K_P) = \mathrm H^{i+1}(P,\pi_* \pi^* K_P) \ \ (i>0).$$ Now the `projection formula' \cite[II.~Ex.~5.2(d)]{Hartshorne} gives $\pi_* \pi^* K_P \cong \pi_*\mathcal O_{\widetilde P} \otimes_{\mathcal O_P} K_P \cong K_P, $ (where we have used $ \pi_*\mathcal O_{\widetilde P} = \mathcal O_P$), and we find 
$$ \mathrm H^{i+1}(P,\pi_* \pi^* K_P) = \mathrm H^{i+1}(P,K_P) = \mathrm H^{d-i-1}(P,\mathcal O_P), $$
by Serre duality. Again, Demazure vanishing gives that $\mathrm H^{d-i-1}(P,\mathcal O_P) = 0$ for $d-i-1>0$. 

Putting everything together, we see that for $i=1,\dots,d-2$, both cohomology groups on the left and right in \eqref{themiddle} vanish, and the result follows for these values of $i$. Finally, for $i=d-1$, we get $\mathrm H^{d-1}(C,\mathcal O_C) \cong \mathrm H^0(P,\mathcal O_P) = \Cc$. 

That the Euler characteristic, as alternating sum of dimensions of the cohomology groups, is zero, follows from the fact that the dimension $d=\varphi(k)$ is even for $k>2$.  Since $\mathrm H^{0,q}(C) = \mathrm H^q(C,\mathcal O_C)$, it follows that $h^{0,q}=0$ for $1<q<d-1$; and we have $h^{q,0}=h^{0,q}$ by the aforementioned symmetry in the Hodge decomposition. 

We have already noticed that $\widetilde P$ is simply connected (this also follows by birational invariance of the fundamental group for regular proper varieties \cite[X, Cor.~3.4]{SGA1}). By the Lefschetz hyperplane theorem, $$\pi_1(C) \cong \pi_1(\widetilde P)$$ for $d \geq 3$ (see, e.g., \cite[Thm.\ 7.4]{Milnor}). With $d=\phi(k)$ and $k=2q^r$, we find this holds if $k \neq 6$. On the other hand, for $k=6$, we have a pencil of cubic curves, that have fundamental group isomorphic to $\Z$. 
\end{proof} 

\begin{remark} 
As Beauville has shown \cite[Prop.\ 2]{Beauville}, a smooth projective variety $M$ of dimension $m$ admits a K\"ahler metric with holonomy exactly equal to $\mathrm{SU}(m)$ (i.e., maximal) precisely if the following purely `algebro-geometric' conditions hold: the canonical bundle $K_M=0$ is trivial, and $h^{0,p}(M')=0$ for $0<p<m-1$ for all finite etale covers $M' \rightarrow M$. Here, we have trivial fundamental group, and, by our result, the conditions on the Hodge numbers holds. Only the condition $K_M=0$ is replaced by a `log' version after the resolution of the base locus. In this sense, Proposition \ref{cohCk} is an algebraic analogue of the requirement of `maximal holonomy' sometimes used to define Calabi-Yau manifolds in the strict sense. 
\end{remark}

\section{The associated algebraic dynamical system; proof of Theorem \ref{mainDyn}} \label{DS} 

\subsection{Construction of the system} \label{constDS} 
By a result of Schmidt, Lind and Ward \cite[Thm.~3.1]{SchmidtLindWard}, the Mahler measure $\m(C_k)$ is the topological entropy of an algebraic dynamical system, consisting of an action of the group $\Z^{d+1}$ on a topological group $X_k$ that is constructed as follows (cf.\ \cite[\S 11]{KitchensSchmidt}). We consider the Laurent series ring in $d+1$ variables $R_{d+1}\coloneqq \Z[x_0^{\pm 1}, x_1^{\pm 1},\dots,x_d^{\pm 1}]$ and its quotient $M_k\coloneqq R_{d+1}/(P_k(x_0,\mathbf{x}))$. We let $X_k=\widehat{M_k}$ denote the Pontrjagin dual of the discrete additive group of $M_k$, i.e.,  $$X_k=\widehat{M_k}\coloneqq \mathrm{Hom}(M_k^+,\T),$$ which is a compact commutative topological group, and the action of $(n_0,\mathbf{n}) \in \Z^{d+1}$ is given as the dual of the multiplication by the monomial $x_0^{n_0} \mathbf{x}^{\mathbf{n}}$ on $M_k$, i.e., a homomorphism $\varphi \colon M_k \rightarrow \T$ is sent to $\varphi( x_0^{n_0} \mathbf{x}^{\mathbf{n}} \cdot \text{--})$. We remark that in this case, the topological entropy also equals the measure-theoretic entropy for the Haar measure on $X_k$ \cite[13.3]{Schmidt}. 

\subsection{The dictionary between dynamical and algebraic properties} 
As in the monograph \cite{Schmidt} of Schmidt, the dynamical properties of this system translate into algebraic properties of $C_k$, so we now run through a list of these, very similar to the table found in \cite[Fig.~1, p.\ xiii]{Schmidt}. We assume $k > 2$, so $d \geq 2$. 

\subsubsection*{Dense periodic points; ergodicity} First of all, notice that the evaluation homomorphism `substitute $x_0$ by $-F_k(\mathbf{x})$' is a surjective ring homomorphism from $R_{d+1}$ onto $R_d$ (where the latter ring is in the Laurent variables $\mathbf x = (x_1,\dots,x_d)$) with kernel $(P_k)$. This gives an isomorphism \begin{equation} \label{mkrd} M_k \cong R_d. \end{equation} Since $R_d$ is a domain, $(P_k)$ is a prime ideal. Since $R_{d+1}$ is a UFD, this implies that $P_k$ is irreducible, a fact which can also be established directly using the shape $x_0+F_k(\mathbf{x})$ of $P_k$. The units of $R_{d+1}$ are the monomials with coefficient $\pm 1$.  

Since $M_k = R_{d+1}/(P_k)$ is Noetherian, the dynamical system satisfies the descending chain condition and periodic points are dense \cite[5.4 \& 5.7]{Schmidt}. Since $d \geq 2$, the fact that $M_k$ is defined by a principal ideal $(P_k)$ implies that the action is ergodic for the Haar measure \cite[6.11]{Schmidt}. 

\begin{remark} 
The system can either be described as a multidimensional shift on a subspace of a $(d+1)$-dimensional regular square grid, or as a more complicated function on a full $d$-dimensional grid. First of all, as in \cite[5.2(2)]{Schmidt}, the system is isomorphic to the full shift restricted to the subspace of $\T^{\Z^{d+1}}$ given by 
$$ X_k \coloneqq  \{ x \in \T^{\Z^{d+1}} \colon \sum_{ \mathbf a \in A_k} x_{(m_0,\mathbf m)+(0,\mathbf a)} = - x_{(m_0,\mathbf m)+(1,\mathbf 0)} \mbox{ mod } 1, \forall (m_0,\mathbf m) \in \Z^{d+1} \}. $$ Secondly, we can use the isomorphism \eqref{mkrd}; this conjugates the action on $X_k$ into the action of $ (n_0,\mathbf n) \in \Z^{d+1}$ on $\widehat{R_d} = \T^{\Z^d}$ given by the dual of multiplication by the Laurent polynomial $F_k(\mathbf x)^{n_0} \mathbf x^{\mathbf n}$. 
\end{remark}

\subsubsection*{Entropy; mixing} We claim that the system has positive entropy, or, equivalently, $\m(C_k)>0$ for $k \geq 2$. On the algebraic side, this is implied by the fact that $(P_k)$ is not generated by a generalized cyclotomic polynomial (meaning a Laurent polynomial of the form $x_0^{m_0} \mathbf{x}^{\mathbf{m}} c(x_0^{m_0} \mathbf{x}^{\mathbf{n}})$ for some $(m_0,\mathbf{m}), (n_0, \mathbf{n}) \in \Z^{d+1}$ with $(n_0, \mathbf{n}) \neq (0,\mathbf{0})$ and $c$ a one-variable polynomial with only roots of unity as roots)---cf.\ the multivariable analogue of Kronecker's theorem (\cite[19.5]{Schmidt}, see also \cite{BoydKronecker}, \cite{SmythKronecker}, \cite{LawtonKronecker}). Since $P_k$ is irreducible monic, this is in fact equivalent to $P_k$ not being generalized cyclotomic itself (and irreducible as Laurent polynomial). Since the entropy equals $\m(C_k)$, the positivity is proven in Lemma \ref{pos}. 
Since $P_k$ is not generalized cyclotomic, we also conclude that the dynamical system is mixing of all orders \cite[6.12 \& 27.2]{Schmidt}. 

Since the system is defined by one prime ideal, it has completely positive entropy \cite[20.8]{Schmidt}, and hence it is Bernoulli \cite[23.1]{Schmidt} (i.e., measurably isomorphic to a Bernoulli shift; this does not imply topological/algebraic conjugacy). 

Finally, the Haar measure is the unique measure of finite maximal entropy \cite[20.15]{Schmidt}.  

So far, we have established all claimed properties in the first part of Theorem \ref{mainDyn}. We now turn to proving the second part. 

\subsection{Non-expansiveness for certain $k$} The system is not expansive if and only if $C_k(\Cc) \cap \T^{d+1} \neq \emptyset$. We believe this is true for all $k>1$ but we prove it under some assumption, as follows. 

\begin{proposition} \label{toricpointsDS} 
If $k$ is a power of $2$, or $k$ is odd with $k \leq  2 \phi(k)$ (for example, $k=q^n$ for some odd prime $q$), or $k=2\ell$ for an odd number $\ell$ with $\ell \leq 2 \phi(\ell)$ (for example, $\ell=q^n$ for some odd prime $q$), then $C_k(\Cc) \cap \T^{d+1} \neq \emptyset$. 
\end{proposition} 

This implies the second part of Theorem \ref{mainDyn} (that the action is not expansive for the indicated values of $k$) and finishes the proof of that theorem. 

\begin{proof} We denote by $\mathbf 1 = (1,\dots,1)$ the all-one vector. 
The proof depends on the observation that $P_k(\mathbf 1) = k+1>0$ whilst $P_k(-\mathbf 1)<0$ for the indicated values. Since $P_k$ is continuous on the  compact {and path-connected} set $\T^{d+1}$, the intermediate value theorem implies the existence of an intersection point in $C_k(\Cc) \cap \T^{d+1}$.  

If $k=2^n$, then $P_k = x_0 + x_1+1/x_1+\dots x_{k/2} + 1/x_{k/2}$, evaluating to $P_k(- \mathbf 1) = -k-1<0$. 

In general, $P_k$ is a sum of the form $P_k = x_0+x_1+\dots+x_d + R_k$, where $d=\phi(k)$ and $R_k$ consists of at most $k-d$ monomials with coefficient $1$ (so all these monomials evaluate to $\pm 1$ at $-\mathbf 1$. In particular, $P_k(- \mathbf 1) = -d-1 + R_k(-\mathbf 1) \leq -d-1+(k-d)$, which is negative if $k-d \leq d$, the stated condition. Notice that for $k=q^n$, $d=\phi(k)=q^n-q^{n-1}$, so the condition holds as soon as $q>2$. 

Similarly, if $k=2\ell$ with $\ell$ odd, then $P_k = x_0+x_1+1/x_1+\dots+x_d+1/x_d + \tilde R_k$ for $d=\phi(2\ell)=\phi(\ell)$ ($\ell$ odd) with $\tilde R_k$ consisting of at most $2 \ell - 2d$ monomials. By a similar reasoning as above, we find that $P_k(-\mathbf 1)<0$ if $2 \ell - 2 d \leq 2 d$, which gives the stated condition. 
\end{proof}

\begin{remark}
An early example of a non-expanding Bernoulli system of this type can be found in \cite{Ward-dots}, corresponding to the polynomial $1+x_0+x_1$ for $R_2$. It is not expansive, since $(e^{2 \pi i/3}, e^{4 \pi i/3})$ is a point on the intersection of the corresponding variety and the torus $\T^2$. This system has the same entropy $L'(\chi_{-3},-1) = 3 \sqrt{3}/(4\pi) L(\chi_{-3},2)$ as $C_2$, so is measurably isomorphic to it, but not topologically/algebraically conjugate to it: the associated prime ideals of the $R_2$-modules are not the same. 
\end{remark}

\begin{remark} \label{remprimorial} 
It is certainly possible that $\phi(k)/k<1/2$; there is a sequence of odd $k$ along which $\phi(k)/k$ is strictly decreasing to zero, for example by choosing $k=(\#n)/2$ with the primorial $\#n$ being product of the first $n$ prime numbers $p_1,\dots,p_n$) and $n>3$: 
$$\phi(\#n/2)/(\#n/2) = \prod_{r=2}^n (1 - 1/p_i) < \phi(\#4/2)/(\#4/2) = 48/105 < 1/2$$ and the product for $n \rightarrow + \infty$ tends to zero in the same way as the truncation of $1/(2\zeta(1)).$ In fact, the set of values of $\phi(k)/k$ is dense in the unit interval \cite[p.~31]{Ribenboim}. On the other hand, the average of $\phi(k)/k$ over $k \leq n$ tends to $6/\pi^2>1/2$ for $n \rightarrow + \infty$.  
\end{remark} 

\begin{remark} 
In terms of the vertices $A_k$ of the Newton polytope, $$P_k(-\mathbf 1) = \#\{ \mathbf a = (a_i) \in A_k \colon \sum a_i \in 2\Z\} - 1 - \#\{ \mathbf a = (a_i) \in A_k \colon \sum a_i \notin 2\Z\}.$$  It is not because $k$ violates the condition in the proposition that $P_k(-\mathbf 1)>0$; in fact, $P_{\#4}(-\mathbf 1)=-71.$ One may wonder whether $P_k(-\mathbf 1)$ is negative for all $k$. 
\end{remark} 

\begin{remark} 
For $k=2\ell$ with $\ell$ an odd prime, the following are explicit points $(x_0,\mathbf x)$ in $C_k(\Cc) \cap \T^{d+1}$. Set $\eta\coloneqq (1+i \sqrt{15})/4$, then 
\begin{enumerate}
\item If $\ell = 3 \mbox{ mod } 4$, $(x_0,\mathbf x)=(-1,-\eta, 1, -1, 1, -1, \dots, 1)$;
\item If $\ell = 1 \mbox{ mod } 4$, $(x_0,\mathbf x)=(1,\eta,\eta, 1, -1, 1, -1, \dots, -1)$.
\end{enumerate} 
Notice that for $k$ of this form, since $F_k$ is reciprocal, all its values at $\mathbf x \in \T^d$ are real, so then any point $(x_0,\mathbf x) \in C_k(\Cc) \cap \T^{d+1}$ forcedly has $x_0 = \pm 1$.  
\end{remark} 

\begin{remark} \label{rem:links} Algebraic dynamical systems such as the ones studied here have been used in the study of links, see \cite{SilverWilliams}. More precisely, one may ask for integer values of $\lambda$ and $k$ for which $\m(C_k(\lambda))$ equals the Mahler measure of the Alexander polynomial $A_L$ of a $d$-component link $L$ in $S_3$. If so, that Mahler measure is the logarithmic growth rate of the torsion homology in the sequence of finite abelian covers of $S^3$ branched over that link. To give two examples,
$$ \m(C_4(-1)) = \m(A_{L5a1}) \mbox{ and } \m(C_6(1)) = \m(A_{L7a6}), $$
where the links are given in the notation of the Thistlethwaite Link Table \cite{Links}. 
\end{remark}

\section{Random walks and asymptotics results for $\m(C_k)$: proof of Theorem \ref{mainRW}} 
\label{sec:random_walks_and_asymp_mCk}

In this section we mainly focus on proving an asymptotic result for $\m(C_k)$ if $k$ is of the form $k = 2q^r$, where $r \geq 1$ is fixed and $q$ is a prime tending to infinity. For this, we use a probabilistic interpretation: given random variables $X_1, \ldots, X_d$ on the complex unit circle, the Mahler measure of $C_k$ can be written as $\int_1^k \log(x) \rho_k(x) \, \d x$ for a probability density function $\rho_k$ corresponding to the random variable $|F_k(X_1,\dots,X_d)|$. We express the function $\rho_k(x)$ (through inversion of a two-sided Laplace transform) in terms of an exponential generating function for the constant terms of powers of $F_k$, and then identify the latter function with an infinite sum of Bessel functions (by direct comparison of Taylor coefficients), which leads to a formula in the style of Kluyver's for the random planar `flight' \cite{kluyver1906local}. Finally, we apply various asymptotic results about Bessel functions to find an asymptotic series for $\rho_k(x)$ in terms of increasing powers of $1/\sqrt{k}$, which will imply our result. After that, we consider the somewhat easier cases $k=2^r$ where $r$ tends to infinity and $k=q$ with $q$ a prime tending to infinity. First, we present some numerics, a result on small values of $\m(C_k)$, and a general upper bound for $\m(C_k)$.

\subsection{Small values of $\m(C_k)$, and an upper bound} In Table \ref{tab:mCk} we listed numerical approximations for the first few values of $\m(C_k)$.  For $k \neq 7$, these were computed as described in Section~\ref{MC6}, while for $k = 7$, we took the result from \cite[Table 2, case $n=8$]{Bailey}. For $k=9$, the method was inconclusive, and there is no result. Recall from Table \ref{cycth} that there are closed forms for $k \leq 4$ and a conjectured closed form for $k=5$, but no conjectured closed forms beyond this. 

 \begin{table}[ht!]
\centering
\begin{tabular}{ll|ll}
\toprule
$k$ 	& $\m(C_k)$ & $k$ 	& $\m(C_k)$  \\ \midrule
1 		& $0$	   & 6 		& $0.6439432099\dots$ \\ 
2 		& $0.3230659472\dots$   & 7		& $0.7668310881\dots$ \\
3 		& $0.4262783988\dots$	 & 8		& $0.7042270121\dots$ \\
4 		& $0.4839979734\dots$  & 9 		& $?$ \\
5 		& $0.6273170748\dots $ & 10 	& $0.7891727197 \dots$ \\
\bottomrule
\end{tabular}
\caption{The values $\m(C_k)$ for $k \leq 10, k \neq 9$ up to ten decimal places.}\label{tab:mCk}
\end{table}

We start with a result about small values of $\m(C_k)$ away from the `cyclotomic' case $k=1$. 

\begin{proposition} 
\label{infmck} 
For $k \geq 2$, the Mahler measure $\m(C_k)$ is bounded away from zero; 
more specifically, we have 
$$ 0.28 \leq \m(x^3-x-1) \leq  \inf_{k \geq 2} \m(C_k) \leq \m(C_2) \leq 0.33. $$
\end{proposition}
\begin{proof}
Recall that $\m(C_k) = \m(x_0 + F_k(x_1, \ldots, x_d))$. 
Define for any positive integers $k$ and $r$ the polynomial $$Q(x) = Q_{k,r}(x) = x^{N_{k,r}}(x + F_k(x^r, x^{r^2}, \ldots, x^{r^d})),$$ where $N_{k,r} \geq 0$ is minimal such that $Q(x) \in \Z[x]$.
Boyd \cite[pp.~118--119]{Boyd} proved that 
\begin{equation*}
\m(x_0 + F_k(x_1, \ldots, x_d)) \geq \lim_{r \to \infty} \m(Q_{k,r}(x)).
\end{equation*}
Smyth \cite{SmythBLMS} proved that any polynomial $P \in \Z[x]\setminus \{0\}$ with 
$$P(x)/P(1/x)  \notin \Z[x,x^{-1}]^\times = \{ \pm x^\ell \colon \ell \in \Z \}$$  has Mahler measure at least $\m(x^3-x-1)$.
Combining Boyd's and Smyth's results, it thus suffices to show that $Q(x) \neq \pm x^{\ell}Q(1/x)$ for any $k \geq 7$ (the result holds for $k \leq 7$ by inspection of Table \ref{tab:mCk}), $r \geq 3$ (or any sufficiently large $r$) and any $\ell$.

Assume to the contrary that $Q(x) = \pm x^{\ell}Q(1/x)$ for some such choice of $k$, $r$, and $\ell$. Note that $N_{k,r}$ is divisible by $r$ and that $Q$ has at least $d = \phi(k) \geq 3$ terms. Therefore at least two of the integers $c_i$ in the expansion
\[
Q(x) = x^{N_{k,r} + 1} + \sum_{i \geq 0} c_i x^{r i}
\]
are non-zero. Hence, comparing coefficients on both sides of the equality $\pm  x^{\ell} Q(1/x) = Q(x)$, we find that $\ell$ is divisible by $r$. The set of exponents modulo $r$ of $\pm x^{\ell} Q(1/x)$  is $\{0,-1\}$ and the set of exponents modulo $r$ of $Q(x)$ is $\{0,1\}$, which do not coincide since $r \geq 3$. This yields the required contradiction. 
Since the upper bound is immediate, this concludes the proof.
\end{proof}
In the other direction, we bound the Mahler measure of $C_k$ from above using the $L^2$-norm, as follows.
\begin{proposition}
	\label{prop:upper-bound}
For all $k$, $ \displaystyle{
		\m(C_k) \leq \frac{\log(k+1)}{2}  = \frac{\log{k}}{2} + \frac{1}{2k} + O\Big( \frac{1}{k^2} \Big).}$
\end{proposition}
\begin{proof}
Let $P = \sum_{\mathbf{n}} \mathbf{a}_{\mathbf{n}} \x^{\mathbf{n}}$ be a polynomial in $\ell$ variables. Then
\begin{equation} \label{l2} 
2 \, \m(P) = \int_{\T^{\ell}} \log(|P|^2) \, \d \mu \leq \log \int_{\T^{\ell}} |P|^2 \, \d \mu = \log \sum_{\mathbf{n}} \mathbf{a}_{\mathbf{n}}^2,
\end{equation}
where we used Jensen’s (concave) inequality and Parseval’s identity.

Apply this to the polynomial $P$ given by multiplying the Laurent polynomial $P_k = x_0 + F_k(\x)$ by its common denominator. Then $\m(C_k)=\m(P)$ and $P$ consists of $k+1$ monomials, all of which have coefficient $1$, so the upper bound for $2\, \m(C_k)$ in \eqref{l2} is $\log(k+1)$. 
\end{proof}

\subsection{Preliminaries using probability theory}
We first collect two general properties of Mahler measures of (reciprocal) Laurent polynomials using probability theory. 

\subsubsection*{Mahler measure and probability density functions} 
The first concerns rewriting Mahler measure in terms of a probability density function. Let $H \in \Cc[x_1^{\pm 1}, \ldots, x_{d}^{\pm 1}]$ be a non-zero Laurent polynomial and $X_1, X_2, \ldots, X_d$ random variables on the complex unit circle $\T$. Then $|H(X_1, \ldots, X_d)|$ is a random variable on the torus $\T^d$, taking values in $\R_{\geq 0}$. The corresponding \emph{probability density function}, when it exists, is an integrable function $\rho_H \colon \R \rightarrow \R$ (determined up to measure zero) for which the distribution function satisfies 
\begin{equation} \label{defrhogen} \int\limits_{\substack{ \mathbf x \in \T^d \\ 0 \leq |H(\mathbf x)| \leq b}} \hspace*{-5mm} \d \mu= \int_0^b \rho_H(t) \, \d t; \end{equation} 
then for every function $g \colon \R \rightarrow \Cc$ such that $g(|H(X_1,\dots,X_d)|)$ admits a probability density function (in particular, if $g \rho_H$ is integrable), we find 
\begin{equation} 
	\label{pdfprope} 
	\int\limits_{\mathbf x \in \T^d} g(|H(\mathbf x)|) \, \d \mu= \int_0^\infty g(t) \rho_H(t) \, \d t.  
\end{equation} 
We first discuss the existence of $\rho_H$.

\begin{proposition} 
	For any non-zero Laurent polynomial $H \in \Cc[x_1^{\pm 1}, \ldots, x_{d}^{\pm 1}]$ an integrable probability density function $\rho_H$ exists if and only if $|H(\mathbf x)|$ is not constant on $\T^d$. 
\end{proposition} 

\begin{proof} 
First assume $|H(\mathbf x)|$ is not constant on $\T^d$.
Let $\lambda_d$ denote the Lebesgue measure on $\R^d$. Define the measure $\nu$ on the $\sigma$-algebra of Lebesgue measurable subsets $A$ of $\R$ by 
$$ \nu(A) \coloneqq \int\limits_{\substack{ \mathbf x \in \T^N \\|H(\mathbf x)| \in A}} \hspace*{-5mm} \, \d \mu \hspace*{5mm} = \hspace*{-5mm} \int\limits_{\substack{ \mathbf t \in U \\|H(\exp(2 \pi i \mathbf t))| \in A}} \hspace*{-5mm}  \d \mathbf t $$
where $U \coloneqq (0,1)^d$. By the Radon--Nikodym theorem, the existence of $\rho_H$ is equivalent to the fact that the measure $\nu$ is absolutely continuous with respect to the Lebesgue measure $\lambda=\lambda_1$ on $\R$, i.e., for all $A$, $\lambda(A)=0$ implies $\nu(A)=0$ (see, e.g., \cite[Thm.\ 31.8]{Bil}). 

So assume $\lambda(A)=0$. Notice that $\nu(A) = \nu(A \cap \R_{\geq 0})$, so we may assume that $A \subseteq \R_{\geq 0}$. Let $A^2$ denote the set of squares of elements of $A$; it follows that $$\lambda(A^2) = 2 \int_{A} y \d y =0$$
since any integrable function integrates to zero over a set of measure zero. 

For any non-zero Laurent polynomial $H$, let $\tilde H \coloneqq H(1/x_1,\dots,1/x_d)$, and set 
$$ \varphi(\mathbf t) \coloneqq (H \tilde H)(\exp(2 \pi i \mathbf t)) = |H(\exp(2 \pi i \mathbf t))|^2. $$ Now $\varphi \colon U \rightarrow \R$ is a smooth map, and, by unfolding the definitions,  
$$ \nu(A) = \lambda_d(\varphi^{-1}(A^2)). $$
We will show that $\varphi$ has the \emph{zero property} in the sense of \cite{Pon}, i.e., the pullback of any set of $\lambda$-measure zero has $\lambda_d$-measure zero. This implies the desired result, since we have $\lambda(A^2)=0$. To show the zero property, we use that $\varphi$ is smooth and apply \cite[Theorem 1]{Pon}, stating that the property is equivalent to the Jacobian $$J \colon U \rightarrow \R^d, \mathbf t \mapsto (\partial_{t_1} \varphi(\mathbf t),\dots,\partial_{t_d} \varphi(\mathbf t))$$ of $\varphi$ being non-vanishing $\lambda$-almost everywhere. We now use the assumption that $|H|$ is not constant on the torus, so that $H \tilde H$ is not constant, and hence there exists $j$ with $\partial_{t_j} \varphi \not \equiv 0$. Since $\partial_{t_j} \varphi$, as a finite sum of exponential functions, is real-analytic on $U$, its zero locus has $\lambda_d$-measure zero (see, e.g., \cite{Mit}), and the result follows. 

On the other hand, if $|H|$ is constant on $\T^d$, say $|H| = w \in \R$, then an integrable probability density function $\rho_H$ does not exist. Indeed, in this case, the left hand side of \eqref{defrhogen} equals the indicator function of $\{ b \geq w\}$, i.e., the Heaviside function $H(b-w)$, and a solution exists only in the sense of distributions, namely, $\rho_H(t) = \delta(t-w)$. 
\end{proof}

\begin{remark} \label{remzetaM}
The probability density function is also, at least in the sense of distributions, the inverse Mellin transform of $Z(H,s-1)$, where $Z(H,s)$ is the \emph{zeta Mahler function} 
 \[
 Z(H;s) \coloneqq  \int_{\T^d} |H(x_1,\dots,x_d)|^s \, \d\mu \ \ (=\int_0^\infty x^{s} \rho_{H}(x) \, \d x); 
 \]
  which is convergent on the right-hand plane $\Re (s) > \mu_0$ for some negative real number $\mu_0$, see \cite[Theorem 8]{Akatsuka}. 
For the existence of the inverse Mellin transform, see, e.g.,  \cite[11.1.2.3]{Mellinex}. For example, the `pathological' case where $|H| \equiv w \in \R$ on $\T^d$ gives zeta Mahler function $Z(H;s) \equiv w$, and this indeed has inverse Mellin transform $\delta(t-w)$. 

  Also notice the identity 
  \begin{equation}  \label{MMandMahlerZeta} 
  \m(H) = \frac{\d}{\d s} Z(H;s) |_{s=0}. 
  \end{equation} 
  for the Mahler measure itself. 
\end{remark} 

The above discussion allows us to rewrite the Mahler measure in the desired form. 
\begin{lemma} For non-zero $H \in \Cc[x_1^{\pm 1}, \ldots, x_{d}^{\pm 1}]$ for which the associated probability density function $\rho_H$ exists, let $L \coloneqq  \max\limits_{\mathbf x \in \T^d} |H(\mathbf x)|$. Then  
\begin{equation}
\label{Jensen}
\m(x_0 + H) = \int_{1}^{L} \log(x) \rho_{H}(x) \, \d x.    
\end{equation}
\end{lemma} 

\begin{proof} 
By the compactness of $\T^d$, the image $|H(\T^d)|$ is contained in the finite interval $[0,L]$, and hence the same holds for the support of $\rho_H$. We have the sequence of equalities
\begin{equation*}
\m(x_0 + H) \overset{(1)}{=}  \int_{\T^d} \log^+{|H(x_1, \ldots, x_d)|}\,  \d \mu \overset{(2)}{=} \int_{0}^{L} \log^+(x) \rho_{H}(x) \, \d x \overset{(3)}{=} \int_{1}^{L} \log(x) \rho_{H}(x) \, \d x, 
\end{equation*}
where $(1)$ is Jensen's formula; $(2)$ is property \eqref{pdfprope} applied to $g=\log^+$ (which is integrable over $[0,L]$, to $1/L-1$), and $(3)$ changes the integration bounds and, correspondingly, $\log^+$ to $\log$.  
\end{proof} 

\subsubsection*{Probability density functions and constant term sequence} The second general result concerns \emph{reciprocal} Laurent polynomials $H$, meaning that $H$ is invariant under all transformations $x_j \mapsto 1/x_j$ for $j=1,\dots,d$. 
Denote by $\text{CT}[g]$ the \emph{constant term} of a Laurent polynomial $g$, and define the exponential generating series
\begin{equation} \label{defeh}
E_H(t) \coloneqq  \sum_{m \geq 0} \frac{\text{CT}[H^m]}{m!} t^m.
\end{equation} 
Since the $\text{CT}[H^m]$ are polynomially bounded in $m$, this defines an entire function in $t$. 
The relation between the density $\rho_H$ and $E_H$ is given in the following lemma, which states that, for reciprocal $H$, the probability density $\rho_H$ is the inverse Laplace transform of the symmetrized exponential generating function $E_H(t) + E_H(-t)$.
\begin{lemma}
  \label{LemmaLaplace}
  For any non-constant reciprocal Laurent polynomial $H \in \R[x_1^{\pm 1}, \ldots, x_d^{\pm 1}]$ 
  and any $x \in [0,L]$, we have 
  \begin{equation} \label{EHEHinv} \rho_H(x) = \frac{1}{2 \pi} \int_{-\infty}^{\infty} (E_H(is) + E_H(-is)) e^{isx} \, \d s. \end{equation} 
  
\end{lemma}
\begin{proof} 
We claim that 
\begin{equation}
\label{ConstantTermAsIntegral}
\text{CT}[H^{2m}] \overset{(a)}{=} \int_{\T^d} H(x_1, \ldots, x_N)^{2m} \, \d \mu \overset{(b)}{=} \int_0^\infty x^{2m} \rho_{H}(x) \, \d x,   
  \end{equation}
  i.e., constant terms of even powers are the corresponding probability theoretic moments. 
  
The equality $(a)$ is a general principle to compute the constant term of a power of a polynomial, noticing that $\int_{|x|=1} x^n \, \d x$ is zero unless $n=0$, when it is $1$. For $(b)$, notice that, since on the torus $\T^d$ the transformations $x_j \mapsto 1/x_j$ describe complex conjugation ($x_j x_j^{-1} = 1 = |x_j| = x_j \overline x_j$), a real reciprocal Laurent polynomial $H$ is real-valued on the torus. Hence for any integer $m$, $H^{2m} = (H \overline H)^m = |H|^{2m}$, and $(b)$ follows by applying \eqref{pdfprope} to the function $g(x)=x^{2m}$.  

Next, we claim that for any $t \in \Cc$,
  \begin{equation} \label{EHEH} E_H(t) + E_H(-t) = \int_{-\infty}^\infty \rho_{H}(|x|) e^{-t x} \, \d x, \end{equation} 
To prove this, we expand the exponential in the integral in \eqref{EHEH} and use \eqref{ConstantTermAsIntegral}, as follows. 
\begin{align*}
\int_{-\infty}^\infty \rho_{H}(|x|) e^{-t x} \, \d x &= \int_0^\infty \rho_{H}(x) \left( e^{t x} + e^{- t x} \right) \, \d x \\ &= \sum_{m \geq 0} \frac{t^m}{m!}  \int_0^\infty \rho_H(x) \left(x^m + (-x)^m \right) \, \d x\\
&=  2\sum_{m \geq 0}\frac{\text{CT}[H^{2m}]}{(2m)!}t^{2m} = E_H(t) + E_H(-t).
\end{align*}
Then equation \eqref{EHEHinv} follows by applying the  inverse two-sided Laplace transform  \cite[Section 4.7.2]{andrews1999integral} 
to \eqref{EHEH} (notice that we apply the result for $x>0$, so $|x|=x$).  
\end{proof}

\subsection{Main result and $\m(C_k)$.} 

In this subsection, we let $k = 2q^r$ for $r \geq 1$ and an odd prime $q$. The goal is to establish the asymptotics for $\m(C_k)$ as $q$ tends to infinity and $r$ is fixed. We first show how to reduce this to an asymptotic property of the corresponding probability density function. 

\begin{notation} 
Let $\rho_k \coloneqq  \rho_{F_k}$, with support contained in  $[0,k]$.
\end{notation} 

The asymptotics in question is the following result. 
\begin{theorem}
\label{MainTHM6}
Let $r \geq 1$ be fixed and suppose $q$ is prime. If $k = 2q^r$, then for sufficiently large $q$, we have
\begin{equation} \label{rholow}
\rho_k(x) = \sqrt{\frac{2}{\pi k}} e^{-x^2/(2 k)} + O(k^{-3/2})
\end{equation} 
where the implied constant is independent of $x$.
\end{theorem}

The result implies that the density $\rho_k$ converges to the density function of a \emph{half-normal distribution} (i.e., the absolute value of the normal distribution) in the large-$q$ limit when $k$ goes through integers of the form $2q^r$ as in the statement of Theorem~\ref{MainTHM6}. The proof of Theorem~\ref{MainTHM6} will occupy the next few subsections. There, we will also indicate how to get higher order terms in the asymptotic expansion \eqref{rholow} and in Theorem \ref{mainRW}. Before we start the proof, we show how the lower order terms in Theorem \ref{mainRW} already follow from Theorem \ref{MainTHM6}.
\begin{proof}[Proof of Theorem \ref{mainRW} assuming Theorem \ref{MainTHM6} for $k=2q^r$, up to $O(1/\sqrt{k})$]

 By \eqref{Jensen}, we have 
 \begin{equation}
\label{JensenCk}
\m(C_k) = \m(x_0 + F_{k}(\mathbf x)) = \int_{1}^{k} \log(x) \rho_{k}(x) \, \d x.    
\end{equation}
Thus, we may rewrite 
\begin{align}
\label{eq:Thm G from rho_k asymptotic}
\m(C_k) = \int_{1}^{\infty} \log(x) \sqrt{\frac{2}{\pi k}} e^{-x^2/(2 k)}  \, \d x - \int_{k}^{\infty} \log(x) & \sqrt{\frac{2}{\pi k}} e^{-x^2/(2 k)}  \, \d x \nonumber \\ &+ O(k^{-3/2}) \int_{1}^{k} \log(x) \, \d x.
\end{align}
In \eqref{eq:Thm G from rho_k asymptotic}, the second integral decays exponentially in $k$ and the third term is $O(\log(k)/\sqrt{k})$. The first integral differs by $O(k^{-1/2})$ from
\begin{equation*}
\int_{0}^{\infty} \log(x) \sqrt{\frac{2}{\pi k}} e^{-x^2/(2 k)}  \, \d x  =  \frac{1}{2}\log(k) - \frac{1}{2}\log(2) - \frac{1}{2}\gamma,
\end{equation*}
where we computed the last integral using \texttt{Mathematica}.
\end{proof}

In the next few subsections, we prepare for the proof of Theorem \ref{MainTHM6} with $k=2q^r$. 

\subsection{From exponential generating functions to probability densities}

Since $F_k(\x)$ is reciprocal if $k$ is even, so in particular when $k = 2q^r$, we will apply Lemma \ref{LemmaLaplace} to the densities $\rho_k$. Write $E_k \coloneqq  E_{F_{k}}$. We first express $E_k$ in terms of modified Bessel functions of the first kind, defined by 
\[
I_\alpha(2x) = \sum_{m \geq 0} \frac{1}{m! \Gamma(m + \alpha + 1)} x^{2m + \alpha}. 
\]

\begin{lemma}
	\label{besselidentity}
	Let $k = 2 q^r$, where $q$ is an odd prime and $r$ a positive integer. For any $t \in \Cc$, we have 
	\[
	E_k(t) \coloneqq  E_{F_k}(t) = \Big(I_0(2t)^q + 2\sum_{j \geq 1} I_{j}(2t)^q\Big)^{q^{r-1}}.\]
\end{lemma}
\begin{proof}
	First assume $r = 1$. We will compute the exponential generating function $E_k$. In this case $F_k(\x)$ is of the form 
	\[
	F_k(\x) = x_1 + \frac{1}{x_1} + \cdots + x_{q - 1} + \frac{1}{x_{q - 1}} + \frac{x_1}{x_2} \frac{x_3}{x_4} \cdots \frac{x_{q - 2}}{x_{q - 1}} + \frac{x_2}{x_1} \frac{x_4}{x_3} \cdots \frac{x_{q - 1}}{x_{q - 2}}.
	\]
	For the $m$-th power we have
	\begin{align*}
		F_k(\x)^m = \sum_{m_1 + \cdots + m_q = m} {m \choose m_{1},m_{2},\ldots ,m_{q}}
		\left(x_1 + \frac{1}{x_{1}}\right)^{m_1} &\cdots \left(x_{q - 1} + \frac{1}{x_{q - 1}}\right)^{m_{q - 1}} \times \\ &\times \left(\frac{x_1}{x_2} \frac{x_3}{x_4} \cdots \frac{x_{q - 2}}{x_{q - 1}} + \frac{x_2}{x_1} \frac{x_4}{x_3} \cdots \frac{x_{q - 1}}{x_{q - 2}} \right)^{m_q}.
	\end{align*}
	Note that $(a+1/a)^m = \sum\limits_{j=0}^m {m \choose j} a^{2j-m}$. 	
	First suppose $m$ is \emph{even}. Taking the constant term yields
	\begin{align*}
		\text{CT}[F_k(\textbf{x})^m] = \hspace*{-7mm} \sum_{\substack{m_1 + \cdots + m_q = m\\ 0 \leq j_q \leq m_q}}   {m \choose m_{1},m_{2},\ldots ,m_{q}} &\binom{m_1}{m_1/2 + m_q/2 - j_q}   \cdots \binom{m_{q -1}}{m_{q-1}/2 - m_q/2 + j_q} \binom{m_q}{j_q},
	\end{align*}
	where we use the convention $\binom{a}{b/2} = 0$ for integers $a$ and \emph{odd} integers $b$. Since $q$ is odd and $m$ is even, we can assume all $m_1, \ldots, m_{q-1}$ are even. 
	Define $h_m(v)$ by 
	\begin{align*}
		h_m(v) &\coloneqq \hspace*{-7mm}  \sum_{\substack{m_1 + \cdots + m_q = m\\ j_q = m_q/2 - v}} {m \choose m_{1},m_{2},\ldots ,m_{q}} \binom{m_1}{m_1/2 + m_q/2 - j_q} \cdots \binom{m_{q -1}}{m_{q-1}/2 - m_q/2 + j_q} \binom{m_q}{j_q} \\
		& \hspace{1mm} =  \hspace*{-7mm} \sum_{m_1 + \cdots + m_q = m} {m \choose m_{1},m_{2},\ldots ,m_{q}}  \binom{m_1}{m_1/2  + v} \cdots \binom{m_{q -1}}{m_{q-1}/2 + v} \binom{m_q}{m_q/2 + v},
	\end{align*}
	and notice that $h_m(v)$ is an even function of $v$. Therefore 
	\[
	\text{CT}[F_k(\textbf{x})^m] = \sum_{\substack{v \in \Z\\|v| \leq \frac{m}{2q}}} h_m(v) = h_m(0) + 2 \sum_{\substack{v \in \Z_{>0}\\1 \leq v  \leq \frac{m}{2q}}}  h_m(v).
	\]
Finally, we observe
	\[
	\frac{h_m(v)}{m!} = \sum_{\widetilde{m_1} + \cdots + \widetilde{m_q} = \frac{m}{2} - q v } \frac{1}{\widetilde{m_1}! \Gamma(\widetilde{m_1} + 2v + 1)} \cdots \frac{1}{\widetilde{m_q}!  \Gamma(\widetilde{m_q} + 2v + 1)}
	\]
	and this is precisely the coefficient of $x^m$ in the series expansion of $I_{2v}(2x)^q$. 
	
	The proof in the second case, where $m$ is \emph{odd}, is entirely analogous. This proves the statement for $r = 1$. If $r > 1$, the result follows directly by combining our formula for $r=1$ with \cite[Proposition~2.3]{Lala} (the integers denoted $U_k(m)$ in that reference coincide with the constant terms $\text{CT}[F_k(\textbf{x})^m]$.)
\end{proof}

\begin{proposition} \label{prop:gtrepr}
For $k=2q^r$ with $q$ an odd prime number, we have 
\begin{equation*} 
	\rho_k(x) =  \frac{1}{\pi} \Re \int_0^\infty \Bigg( \Big(I_0(2 i s)^q + 2\sum_{j \geq 1} I_{j}(2 i s)^q\Big)^{q^{r-1}} \hspace*{-5mm} +  \Big(I_0(-2 i s)^q + 2\sum_{j \geq 1} I_{j}(-2 i s)^q\Big)^{q^{r-1}}\Bigg) e^{isx} \, \d s. 
\end{equation*}    
\end{proposition} 

\begin{proof} Notice that for such values of $k$, the polynomial $F_k$ is reciprocal. Hence we can apply \eqref{EHEHinv} 
to find $$\rho_k(x) = \frac{1}{2 \pi} \int_{-\infty}^{\infty} (E_k(is) + E_{k}(-is)) e^{isx} \, \d s.$$ All that is left is to plug in the representation from Lemma 
\ref{besselidentity}, and the result follows. 
\end{proof} 

\begin{remark} \label{remkluif} 
Kluyver  \cite{kluyver1906local} gave a formula for the probability density function of the `linear' walk $|X_1 + \cdots + X_N|$ (no dependencies) in terms of integrals of Bessel functions. In our setting, a new feature is that an infinite sum of Bessel functions appears. 
\end{remark}

\begin{remark}
Despite being an explicit formula for the density $\rho_k$,  Proposition \ref{prop:gtrepr}  is of no use to compute values of $\rho_k(x)$ to high precision for a fixed $k$ and $x$, the main bottleneck being that we do not know the asymptotics for the integrand for large $s$. We will return to this computational question in the next section.  
\end{remark}

\subsection{Large $q$ asymptotics of $\rho_k$} We use the integral representation from Proposition \ref{prop:gtrepr} to study the large $q$ asymptotics for $\rho_k$ with $r$ fixed. It will turn out that the integral for $0 \leq s \leq 1$ of the term $I_0(2 i s)^{q^r}$ in the integrand is dominant for large $q$.

\begin{lemma}
	\label{Lem:725}
	As $q \rightarrow +\infty$, we have  
	\[
	\rho_{k}(x) = \frac{2}{\pi} \Re \int_{0}^{\infty} I_{0}(2 i s)^{q^r} e^{i s x} \, \d s + O(\kappa^{-q})
	\]
	for some $\kappa>1$, where the implied constant is independent of $x$.
\end{lemma}
\begin{proof}
	First of all, the Bessel function $I_j$ is an odd function for odd $j$ and an even function for even $j$. Using this and the binomial theorem, we find
	\begin{align*}
		&\rho_{k}(x) - \frac{2}{\pi} \Re \int_{0}^{\infty} I_{0}(2 i s)^{q^r} e^{i s x} \, \d s \\&= 
		\frac{1}{\pi} \Re \bigintssss_{\, 0}^{\infty} \sum_{w = 1}^{q^{r-1}} \binom{q^{r-1}}{w} \Bigg(\Big(2\sum_{j \geq 1} I_{j}(2 i s)^q \Big)^w + \Big(2\sum_{j \geq 1} (-1)^j I_{j}(2 i s)^q \Big)^w\Bigg) I_{0}(2 i s)^{q^r - qw} e^{i s x} \, \d s.
	\end{align*}
	We find that
	\begin{align*}
		&\left| \rho_{k}(x) - \frac{2}{\pi} \Re \int_{0}^{\infty} I_{0}(2 i s)^{q^r} e^{i s x} \, \d s \right| \\& \leq  
		\frac{1}{\pi} \sum_{w = 1}^{q^{r-1}} 2^w \binom{q^{r-1}}{w}  \bigintsss_{\, 0}^{\infty} \Bigg( \Big|\sum_{j \geq 1} I_{j}(2 i s)^q \Big| ^w + \Big|\sum_{j \geq 1} (-1)^j I_{j}(2 i s)^q \Big|^w\Bigg)  \, \d s
	\end{align*}
	using the bound $|I_0(2 i s)| \leq 1$ (see \cite[10.14.1]{NIST}) and the triangle inequality. 
	By the triangle inequality for the $L^{w}$-norm we find
	\begin{equation}
		\label{eq:ineqBess}
		\int_{0}^{\infty} \left( \left|\sum_{j \geq 1} I_{j}(2 i s)^q \right| ^w + \left|\sum_{j \geq 1} (-1)^j I_{j}(2 i s)^q \right|^w\right)  \, \d s \leq 2 \left( \sum_{j \geq 1} \left(\int_{0}^\infty |I_{j}(2 i s)|^{qw} \, \d s \right)^{1/w} \right)^w.
	\end{equation} 
	For the Bessel function we have the bounds
	\[
	|I_{j}(2 i s)| \leq \min \{ b (j)^{-1/3}, c |2s|^{-1/3}\}
	\]
	for $j > 0$ and real numbers $b = 0.674885\dots$ and $c = 0.7857468704 \dots$; this is the main result in \cite{Landau}. Using these bounds, we find that
	\begin{align*}
		\int_0^\infty |I_{j}(2 i s)|^{qw} \, \d s= \left( \int_0^{j} + \int_{j}^{\infty} \right) |I_{j}(2 i s)|^{qw} \, \d s &\leq b^{qw} j^{1 -qw/3} + c^{qw} \int_{j}^\infty  (2s)^{-qw/3} \, \d s \\ &= b^{qw} j^{1 -qw/3} + \frac{c^{qw} 2^{-qw/3} j^{1- qw/3}}{qw/3 - 1}\\
		& \leq 2 j^{1-qw/3} b^{qw}
	\end{align*}
	if $q \geq 6$; hence 
	\[
	\left( \int_0^\infty |I_{j}(2 i s)|^{qw} \, \d s \right)^{1/w} \leq 2^{1/w} j^{1/w-q/3} b^{q} \leq 2 j^{1-q/3} b^{q}, 
	\]
	which we can sum, using the Riemann zeta function $\zeta(s)$, to obtain
	\[
	\sum_{j \geq 1} \left(\int_0^{\infty} |I_{j}(2 i s)|^{qw} \, \d s\right)^{1/w} \leq 2 \zeta(q/3 - 1) b^q \leq 2 \frac{3}{2} b^q = 3 b^q
	\]
	if $q \geq 11$. Combining this with \eqref{eq:ineqBess}, we find 
	\begin{align}
		\label{differenceRHOk}
		\left| \rho_{k}(x) - \frac{2}{\pi} \Re \int_{0}^{\infty} I_{0}(2 i s)^{q^r} e^{i s x} \, \d s \right| \leq
		\frac{1}{\pi} 2 \sum_{v = 1}^{q^{r-1}}\binom{q^{r-1}}{v} (6 b^q)^v =\frac{2}{\pi} \left( (1 + 6 b^q)^{q^{r-1}} - 1 \right).
	\end{align}
	Finally, the inequalities $1+y \leq e^y \leq 2y + 1$ for $0 \leq y \leq 1$ give
	\[ (1 + 6 b^q)^{q^{r-1}} - 1 = e^{q^{r-1} \log(1 + 6 b^q )} - 1 \leq  e^{6 q^{r-1}  b^q} - 1 \leq 12 q^{r-1}  b^q + 1 - 1 = 12 q^{r-1}  b^q\]
	for $q$ sufficiently large (depending on $r$). The result follows. 
\end{proof}

Next, we show that we can also discard the remaining integral over $s \geq 1$. 

\begin{lemma}
	\label{lem:726}
	The integral
	\[
	\Re \int_{1}^{\infty} I_{0}(2 i s)^{q^r} e^{i s x} \, \d s
	\]
	is exponentially small in $q$, independent of $x$.
\end{lemma}
\begin{proof}
	Again using the main estimates in \cite{Landau}, we have
	\[
	|I_{0}(2 i s)| \leq \frac{1}{(s \pi)^{1/2}}
	\]
	for $s > 0$.
	This implies
	\[
	\left| \int_{1}^{\infty} I_{0}(2 i s)^{q^r} e^{i s x} \, \d s \right| \leq \pi^{-q^r/2} \int_1^\infty s^{-q^r/2} \, \d s = \pi^{-q^r/2} \frac{1}{q/2 - 1} \leq  \pi^{-q^r/2}
	\]
	if $q \geq 5$, proving the statement.
\end{proof}

For the remaining integral, we first study the asymptotics of the integrand, suitably rescaled. 
\begin{lemma}
	\label{LemI0exp}
	For $s \in [0, \sqrt{v}]$ and $M \geq 1$, 
	\begin{align*}
		I_0(2 i s/\sqrt{v})^{v} =\sum_{j = 0}^{M-1} e^{-s^2} \cdot \frac{a_j(s)}{v^j} + A_M(s,v),
	\end{align*}
	where $a_j(s)$ is a polynomial in $s$ and 
	\[ |A_M(s,v)| \leq e^{-s^2} \cdot\frac{W_M(s)}{v^{M}},\]
	for a polynomial $W_M$ only depending on $M$.
\end{lemma}
\begin{proof}
	We claim that $I_0(2 i y) \leq e^{-y^2}$ in the interval $[0, 1]$. Since $I_0(0) = 1$, it suffices to shows that the product $I_0(2 i y)  e^{y^2}$ is a decreasing function. For this it is enough to check that its derivative $$-2 e^{y^2} y I_2(2 i y)$$ is a non-positive function, which is easy to see as  $I_2(2 i y)$ has no zeros in this region \cite[p. 409]{Abramowitz}.
	Moreover, since $I_0(2 i y)$ has no zeros on $[0,1]$ (loc.\ cit.), $\log(I_0(2 i y))$ is a smooth function on this interval. Its Taylor series is of the form 
	\begin{equation}
		\label{LogI0Taylor}
		\log(I_0(2 i y)) = -y^2 - \frac{y^4}{4} -\frac{y^6}{9}-\frac{11 y^8}{192} + \cdots + O(y^{2L}) = -y^2 + p_{L}(y) + O(y^{2L})    
	\end{equation}
	where $p_L(y)$ is a polynomial of degree at most $2L-2$. 
	By Taylor's theorem and the fact that $\log(I_0(2 i y))$ is smooth, there exists a constant $r_L$ only depending on $L$ such that
	\begin{equation}
		\label{LogI0TaylorError}
		\left| \log(I_0(2 i y)) + y^2 -p_L(y) \right| \leq r_L y^{2L}.
	\end{equation}
	Applying Taylor's theorem to the exponential function, we have for $z \leq 0$
	\begin{equation}
		\label{exponent_taylor}
		\left|e^{z} - \sum_{m = 0}^{M-1} \frac{z^m}{m!} \right| \leq |z|^{M}.    
	\end{equation}
	We apply \eqref{exponent_taylor} to $z = \log(I_0(2 i y)) + y^2$, resulting in	
	\begin{align}
		I_0(2 i y)^v &= e^{-v y^2} e^{v(\log(I_0(2 i y)) + y^2)}\\
		&=e^{-vy^2}\sum_{m = 0}^{M-1}\frac{1}{m!} v^m (\log(I_0(2 i y)) + y^2)^m + e^{-v y^2}\epsilon_{M}(y), 
	\end{align}
	where, using the inequalities \eqref{LogI0TaylorError} and \eqref{exponent_taylor}, we have 
	\[
	|\epsilon_M(y)| \leq  v^{k} |\log(I_0(2 i y)) + y^2|^{M} \leq v^{M} r_4^{M} y^{4 M}.
	\]
	By \eqref{LogI0TaylorError} we can write
	\[
	\left| (\log(I_0(2 i y)) + y^2)^m - p_{2M}(y)^m \right| \leq u_M y^{4M}
	\]
	for some constant $u_M$, independent of $y$.
	Combining all this gives
	\begin{align}
		I_0(2 i y)^v     &=e^{-vy^2}\sum_{m = 0}^{M-1}\frac{1}{m!} v^m p_{2k}(y)^m + e^{-v y^2}\epsilon'_{M}(y),
	\end{align}
	where
	\begin{equation}
		|\epsilon'_{M}(y)| \leq w_M v^M y^{4M}    
	\end{equation}
	for some constant $w_M$.  We now substitute $y = s/\sqrt{v}$ for $s \in [0, \sqrt{v}]$, to show that we have the following expression for $I_0(2 i s/\sqrt{v})^v$ in terms of polynomials $a_j(s)$ in $s$:
	\begin{align}
		I_0(2 i s/\sqrt{v})^v 
		&=e^{-s^2}\sum_{m = 0}^{M-1}\frac{1}{m!} v^m p_{2M}(s/\sqrt{v})^m + e^{-s^2}\epsilon'_{M}(s/\sqrt{v})\\
		&=e^{-s^2} \sum_{j = 0}^{(M-1)(2M - 3)} \frac{1}{v^j} a_j(s)   
		+ e^{-s^2}\epsilon'_{M}(s/\sqrt{v})\\
		&= e^{-s^2} \sum_{j = 0}^{M-1} \frac{1}{v^j} a_j(s) + e^{-s^2} \sum_{j = M}^{(M-1)(2M - 3)} \frac{1}{v^j} a_j(s)  
		+ e^{-s^2}\epsilon'_{M}(s/\sqrt{v}).
	\end{align}
	
	Next we bound the terms $1/v^j$ for $j \geq M$ trivially by $1/v^M$. It follows that
	\[
	\left|e^{-s^2} \sum_{j = M}^{(M-1)(2M - 3)} \frac{1}{v^j} a_j(s)  
	+ e^{-s^2}\epsilon'_{M}(s/\sqrt{v}) \right| \leq \frac{1}{v^M} W_M(s)
	\]
	for a polynomial $W_M$ depending only on $M$. This concludes the proof.
\end{proof}
\begin{remark}
The proof of Lemma \ref{LemI0exp} provides an explicit method to compute the polynomials $a_j(s)$. For instance
\begin{align*}
	a_0(s) &= 1; \quad
	a_1(s) = -\frac{s^4}{4}; \quad
	a_2(s) = \frac{s^6}{288} \left(9 s^2-32\right); \quad
	a_3(s) = -\frac{s^8}{1152}  \left(3 s^4-32 s^2+66\right);\\
	a_4(s) &= -\frac{s^{10}}{4147200}\left(675 s^6-14400 s^4+85000 s^2-131328\right);
	\\a_5(s) &= -\frac{s^{12}}{16588800}\left(-135 s^8+4800 s^6-55300 s^4+236928 s^2-302720\right). 
\end{align*}
These polynomials have only real roots, and one wonders whether this is true for all $a_j(s)$.    
\end{remark}

\begin{proof}[Proof of Theorem \ref{MainTHM6} for $k=2q^r$]
From Lemmas \ref{Lem:725} and \ref{lem:726}, we have 
\begin{equation}
	\label{RHOKfinal}
	\rho_k(x) \sim \frac{2}{\pi} \Re \int_{0}^{1} I_{0}(2 i s)^{q^r} e^{i s x} \, \d s = \frac{2 \sqrt{2}}{\pi \sqrt{k} }\Re \int_{0}^{\sqrt{k/2}} I_{0}(2 i s/\sqrt{k/2})^{k/2} e^{i s x/ \sqrt{k/2}} \, \d s     
\end{equation}
	to all orders in $q$ as $q \to \infty$. Applying Lemma \ref{LemI0exp} with $M = 1$ gives 
	\begin{equation} \label{sumsum} 
	I_{0}(2 i s/\sqrt{k/2})^{k/2} = e^{-s^2} + A_1(s,k/2) 
	\end{equation} 
	and integrating the second term yields
	\begin{equation} \label{refber} \frac{1}{\sqrt{k/2}}\left|\Re \int_0^{\sqrt{k/2}} A_1(s,k/2) e^{i s x/\sqrt{k/2}} \, \d s \right| \leq  \frac{2\sqrt{2}}{k^{3/2}}
	\int_0^\infty |W_1(s)| e^{-s^2} \, \d s
	\end{equation} 
	for some polynomial $W_1(s)$ in $s$. The right hand side in \eqref{refber} is $O(1/k^{3/2})$, since the integral occurring in it is bounded: for all $j \geq 0$, we have 
	\[
	\int_0^\infty s^j e^{-s^2} \, \d s = \frac{1}{2} \Gamma \left( \frac{1 + j}{2} \right). 
	\]
	To integrate the first term in \eqref{sumsum}, we notice that 
	\[
	\int_0^{\sqrt{k/2}} e^{-s^2} e^{i s x/ \sqrt{k/2}} \, \d s \sim \int_0^{\infty} e^{-s^2} e^{i s x/ \sqrt{k/2}} \, \d s
	\]
	to all orders in $q$, 
	and
	\[
	\Re \int_0^{\infty} e^{-s^2} e^{i s x/ \sqrt{k/2}} \, \d s =
	\int_0^{\infty} e^{-s^2} \cos(s x/ \sqrt{k/2}) \, \d s = \frac{\sqrt{\pi}}{2} e^{-x^2/(2k)}\]
	(see \cite[p. 488]{GR}). 
	Combining everything, we find that 
	\[
	\rho_k(x) = \sqrt{\frac{2}{\pi k}} e^{-x^2/(2k)} +O\Big(\frac{1}{k^{3/2}}\Big)
	\]
as $q \to \infty$, independently of $x$, as was to be shown. \end{proof}

\subsection{Higher order asymptotics of the density function} 
Theorem \ref{MainTHM6} can be generalized to find more terms of an asymptotic expansion of $\rho_k(x)$ as $q \to \infty$, 
 for $k=2q^r$ with $q$ an odd prime, leading to the higher order terms in Theorem \ref{mainRW}. Following the proof of 
Theorem~\ref{MainTHM6}, we have to compute the integrals
\[
\int_0^{\infty} e^{-s^2} a_j(s) \cos(s x/\sqrt{k/2}) \, \d s.
\]
Since $a_j(s)$ is an even polynomial, one only has to compute integrals of the form
\[
\int_0^\infty s^{2 n} e^{-s^2} \cos(s u) \, \d s
\]
and these can be dealt with by differentiating the function 
\[
\int_0^\infty e^{-s^2} \cos(s u) \, \d s = \frac{\sqrt{\pi}}{2}e^{-u^2/4}
\]
(see \cite[p. 488]{GR})
$2n$ times with respect to $u$.
For instance, we find
\begin{align}
	\label{RHOkM2}
	\rho_k(x) &= e^{-x^2/(2k)} \Bigg (\frac{8 k^3-3 k^2+6 k x^2-x^4}{4 \sqrt{2 \pi } k^{7/2}}\Bigg)+ O\left(\frac{1}{k^{5/2}}\right),
\end{align}
where the implied constant is independent of $x$. This follows by performing the computation in the proof of Theorem \ref{MainTHM6} using $M = 2$.
This refined asymptotic expression \eqref{RHOkM2} allows us to give a better approximation to the Mahler measure. 

\begin{proof}[Proof of Theorem \ref{mainRW} for $k=2q^r$] Using the refinement \eqref{RHOkM2} we find
\begin{align*}
	\m(C_k) &= \int_1^{k} \rho_k(x) \log(x) \, \d x  \\
	&=\bigintssss_{\, 1}^k e^{-x^2/(2k)} \Bigg (\frac{8 k^3-3 k^2+6 k x^2-x^4}{4 \sqrt{2 \pi } k^{7/2}}\Bigg) \log(x) \, \d x + O\left(\frac{\log k}{k^{3/2}}\right)\\
	&\sim \bigintssss_{\, 1}^\infty e^{-x^2/(2k)} \Bigg (\frac{8 k^3-3 k^2+6 k x^2-x^4}{4 \sqrt{2 \pi } k^{7/2}}\Bigg) \log(x) \, \d x + O\left(\frac{\log k}{k^{3/2}}\right)
\end{align*}
to all orders in $q$.
We split the integral into two parts: $\int_0^1$ and $\int_0^\infty$. For the integral $\int_0^1$ we expand the integral into a series about $k = \infty$, which gives
\begin{align*} 
\bigintssss_{\, 0}^1 e^{-x^2/(2k)} \Bigg (\frac{8 k^3-3 k^2+6 k x^2-x^4}{4 \sqrt{2 \pi } k^{7/2}}\Bigg) \log(x) \, \d x & = \bigintssss_{\, 0}^1  \left( \sqrt{\frac{2}{\pi k }} + O\left(\frac{1}{k^{3/2}} \right)\right) \log(x) \, \d x \\ &= -  \sqrt{\frac{2}{\pi k }} +O\left(\frac{1}{k^{3/2}} \right).
\end{align*} 
For the integral $\int_0^\infty$, we differentiate 
\[
\int_0^\infty e^{-y^2} y^{s} \, \d y = \frac{1}{2}\Gamma \left(\frac{1+s}{2}\right)
\]
with respect to $s$ (see \cite[p. 337]{GR}), to find that 
\[
\int_0^\infty e^{-y^2} y^{2n} \log(y) \, \d y=  \frac{1}{4} \Gamma \left(n+\frac{1}{2}\right) \psi \left(n+\frac{1}{2}\right)
\]
for positive integers $n$, where $\psi$ is the digamma function. This leads to 
\[
\int_0^\infty e^{-x^2/(2k)} \Bigg (\frac{8 k^3-3 k^2+6 k x^2-x^4}{4 \sqrt{2 \pi } k^{7/2}}\Bigg) \log(x) \, \d x = \frac{\log k}{2} - \frac{\gamma}{2} - \frac{\log 2 }{2} + \frac{1}{4 k}.
\]
The conclusion is 
\[
\m(C_k) = \frac{\log k }{2} - \frac{\gamma}{2} - \frac{\log 2 }{2} + \sqrt{\frac{2}{\pi k }} + \frac{1}{4 k} + O\left(\frac{\log k }{k^{3/2}} \right), 
\] 
which finishes the proof of this case of Theorem \ref{mainRW}. 
\end{proof}

\begin{remark} 
A similar computation in \texttt{Mathematica} using $M = 10$ gives (with factored numerators and denominators)
\begin{align}
	\m(&C_k) =\frac{\log k}{2} - \frac{\gamma}{2} - \frac{\log 2 }{2} + 2 \frac{1}{(2\pi k)^{1/2}} + \frac{1}{2^2} \frac{1}{k}  -\frac{31}{2^2\cdot 3^2} \frac{1}{(2 \pi k)^{3/2}}+\frac{5}{2^3\cdot 3^2} \frac{1}{k^2}
	\nonumber \\&+\frac{13\cdot 71 }{2^6\cdot 3\cdot 5^2}\frac{1}{(2 \pi k)^{5/2}}-\frac{1}{2^3\cdot 3}\frac{1}{k^3} +\frac{423469}{2^9\cdot 3^3\cdot 5\cdot 7^2 } \frac{1}{(2 \pi k)^{7/2}}-\frac{29\cdot 59}{2^4\cdot 3^3\cdot 5^2}  \frac{1}{k^4}
	\nonumber \\&+\frac{13\cdot 5510303}{2^{14}\cdot 3^5\cdot 5\cdot 7} \frac{1}{(2 \pi k)^{9/2}}-\frac{101}{2^3\cdot 3^3\cdot 5}\frac{1}{k^5}-\frac{23\cdot 29\cdot 4202993}{2^{17}\cdot 3^2\cdot 5\cdot 7\cdot 11^2} \frac{1}{(2 \pi k)^{11/2}} 
	\nonumber \\&+\frac{153841}{2^4\cdot 3^5\cdot 7^2}\frac{1}{k^6}-\frac{1013699\cdot 3491520091}{2^{21}\cdot 3^5\cdot 5^3\cdot 7\cdot 11\cdot 13^2  } \frac{1}{(2 \pi k)^{13/2}}+\frac{118387}{2^4\cdot 3^4\cdot 5\cdot 7}\frac{1}{k^7}
	\nonumber \\&-\frac{2269\cdot 929444980559}{2^{24}\cdot 3^5\cdot 5^3\cdot 7\cdot 11\cdot 13}\frac{1}{(2 \pi k)^{15/2}} -\frac{95944159}{2^5\cdot 3^6\cdot 5^3\cdot 7}\frac{1}{k^8} \nonumber \\ &+\frac{6163\cdot 9656111\cdot 798140689}{ 2^{30} \cdot 3^6 \cdot 5^2 \cdot 11 \cdot 13 \cdot 17^2}\frac{1}{(2 \pi  k)^{17/2}}
	-\frac{2971\cdot 11681}{2^3\cdot 3^4\cdot 5^3\cdot 7}\frac{1}{k^9}+O\left(\frac{\log k}{k^{19/2}}\right). \label{moreterms}
\end{align}
The corresponding approximation of $\rho_k$ is too big to write down.

\end{remark}

\subsection{Asymptotics of $\m(C_k)$ if $k = q$ or $k = 2^r$} \label{subsecas}

In this subsection, we study the asymptotics of $\m(C_k)$ when $k = q$ is an odd prime tending to infinity, or when $k = 2^r$ for $r$ tending to infinity. Since the method is very similar to that of the previous subsections we will omit some details. First, we relate $\m(C_k)$ to more familiar linear and symmetrized linear Mahler measures, that have an easy probability density function. 

\begin{lemma} \label{lemSd} For $d \geq 1$, define $L_d \coloneqq x_1+\dots+x_d$ and $S_d \coloneqq x_1+1/x_1+\dots+x_d+1/x_d$. Then 
\begin{enumerate}
\item \label{lemSd1} for all $r \geq 1$, $\m(C_{2^r}) = \m(x_0+S_{2^{r-1}})$; 
\item \label{lemSd2} for all odd primes $q$, $\m(C_q) = \m(S_{q+1}) = \m(L_{q+1})$; 
\item \label{lemSd3} for all $d \geq 1$, 
$$\rho_{S_d}(x) = \Re  \frac{2}{\pi}\int_0^\infty I_0(2 i s)^d e^{i s x} \, \d s = e^{-x^2/(4 d)} \left( \frac{64 d^3-12 d^2+12 d x^2-x^4}{64 \sqrt{\pi } d^{7/2}} \right) + O\left(\frac{1}{d^{5/2}}\right)$$

\end{enumerate}
\end{lemma}

\begin{proof} 
Statement \eqref{lemSd1} follows from the explicit expression for $C_{2^r}$ in Example \ref{N2k}. For \eqref{lemSd2}, we substitute   $y_i =  x_i/x_0$ for $i=1,\dots,q-1$ and $y_q = 1/(x_0 x_1 \cdots x_{q-1})$ in the explicit expression from Example \ref{Nq}, to find 
$$ \m(C_q) = \m(1+y_1+\dots+y_q) = \m(L_{q+1}). $$
On the other hand, substituting $x_i=y_i/y_{q+1}$ for $i=1,\dots,q$ and $x_{q+1}=y_{q+1}$ in $S_{q+1}$ gives 
\begin{align*} \m(S_{q+1}) & = \m\Big( \frac{1}{y_{q+1}}(1+y_1+\dots+y_q) + y_{q+1}(1+\frac{1}{y_1}+\dots+\frac{1}{y_q})\Big) \\ 
& = \m\Big( \frac{1}{y_{q+1}}(1+y_1+\dots+y_q)(y^2_{q+1} + G)\Big) \mbox{ with } G\coloneqq  \frac{1+y_1+\dots+y_q}{1+y_1^{-1}+\dots+y_q^{-1}} \\ 
& = \m(L_{q+1}) + \m(y_{q+1}^2 + G). 
\end{align*}
Now we substitute $Y=y_{q+1}^2$ and use Jensen's formula to find 
$$ \m(y_{q+1}^2 + G) = \frac{1}{2} \bigintssss_{\T^{q+1}} \log |Y + G| \frac{\d Y}{Y} \, \d \mu_{\T^q} = \bigintssss_{\T^q} \log^+ |G| \, \d\mu_{\T^q} = 0, $$
since for $y_i \in \T$, $\overline{1+y_1+\dots+y_q} =1+\overline{y_1}+\dots+\overline{y_q} = 1+y^{-1}_1+\dots+y^{-1}_q$
and hence $|G| \equiv 1$ on $\T^q$. 

Statement \eqref{lemSd3} follows by taking the inverse Laplace transform of the exponential generating function \begin{equation} \label{ESd} E_{S_d}(t) = I_0(2t)^d. \end{equation} This last identity \eqref{ESd} is easy to show by following the proof of Lemma \ref{besselidentity}. The asymptotics then follows from Lemma~\ref{LemI0exp}, in a manner similar to the derivation of the identity \eqref{RHOkM2}.
\end{proof} 

\begin{proof}[Proof of Theorem \ref{mainRW} for $k=2^r$ and $k=q$]
If $k = 2^r$, we combine \eqref{Jensen} and Lemma \ref{lemSd}\eqref{lemSd1} to find 
\[
\m(C_k) = \int_1^{2^r} \log(x) \rho_{S_{2^{r-1}}}(x) \, \d x.
\]
On the other hand, for $k = q$, we find, using Lemma \ref{lemSd}\eqref{lemSd2} 
\[
\m(C_k) = \int_0^{2(q+1)} \log(x) \rho_{S_{q+1}}(x) \, \d x.
\]
The results follow by combining these integral expressions and the asymptotics of the density function from Lemma \ref{lemSd}\eqref{lemSd3}. \end{proof}

\begin{remark}
The leading term in the asymptotics of $\m(C_k)$ is $(\log k)/2$, compatible with a half-normal distribution. The higher order asymptotics of $k=2q^r$ (with $r$ fixed and $q \rightarrow + \infty$) and of $k=2^r$ (with $r \rightarrow + \infty$) are identical to all orders, but differ from that of the case where $k=q$ (with $q \rightarrow +\infty$), already in the constant term. 
\end{remark} 

\begin{remark}
The integral expression in Lemma \ref{lemSd}\eqref{lemSd3} allows for a good numerical approximation of the Mahler measure if $k = 2^r$ or $k = q$. An inequality for linear Mahler measure for `finite' $d$ of the form  $ \log d - \gamma/2 - 2 \leq \m(L_d) \leq \log d$ has been established in \cite[Thm.\ 1]{RVTV} (see also \cite{Laila}), and asymptotic estimates up to the constant terms can be found in \cite{MyersonSmyth} (corrected from \cite{SmythMP}) and \cite[Thm.\ 2]{RVTV}. 
\end{remark}

\begin{remark} 
The proof of Lemma \ref{lemSd} shows that 
	\begin{equation}
		 \label{eqmsym} 
		 \m(L_d) = \m(S_d) 
	\end{equation}  
for any $d>0$. 
This also follows from the more general identity
\begin{equation}
	\label{idZetaMahler} 
	Z(S_d; s) = \binom{s}{s/2}\, Z(L_d;s) \mbox{ with }  \binom{s}{s/2} = \frac{\Gamma(s+1)}{\Gamma(s/2+1)^2}
\end{equation} 
relating the corresponding zeta Mahler functions (and hence, probability density functions via the inverse Mellin transform), see Remark \ref{remzetaM}. Indeed, the equality \eqref{eqmsym} follows from \eqref{idZetaMahler} by taking derivatives with respect to $s$ and substituting $s=0$ (cf. \eqref{MMandMahlerZeta}), and noticing the power series expansion $\frac{\Gamma(s+1)}{\Gamma(s/2+1)^2} = 1 +O(s^2)$ in a neighbourhood of $s=0$. 

To prove the identity \eqref{idZetaMahler}, we will use the notation $\tilde H \coloneqq H(1/x_1,\dots,1/x_d)$ for any non-zero Laurent polynomial $H$, so that in particular $S_d = L_d + \tilde L_d$. In general, we have
\begin{equation} 
	\label{ZetaMahlerCT}
  Z(H;2n) = \text{CT}[(H \tilde H)^n]
  \end{equation} 
 for any $H \in \R[x_1^{\pm 1}, \ldots, x_d^{\pm 1}]$ with real coefficients. Indeed, for values $x \in \T$, we have that $1/x = \overline x$ is the complex conjugate of $x$. Therefore, we can rewrite 
 $$ Z(H,2n) = \int_{\T^d} |H|^{2n} \, \d \mu = \int_{\T^d} H^n \overline H^n \, \d \mu =  \int_{\T^d} (H \tilde H)^n \, \d \mu = \text{CT}[(H \tilde H)^n], $$
 where the last identity follows by expanding the Laurent polynomial into monomials. 

For notational simplicity, write $P=L_d$. We first prove the identity \eqref{idZetaMahler} for non-negative even integers $s=2n$. By \eqref{ZetaMahlerCT}, we find 
$$ Z(P + \tilde P; 2n) = \text{CT}[(P+\tilde P)^{2n}] = \text{CT}[\sum_{j=0}^{2n} \binom{2n}{j} P^j \tilde P^{2n-j}].$$
Now, since $P^j$ is a homogeneous polynomial of degree $j$ in $x_1,\dots,x_d$ and $\tilde P^{2n-j}$ is a homogeneous polynomial in $1/x_1,\dots,1/x_d$ of degree $2n-j$, the constant term in the last expression is given by just the term with $2n-j=j$, i.e., $j=n$, so that
$$ Z(P + \tilde P; 2n) = \text{CT}[\binom{2n}{n} P^n \tilde P^{n}],$$
which is the claimed identity.

Finally, to prove the identity \eqref{idZetaMahler} in full generality, consider the function 
$$ g(s) \coloneqq  Z(P+\tilde P; s) - \binom{s}{s/2}\, Z(P;s). $$
We aim to show that $g(s)$ is identically zero by using Carlson's theorem in the sharper version \cite[\S 5.8]{Tit}. For this, notice that by compactness of $\T^d$, for any Laurent-polynomial, $L_H=\max_{\mathbf x \in \T^d} |H(\mathbf x)|$ exists, so we get (A) the inequality $|Z(H;s)| \leq L_H^{\Re(s)}$. Next, we have (B) the inequality $|\frac{\Gamma(s+1)}{\Gamma(s/2+1)^2} | \leq 4^{\Re(s)}$, which follows by expanding $4^s = (1+1)^{2s}$. 
After these preparations, we can check the three conditions for Carlson's theorem (since $Z(H;s)$, so a fortiori $g(s)$, is analytic in $\Re(s) \geq 0$): 
(a) $g(n)=0$ for all $n \in \Z_{>0}$, which we have shown above; (b) $|g(s)|\leq c_1 e^{c_2|s|}$ for all $s$ with $\Re(s)>0$  for some positive real constants $c_1,c_2$, which follows by combining the two estimates (A) and (B) just obtained; and (c) $|g(iy)| \leq c_3 e^{c_4 |y|}$ for all real $y$ and some positive real constants $c_3,c_4$ with $c_4<\pi$, which follows from substituting $\Re(s)=0$ in the above two estimates (A) and (B).  
\end{remark}

\section{Accurate and fast numerical computation of $\rho_k$ and $\m(C_k)$ for small $k$} 
\label{MC6} 

The approximation in \eqref{RHOKfinal} is quite poor for small values of $k$. In this section, we focus on methods to compute both the probability density function $\rho_k$, as well as the value $\m(C_k)$ of the Mahler measure of the cyclovariety, for small, fixed $k$, to high accuracy. Our approach is based on using linear differential equations satisfied by the generating series for constant terms of powers of $F_k$.  

\subsection{The graph of $\rho_k$}
We start with the computation of $\rho_k$ for small $k$. Although the distribution converges to a half-normal distribution quite fast, for small $k$, the graphs of $\rho_k$ looks quite different. As a running example, we will take $\rho_6$, which is non-smooth; in fact, it is closely related to the `density of states' function for the tight-binding model on the triangular lattice in condensed matter physics, and the function has what in the physics literature is called `van Hove singularities'. 

Let $m \in \Z_{\geq 0}$. Set
\[
H_k(t) = \sum_{m \geq 0} \mathrm{CT}[F_k(\mathbf{x})^m] t^m
\]
for the ordinary generating function of the constant terms sequence of powers of $F_k$ (recall that in the previous section we instead used the exponential generating function \eqref{defeh}). Since $|F_k(\mathbf{x})| \leq k$ and the coefficients of $F_k$ are all positive, the series $H_k(t)$ has radius of convergence at least $1/k$. In the terminology of \cite{Lipshitz}, the series $H_k$ is the \emph{(complete) diagonal} of the rational function 
\[
\frac{1}{1-t x_1 \cdots x_{d} F_k(\mathbf{x})},
\]
and it follows that $H_k(t)$ satisfies a linear differential equation in $t$ \cite[Theorem~1]{Lipshitz}. This implies that $H_k(t)$ can be continued analytically to any point in the complex plane, except for the singularities of the differential equation.

We denote a corresponding annihilating differential operator by $\mathcal{L}_k$. 
To compute such an operator $\mathcal L_k$, we can use creative telescoping for holonomic functions, see \cite{Koutschan}. In theory, this algorithm always terminates, but in practice we managed to do this only for $k \leq 10$ with $k \neq 7,9$.

\begin{remark} \label{remPF} 
The function $H_k(t)$ is also called the \emph{fundamental period} of $F_k$ and hence it satisfies a Picard--Fuchs equation \cite{Christol} corresponding to the family of algebraic varieties $t x_1 \cdots x_{d} F_k(\mathbf{x})=1$, with $t$ as parameter. The Picard--Fuchs operator is a possible $
\mathcal L_k$. To compute it (or a multiple of it), one may use Lairez generalisation of the Griffiths--Dwork algorithm \cite{Lairez} (implemented in \texttt{Magma}). 

The fact that $H_k(t)$ is such a period also implies that, if the corresponding cyclopolytope is reflexive, the coefficients $\mathrm{CT}[F_k(\mathbf{x})^m]$ satisfy higher versions of Lucas congruences, see \cite[Theorem 3.3]{SamolvStraten}.
\end{remark} 

\begin{remark} 
The constant term sequence $\text{CT}[F_k(\mathbf x)^m]$ has the following combinatorial interpretation: it is the number of closed walks of length $m$ in the complex plane, where all steps are given by $k$-th roots of unity, see \cite[Thm.\ 2.2]{Lala}. 
\end{remark}

The main theoretical result underlying our algorithm to compute $\rho_k(x)$ is the following.

\begin{proposition} 
\label{prop: analytic_H_k}
%
Suppose $F_k$ is reciprocal. Then $H_k$ admits an analytic continuation to $\Cc \setminus ( (-\infty, -1/k] \cup [1/k, \infty ))$.
Write $H_k^{+}(t)$ for the analytic continuation of $H_k(t)$, where the analytic continuation is done via the complex \emph{upper} half plane and similarly write $H_k^{-}(t)$ for the analytic continuation of $H_k(t)$, where the analytic continuation is done via the complex \emph{lower} half plane. Then 
 for $x > 0$,
\begin{equation}
  \label{pdfidentity}
  \rho_k(x) = \frac{1}{2 \pi i x} \left( H^{+}_k(1/x) + H^{-}_k(-1/x) - H^{-}_k(1/x) - H^{+}_k(-1/x)\right).    
\end{equation}
\end{proposition} 

\begin{proof} 
The constant terms have the integral representation
\[
\mathrm{CT}[F_k(\mathbf{x})^{m}] = \int_{\T^d} F_k(\mathbf{x})^m \, \d \mu,
\]
which implies
\begin{equation}
\label{eq:H_k_asintegral}
H_k(t) = \frac{1}{t}\int_{\T^d} \frac{1}{1/t - F_k(\mathbf{x})} \, \d \mu.
\end{equation}
Since $F_k$ is a reciprocal Laurent polynomial with integer coefficients, it is real-valued on the torus. Together with the bound $|F_k| \leq k$, this implies that the integral expression \eqref{eq:H_k_asintegral} is analytic on $\Cc \setminus ( (-\infty, -1/k] \cup [1/k, \infty ))$. This proves the first part.

Recall from \eqref{ConstantTermAsIntegral} the relation 
\[
\mathrm{CT}[F_k(\mathbf{x})^{2m}] = \int_0^k x^{2m} \rho_k(x) \, \d x
\]
between constant terms and even moments, which implies that we can write 
\begin{align*} 
  \frac{1}{2 z} \left( H_k(1/\sqrt z) + H_k(-1/\sqrt z) \right)  & =     \sum_{m \geq 0} \frac{ \mathrm{CT}[F_k(\mathbf{x})^{2m}]}{z^{m+1}} \\ & 
  = \frac{1}{z}\sum_{m \geq 0} \frac{1}{z^m} \int_{-\infty}^{+\infty} x^{2m} \rho_{k}(x) \, \d x =  \int_{- \infty}^{+\infty} \frac{\rho_k(x)}{z-x^2} \, \d x\\
  &= \int_{- \infty}^{+\infty} \frac{\rho_k(\sqrt{t})}{2\sqrt{t}(z-t)} \, \d t.
\end{align*}
We recognize the right-hand side as the Stieltjes transform of the density function $$f(t) \coloneqq  \rho_k(\sqrt{t})/2\sqrt{t}  \colon (0,k^2] \rightarrow \R$$ using as integration contour the interval in the real line where the integrand is non-zero and defined. 
Define 
\[
\Psi_k(z) \coloneqq \frac{1}{2 z} \left( H_k(1/\sqrt z) + H_k(-1/\sqrt z) \right), 
\]
which is an analytic function on $\Cc \setminus [0 , k^2] $.
The Stieltjes inversion formula enables us to express $f(t)$ in terms of $\Psi_k(z)$, see \cite[page 40]{Nica_Speicher_2006}: 
for all $t \in (0, k^2)$, we have
 \[
2 \pi i f(t) = \lim_{b \to 0^{-}} \Psi_k(t+i b) - \lim_{b \to 0^{+}} \Psi_k(t+i b). 
\]
Since $H_k$, and hence $\Psi_k$, has only finitely many singularities, the limit on the right-hand side of the above equation exists for all but finitely many values $t \in (0, k^2)$.
The result now follows from 
\begin{align*}
  2 \pi i \cdot \frac{\rho_k(\sqrt{t})}{2 \sqrt{t}} &=  \lim_{b \to 0^{-}} \Psi_k(t+i b) - \lim_{b \to 0^{+}}  \Psi_k(t+i b) \\
  &= \frac{1}{2 t} \left(H^{+}_k(1/\sqrt{t}) + H^{-}_k(-1/\sqrt{t}) - H^{-}_k(1/\sqrt{t}) - H^{+}_k(-1/\sqrt{t}) \right).
\end{align*}
by substituting $x=\sqrt{t}>0$. 
\end{proof} 

\begin{remark}
If $F_k$ is \emph{non-reciprocal}, then $\rho_k$ can be recovered from the generating function of the constant terms of the \emph{reciprocal} Laurent polynomial $F_k \tilde F_k$. We will briefly outline this procedure. 

Define $\hat{H}_k = \sum_{m \geq 0} \mathrm{CT}[(F_k \tilde{F_k})^m] t^m$ and define $\hat{\rho}_k \coloneqq \rho_{F_k \tilde F_k}$.
	Then, using the notation of Proposition \ref{prop: analytic_H_k}, we have
	\begin{equation}
		\hat{\rho}_k(x) = \frac{1}{2 \pi i x} \left( \hat{H}^{+}_k(1/x) + \hat{H}^{-}_k(-1/x) - \hat{H}^{-}_k(1/x) - \hat{H}^{+}_k(-1/x)\right).    
	\end{equation}
Following Remark \ref{remzetaM}, we obtain for $\Re(s) \geq 0$
\[
\int_0^\infty x^s \hat{\rho}_k(x) \, \d x = Z(F_k \tilde F_k; s) = Z(F_k; 2s) = \int_0^\infty x^{2s} \rho_k(x) \, \d x.
\]
Comparing the outer expressions of this equality, it follows that $\rho_k(x) = 2 x \hat{\rho}_k(x^2)$.
\end{remark}

\begin{example}[$k=2$]
\label{excompk2} As a warm-up, we first show an example where the analytic continuation of solutions to $\mathcal L_k$ is easy to do in closed form. We consider $k=2$, i.e., $F_2 =  x_1 + x_1^{-1}$. In this case 
$
    H_2(t) =  \sum_{m \geq 0} \binom{2m}{m} t^{2m} = 1/\sqrt{1-4t^2}, 
 $
which converges for $|t|<1/2$ and is annihilated by the operator $\mathcal L_2 = (4t^2-1)\d+4t$ with $\d = {\d}/{\d t}$.
Note that $H_2^{+}(t) = - H_2^-(t) = -\frac{1}{\sqrt{1-4t^2}}$, so 
\[
\rho_2(x) = \frac{2}{\pi} \frac{1}{\sqrt{4-x^2}}
\mbox{ 
and 
} 
\m(C_2) = \frac{2}{\pi}\int_{1}^2  \frac{\log(x)}{\sqrt{4-x^2}} \, \d x= 0.323066 \ldots.
\]
\end{example}

In general we might not be able to find an explicit description for the analytic continuation of $H_k(t)$; however, we can \emph{numerically} continue analytically to arbitrary precision, as in the following example. 

\begin{example}[$k = 6$]
  In this case
  \[
  F_6(\mathbf{x}) = x_1 + 1/x_1 + x_2 + 1/x_2 + x_1/x_2 + x_2/x_1
  \]
  and the sequence of constant terms $\mathrm{CT}[F_k(\mathbf{x})^m]$ is recorded in \cite[A002898]{oeis}. A closed form for the solution \emph{is} known in this case--see below--, but we will not use it for the computation. The ordinary generating function $H_6(t)$ is annihilated by the second order linear differential operator (found by creative telescoping) 
  \[  \mathcal L_6 = t(2t+1)(3t+1)(6t-1)\d^2+(6t+1)(24t^2+8t-1)\d + 24t(3t+1),
  \]
  where $\d = \frac{\d}{\d t}$. The singular points of the differential equation are at $t = 0, -\frac{1}{2}, -\frac{1}{3}$ and $\frac{1}{6}$, so $H_6(t)$ admits analytic continuation to $\Cc \setminus ( (-\infty, -1/6] \cup [1/6, \infty ))$.
  
The support of $\rho_6(x)$ is contained in $[0,6]$. Formula \eqref{pdfidentity} can be used to compute any value to arbitrary precision;  we explain the example of $\rho_6(4)$. To employ \eqref{pdfidentity}, we need to compute the four limits  $\lim_{b \to 0^{\pm}} H_6(\pm 1/4 + i b).$ For example, choosing both signs positive, we need to analytically continue $H_6(t)$ from $t = 0$ to $t = \frac{1}{4}$, where the path is in the complex \emph{upper} half plane.  This can be done, e.g., with the \texttt{OreAlgebra} package \cite{OREALG} in \texttt{SageMath}, using the routine
  \begin{center}
    \begin{lstlisting}[language=Python, label=compX]
      from ore_algebra import *
      Pol.<t> = QQ[]
      Dop.<Dt> = OreAlgebra(Pol)
      dop = (36*t^4 + 24*t^3 + t^2 - t)*Dt^2 + (144*t^3 + 72*t^2 + 2*t - 1)*Dt + (72*t^2 + 24*t)
      dop.numerical_solution([0,1], [0,I,1/4])\end{lstlisting}
  \end{center} 
 (in the last line of the code, $[0,1]$ indicates that we pick the unique holomorphic solution around $t = 0$ with constant term $1$, and $[0, I ,1/4]$ indicates the path from $0$ to $i$ and $i$ to $1/4$ connected by straight lines). We find
\begin{equation} \label{hb6}
  \lim_{b \to 0^{+}} H_6(1/4 + i b) = 0.76072185063246146  \ldots + 0.69719958972491171  \ldots \cdot i.
  \end{equation} 
Repeating the procedure, we find the same value \eqref{hb6} for   $\lim\limits_{b \to 0^{-}} H_6(1/4 + i b)$, and 
  \[  \lim_{b \to 0^{-}} H_k(-1/4 + i b)  =   \lim_{b \to 0^{+}} H_k(-1/4 + i b) = 1.4477785313398085  \ldots\]   Using \eqref{pdfidentity} we immediately find 
  $
  \rho_6(4) = 0.05548138051318054 \dots.
  $
  Figure \ref{Fig:rho6} displays the graph of $\rho_k(x)$ based on running this computation at $1000$ sample points in the interval $[0,6]$ (avoiding the finitely many singular points). Given $\mathcal L_6$, the computation of this list of values took under $9$ minutes on a standard user-grade laptop  (Intel i5, 16GB  RAM). 
  
  \begin{figure}
    \centering
    \includegraphics[scale = 0.55]{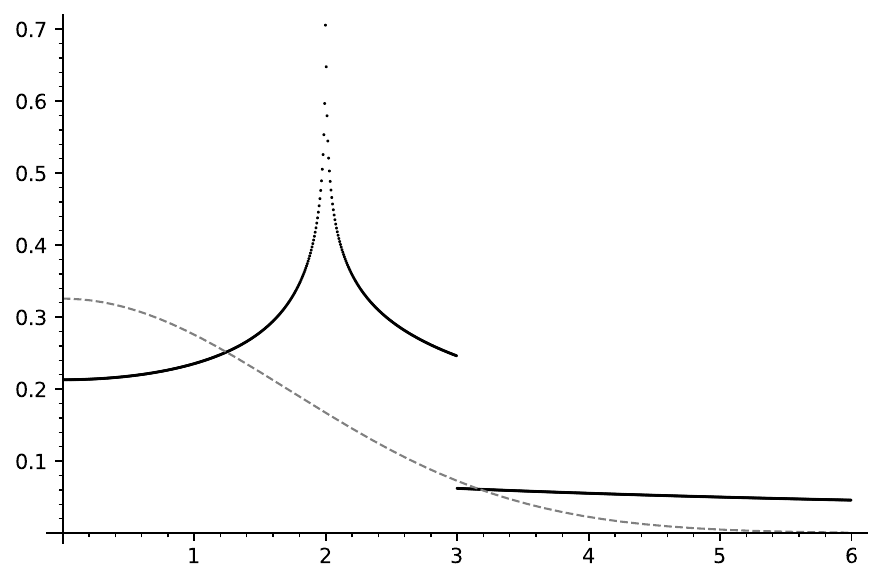}
    \caption{The graph of $\rho_6(x)$, with the corresponding half-normal distribution superimposed (dashed line). The points of discontinuity at $x= 2$ and $x = 3$ can be explained by the fact that $\mathcal{L}_6$ has singularities at the points $-\frac{1}{2}$ and $-\frac{1}{3}$.}
    \label{Fig:rho6}
  \end{figure}
\end{example}

\begin{remark} There are closed formulas for $H_6(t)$ in terms of hypergeometric functions or elliptic integrals, leading to closed formulas for the probability density function $\rho_6(x)$. The function $\rho_6(x)$ is related to the density of states function $\rho_{\triangle}$  for the tight-binding model on a triangular lattice (and $H_6(t)$ is related to the Green's function for that lattice) in condensed matter physics, through the identity
\begin{equation} \label{sumrhotriangle}
  \rho_6(x) = \rho_{\triangle}(x) + \rho_{\triangle}(-x).
  \end{equation} 
We took the explicit expression for $\rho_\triangle$ from \cite{Kogan} (the expression of the corresponding Green's function in terms of elliptic functions were first stated by Hobson and Nierenberg \cite{HN}). 
By substituting this into \eqref{sumrhotriangle}, we find an explicit expression for $\rho_6(x)$ in terms of elliptic integrals, as follows. 
Let $K$ denote the complete elliptic integral of the first kind using the convention
$$ K(x) \coloneqq  \bigintssss_0^{\pi/2} \frac{\d\theta}{\sqrt{1-x^2 \sin^2 \theta}} = \frac{\pi}{2}\,  \hypgeo{2}{1}\Big(\frac12, \frac12,1;x^2\Big),$$ and define the two functions 
\begin{align*} 
F_1(x) &\coloneqq \frac{2}{\pi^2 \sqrt{(3-\sqrt{x+3})(1+\sqrt{x+3})^3)}}\, K \Bigg(\frac{4\sqrt[4]{x+3}}{\sqrt{(3-\sqrt{x+3})(1+\sqrt{x+3})^3)}} \Bigg), \\[4mm] 
F_2(x) &\coloneqq \frac{1}{2 \pi^2 \sqrt[4]{x+3}}\, K \Bigg( \frac{\sqrt{(3-\sqrt{x+3})(1+\sqrt{x+3})^3)}}{4\sqrt[4]{x+3}} \Bigg).  
\end{align*} 
Then 
\[
\rho_6(x) = \left\{ \begin{array}{ll} 
F_2(x)+F_2(-x) & \mbox{ for } 0 \leq x < 2, \\
F_2(x)+F_1(-x) & \mbox{ for } 2 <x < 3, \\
F_2(x)& \mbox{ for } 3 < x \leq 6.
\end{array} \right. 
\]
If one plots $\rho_6(x)$ using this analytic expression, one finds an exact fit with Figure \ref{Fig:rho6}. 
\end{remark}

\begin{remark} 
Cousins of the functions $\rho_2$ and $\rho_6$ have also occurred in the engineering literature, in the study of the spread of the spectrum of the received signal in a noisy communication between a fixed base station and a mobile receiver, and a mobile sender and receiver, leading to $\rho_2$ and $\rho_6$ as representing the power spectrum of the complex envelope of the received signal, respectively, see, e.g., \cite[p.~6, 2nd column]{Haber}. 
\end{remark}

\begin{figure}
  \begin{subfigure}[h]{0.45\linewidth}
    \includegraphics[width=\linewidth]{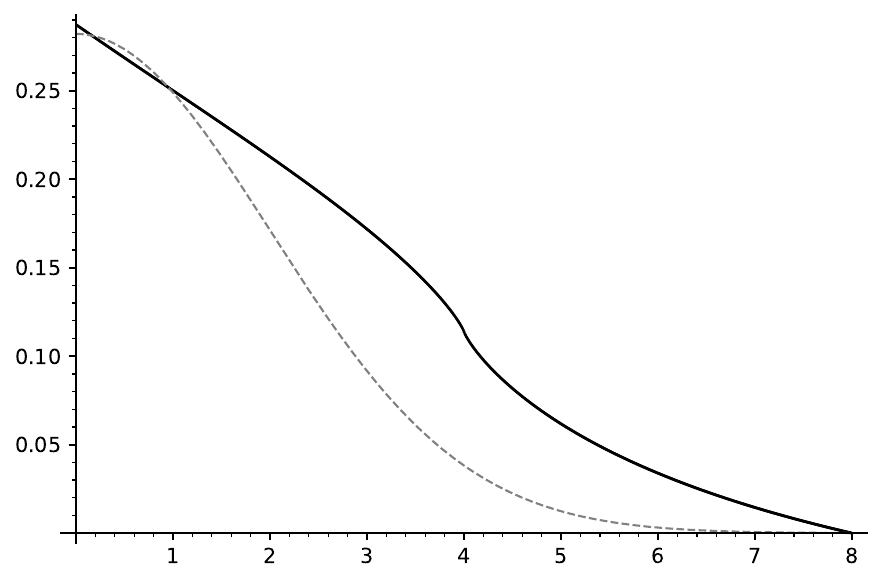}
    \caption{The graph of $\rho_8(x)$}
  \end{subfigure}
  \hfill
  \begin{subfigure}[h]{0.45\linewidth}
    \includegraphics[width=\linewidth]{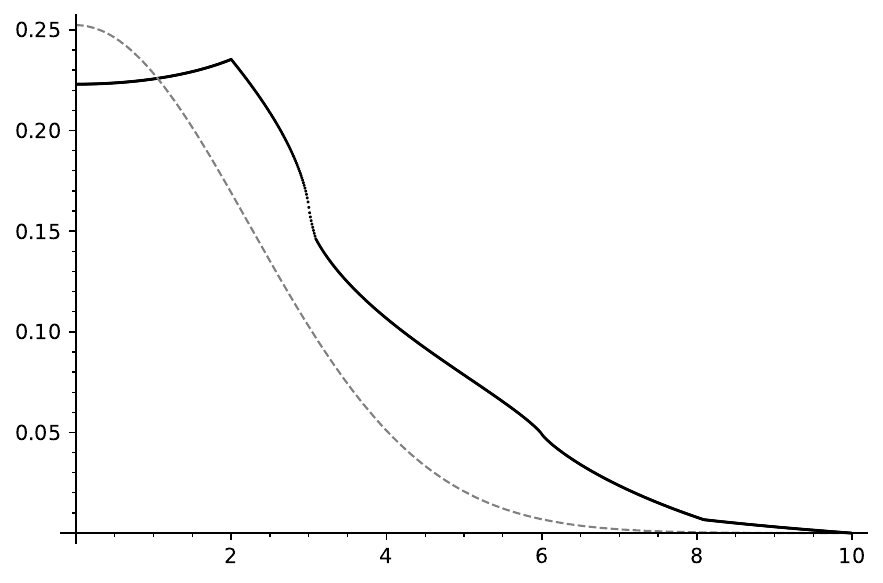}
    \caption{The graph of $\rho_{10}(x)$}
  \end{subfigure}%
  \caption{Graphs of two more probability density functions computed by numerical analytic continuation, with the corresponding half-normal distribution superimposed (dashed line).} \label{Fig:rho810}
\end{figure}

\begin{example}
By following the same strategy, computing $\mathcal L_k$ and using computer algebra to determine values of analytic continuations of the solution $H_k$ for $k=8, 10$, we can also plot the graphs of $\rho_8$ and $\rho_{10}$, see Figure \ref{Fig:rho810}. Note that the graph of $\rho_{10}$ already quite closely resembles the half-normal distribution $\sqrt{\frac{2}{10 \pi}} e^{-x^2/20}$, as predicted by Theorem \ref{MainTHM6}. Again writing $d = \d/\d t$, in these cases the differential operators are given by 
\begin{align*} \mathcal L_8 =\;& t^3(4t-1)(4t+1)(8t-1)(8t+1)d^4 + 2t^2(64t^2-3)(112t^2-1)d^3 \\ &+ t(55296t^4 - 2048t^2 + 7)d^2 + (61440t^4 - 1344t^2 + 1)d + 128t (96t^2-1) \end{align*}
and
\begin{align*}
	\mathcal{L}_{10} =\;& 
	t^3 (2t - 1)(3t - 1)(6t + 1)(10t - 1)(90t - 23)^2 (25t^2 - 5t - 1) d^4 \\
	&+ 2t^2 (90t - 23) 
	(6480000t^7 - 7344000t^6 + 2610600t^5 - 152120t^4 \\
	&\quad\quad - 83118t^3 + 10880t^2 + 548t - 69) d^3 \\
	&+ t (5248800000t^8 - 6823440000t^7 + 3244068000t^6 - 607966800t^5 \\
	&\quad\quad - 2477080t^4 + 14831484t^3 - 1231504t^2 - 44114t + 3703) d^2 \\
	&+ (6998400000t^8 - 8475840000t^7 + 3805992000t^6 - 710474400t^5 \\
	&\quad\quad + 20154160t^4 + 9778296t^3 - 889272t^2 - 6394t + 529) d \\
	&+ 1749600000t^7 - 1963440000t^6 + 822024000t^5 - 147444000t^4 \\
	&\quad\quad + 6697440t^3 + 1010160t^2 - 84640t.
\end{align*}
Notice that $t=23/90$ is an (in fact, the only) apparent singularity of $
\mathcal L_{10}$. 
\end{example}

\begin{remark} \label{realsing} We observe that in all our examples, the singularities of the operator $\mathcal L_k$ are real, though not necessarily rational. We do not know whether it is true in general that an annihilating differential operator for $H_k$ exists for which all singularities are real. 

It is not true that the Picard--Fuchs operator of \emph{any} reciprocal polynomial has real singularities only, as the following example shown to us by Frits Beukers demonstrates: the Picard--Fuchs operator associated to the reciprocal polynomial \[R(x_1,x_2,x_3) \coloneqq x_1+x_2+x_3+x_1x_2x_3/2+1/x_1+1/x_2+1/x_3+1/(2x_1x_2x_3)\] is 
\begin{align*}
&t^4(1-t^2)(1-49t^2)(1+t^2)(-4+21t^2) d^4 + t^3(17493t^8-3976t^6-7147t^4+1806t^2-16)d^3 \\
   &+ t^2(83349t^8-20580t^6-9466t^4+3659t^2-4)d^2 \\ &+ t(117306t^8-32004t^6-682t^4+1557t^2+4)d+ 2t^4(15435t^4-4746t^2-16),
\end{align*}
which has a non-apparent singularity at $\pm i$. Note that $x_0+R(x_1,x_2,x_3)=0$ defines a Laurent hypersurface of dimension $3$; cyclovarieties, on the other hand, have dimension lying in the range of the Euler totient function, which does not contain $3$.  

The singularities have an algebro-geometric interpretation. As observed in Remark \ref{remPF}, we can choose $\mathcal L_k$ as the Picard--Fuchs operator annihilating the constant term generating series $H_k(t)$, seen as the period of a holomorphic differential form on the fibers in the family $C_k(\lambda) \rightarrow \PP^1_\lambda$ (but with a change of variables $t=-\lambda^{-1}$). The general theory of regularity for Picard--Fuchs equations arising from periods \cite[\S 2]{vS} implies that $\lambda$ for which $C_k(\lambda)$ is singular correspond to some non-apparent singularities of $\mathcal L_k$. As a side-remark, by \eqref{pdfidentity}, the density function $\rho_k(x)$ is non-smooth at most at values $x$ corresponding to $|1/t|$ for $t$ singular points of $\mathcal L_k$ (since the support is in the positive real numbers). 

For example, if $k=6$ the singular fibers are at $\lambda=2,3,-6,\infty$, and the singular points of $\mathcal L_6$ are, correspondingly, at $-1/2,-1/3,1/6$ and $0$, while $\rho_6(x)$ is not smooth at $x=0,2,3,6$. For $\mathcal L_{10}$, the fiber $C_{10}(-90/23)$ corresponding to the non-apparent singularity $23/90$ of $\mathcal L_{10}$ is not singular.  
 
\end{remark} 

\subsection{Computing Mahler measures; the value $\m(C_6)$}

Computing Mahler measures numerically up to high accuracy directly from the definition can be excruciatingly slow; the main reason is that the logarithm to be integrated blows up near the zeros of the defining polynomial on the torus $\T^d$. For $\m(C_k)$, we know this issue appears as well, see Proposition \ref{toricpointsDS}. A way around the problem is to rely on the techniques from the previous subsection. More precisely, 
we illustrate how to use the differential operator $\mathcal L_k$ to compute $\m(C_k)$ for even $k$ fast and to high accuracy. 

Write $G_k$ for the cumulative density function 
\[
G_k(x) = \int_{0}^x \rho_k(y) d y, 
\]
and note that it has support in $[0,k]$ and total integral $1$. 
By using integration by parts, for the Mahler measure we have
\begin{align}
	\label{MMIntegrationByParts}
\m (C_k) & = \int_1^{k} \rho_k(x) \log(x) d x = [G_k(x) \log(x)]_{1}^k - \int_1^{k} G_k(x) \frac{d x}{x} \nonumber \\ &= \log(k) - \int_1^{k} G_k(x) \frac{d x}{x}.
\end{align}
Write \[
G_k^{\pm}(x) = \int_{1/x}^{\infty} H_k^{\pm}(u) \frac{d u}{u},
\mbox{ 
so that } 
G_k(x) =\frac{1}{2 \pi i} \left( G^{+}_k(x) + G^{-}_k(-x) - G^{-}_k(x) - G^{+}_k(-x) \right)
\]
by \eqref{pdfidentity}, and define
\[
B^{\pm}_k(t) \coloneqq  \int_{t}^\infty G_{k}^{\pm}(1/v) \frac{d v}{v}
\mbox{ 
and }
B_k(t) \coloneqq  B^{+}_k(t) + B^{-}_k(-t) - B^{-}_k(t) - B^{+}_k(-t).
\]
After rewriting everything using \eqref{MMIntegrationByParts}, we conclude that
\begin{equation}
	\label{MahlerNumerical}
	\m (C_k) =  \log(k) - \frac{1}{2 \pi i} \left(B_k(1/k) - B_k(1)\right). 
\end{equation}

Write $\theta$ for the differential operator $t \frac{d}{d t}$, and 
observe that 
$
\theta^2 \, B_k^{\pm}(t) = \theta \, G_k^{\pm }(1/t) = H_k^{\pm}(t).
$
This implies that $B_k^{\pm}(t)$ is annihilated by the product of differential operators $\mathcal{L}_k \theta^2$. Since $\ker(\theta^2-1) = \{ \log x, (\log x)^2 \}$, locally around $t = 0$, the function $B_k^{\pm}(t)$ looks like
\begin{equation}
	\label{Bklocalexp}
	B_k^{\pm}(t) = c_0 + c_1 \log(t) + \frac{1}{2}\log(t)^2 + \sum_{m \geq 1} \frac{\mathrm{CT}[F_k(\mathbf x)^m]}{m^2} t^m
\end{equation}
for some constants $c_0, c_1$. As for positive $a$, \[\lim_{\epsilon \to 0^{+}} \left( \log(a+\epsilon i) + \log(-a-\epsilon i)-\log(a- \epsilon i) - \log(-a + \epsilon i) \right) = - 2 \pi i
\]
is independent of $a$ whatever value is chosen for $c_0,c_1$, the result in \eqref{MahlerNumerical} is independent of such a choice, and hence we may and will choose them to be zero. 
\begin{example}[$k = 6$]
	We use the differential operator \emph{dop} from the Code \ref{compX}. The routine
	
	\begin{center}
		\begin{lstlisting}[language=Python, label=compY]
			theta = t*Dt
			dopmahler = dop*theta^2
			print(dopmahler.local_basis_monomials(0))\end{lstlisting}
	\end{center} 
	returns 
	\[
	[1/6 \log(t)^3, 1/2 \log(t)^2, \log(t), 1].
	\]
	It follows that $B_k^{\pm}(t)$ (with our choice $c_0 = c_1 = 0$) is the unique solution of $\mathcal{L}_6 \theta^2$ that looks like $\frac{1}{2} \log(t)^2 + O(t)$ around $t=0$.  We implement the function $B(t)$ as follows, where \emph{prec} indicates the desired precision:
	\begin{center}
		\begin{lstlisting}[language=Python, label=compZ]
			Init = [0,1,0,0] 
			def B(x, prec):
			v1 = dopmahler.numerical_solution(Init, [0,I,x],prec)
			v2 = dopmahler.numerical_solution(Init, [0,-I,-x],prec)
			v3 = dopmahler.numerical_solution(Init, [0,-I,x],prec)
			v4 = dopmahler.numerical_solution(Init, [0,I,-x],prec)
			return v1 + v2 - v3 - v4\end{lstlisting}
	\end{center} 
After these preparations, we can use \eqref{MahlerNumerical} to compute the Mahler measure $\m(C_k)$ (with $k=6$), as follows: 
	\begin{center}
		\begin{lstlisting}[language=Python, label=compB]
			prec = 10^(-100)
			(log(k) - 1/(2*pi*I)*(B(1/k,prec)-B(1,prec))).real().n(200)
		\end{lstlisting}
	\end{center} 
This returns 
	\[
	\m(C_6) = 0.64394320995350600341722247622009734279536970069336928444715.
	\]
Given $\mathcal L_6$, the computing time to find the first 100 decimals of $\m(C_6)$ was around $2$ seconds on a standard user-grade laptop (Intel i5, 16GB  RAM). 
\end{example}

\begin{remark}
As a word of caution, we would like to remind the reader that we are computing the Mahler measure of $x_0+F_k(\mathbf x)$ with $x_0$ \emph{included in the set of torus integration variables}, not as a parameter. One can find in the literature closed hypergeometric expressions for $\m(\lambda+F_6(x_1,x_2))$ as a function of a parameter $\lambda$ (for sufficiently large $|\lambda|$, see, e.g., \cite[(2-37)]{LalinRogers}), but this is not the value of $\m(C_6)$, which is given by a further integration (in $\lambda$) of those expressions. 
\end{remark}

\section{An algorithm to compute the minimal Mahler measure $m_q$; proof of Theorem~\ref{mainMM}} \label{sectalg} 

Let $\mathfrak C_n$ be the set of positive Mahler measures attained by degree $n$ cyclotomic integers with cyclic Galois group. Recall from \eqref{eq:m_n-t_n-definition} that $m_n$ is the minimum of $\mathfrak C_n$.

Let $q$ be an odd prime. In this section, we devise an algorithm (see \S \ref{subsec:The algorithm}) that searches for $m_q$. 
We implement this algorithm to prove Theorem~\ref{mainMM}, which concerns the cases $q=3$, $5$, and $7$. As $q$ increases, more and more sophisticated techniques are necessary to make the search feasible, and we group these into the first four subsections below. As will become clear below, the approach is to derive bounds on the conductor, as well for individual conductor exponents and congruences for prime divisors of the conductor. 
 In this section, it is more practical to work with the exponential Mahler measure $\MM(f) = \exp \m(f)$ instead of the logarithmic one. 

\subsection{First tool: naive coefficient bound} 

First of all, there exists an algorithm that, given $n$, will theoretically find $m_n$ in finite time: this is because the set of polynomials of degree $q$ with bounded Mahler measure is finite, as the following elementary lemma establishes. 

\begin{lemma}[{\cite[Lemma~1.6]{BrunaultZudilin}}]
\label{lem:Mahler coefficient bound}
Suppose $f(x) = a_0 + \dots + a_n x^n \in \Z[x]$ and $\MM(f) \leq B$. Then $|a_{j}| \leq \binom{n}{j} B$ for each $j=0,\dots,n$. \hfill \qed
\end{lemma}

We apply this with $B=\MM(\alpha_n)$. For example, the minimality of the Mahler measures of $\alpha_3$ and $\alpha_5$ can easily be verified with a naive search based on Lemma~\ref{lem:Mahler coefficient bound} applied to their minimal polynomial. 

In degree $7$, the minimal polynomial of $\alpha_7 = \zeta_{29}+ \zeta^{12}_{29} + \zeta^{17}_{29} + \zeta^{28}_{29}$ is
\begin{equation*}
x^7 + x^6 - 12x^5 - 7x^4 + 28x^3 + 14x^2 - 9x + 1,
\end{equation*}
which has exponential Mahler measure $ \MM(\alpha_7) \approx 24.21657$ (and Mahler measure $\m(\alpha_7) \approx 3.187$). We only need to consider monic polynomials, and after replacing $f(x)$ by $-f(-x)$ we may also assume that the constant coefficient of $f$ is positive. Even with these additional reductions, it is computationally infeasible to check all remaining $\approx 8 \cdot 10^{18}$ polynomials satisfying Lemma~\ref{lem:Mahler coefficient bound} with input $B = \MM(\alpha_7)$. To prove Theorem~\ref{mainMM}, we need to drastically curb the search space. Doing this by a factor around $10^{9}$ is a task that will keep us busy for the remainder of this section.  

We will use the bound from Lemma \ref{lem:Mahler coefficient bound} in steps \eqref{alg-ii}, \eqref{alg-iii} and \eqref{alg-iv} of the algorithm.  

\subsection{Second tool: the discriminant} 
 In this subsection, we state two results regarding the zeros and the discriminant of polynomials $f \in \Z[x]$ with cyclic Galois group of odd order, and an inequality relating the latter to the (exponential) Mahler measure of $f$. 
 
 \begin{lemma}
\label{lem:Cn real roots}
Let $f \in \Z[x]$ be a polynomial of odd degree $n$ with Galois group $\Gal(f) = \Z/n\Z$. Then $f$ has only real zeros and square discriminant.
\end{lemma}
\begin{proof}
The group $\Z/n\Z$ is of odd order, so that the automorphism corresponding to complex conjugation -- generically of order $2$ -- must act trivially on the zeros of $f$. Hence the zeros of $f$ are real. Furthermore, by transitivity (since $f$ is irreducible), the group $\Gal(f)$ is generated by an $n$-cycle. This is an even permutation, so $\Gal(f)$ is contained in the alternating group on $n$ letters. 
\end{proof}

Mahler \cite{Mahler} showed that $\MM(f)^{2n-2} \geq \abs{\Delta(f)}/n^n$ for a polynomial $f \in \Z[x]$ of degree $n$ with discriminant $\Delta(f)$. When $f$ has all its roots in $\R \setminus \{0, \pm 1\}$, Schinzel \cite{Schinzel} proved that $\MM(f) \geq ((1+\sqrt{5})/2)^{n/2}$; A simple proof was later provided by 
H\"ohn and Skoruppa \cite{HS}, and H\"ohn also gave a short proof of a generalization \cite{Hohn}. The next lemma is an alternative lower bound in this setting relating the Mahler measure and discriminant of $f$.

\begin{lemma}
\label{lem:disc Mahler bound}
Suppose $f$ is a monic, squarefree polynomial of degree $n$ all of whose roots are real and such that $f$ does not vanish at $1$ and at $-1$. Then $0 < \abs{f(1)f(-1)}\Delta(f) < \MM(f)^{2n}$.
\end{lemma}
\begin{proof}
From the definition of the discriminant, it follows that $\Delta(f) > 0$. Denote the roots of $f$ by $\alpha_1, \ldots, \alpha_n$. Consider the Vandermonde matrix
\begin{equation*}
V = \begin{pmatrix}
1 & \alpha_1 & \alpha_1^2 & \dots & \alpha_1^{n-1}\\
1 & \alpha_2 & \alpha_2^2 & \dots & \alpha_2^{n-1}\\
\vdots & \vdots & \vdots & \ddots &\vdots \\
1 & \alpha_n & \alpha_n^2 & \dots & \alpha_n^{n-1}\\
\end{pmatrix}
\end{equation*}
and recall that $\Delta(f) = \det(V)^2 = \det(V^{\mathsf{T}})^2$. Hadamard's inequality (see, e.g., \cite[Lemma 2.7]{BrunaultZudilin}) thus yields
\begin{equation*}
0 < \Delta(f) = \det(V^{\mathsf{T}})^2 \leq \prod_{j=1}^n 
\Big( \sum_{k=1}^n \big|\alpha_j^{k-1}\big|^2 \Big)  
= \prod_{j=1}^n \Big( \sum_{k=1}^n \alpha_j^{2(k-1)} \Big) = \prod_{j=1}^n \frac{1-\alpha_j^{2n}}{1-\alpha_j^2},
\end{equation*}
where the last equality holds since we assumed that $f(1)f(-1) \neq 0$.
Since $$1-\alpha_j^2 = (1-\alpha_j)(1+\alpha_j)$$ we have $$\prod\limits_{j=1}^n \big|{1-\alpha_j^2}\big| = \abs{f(1)f(-1)}.$$ We conclude that 
\begin{equation*}
\abs{f(1)f(-1)}\Delta(f) \leq \prod_{j=1}^n \big|\alpha_j^{2n}-1\big| \overset{(\ast)}{<} \prod_{j=1}^{2n} \max \{1, \alpha_j^{2n}\} = \prod_{j=1}^{2n} \max \{1, \abs{\alpha_j}\}^{2n} = \MM(f)^{2n},
\end{equation*}
where the inequality $(\ast)$ follows upon separating the cases $\alpha_j^{2n} \geq 2$ {(so that $|\alpha_j^{2n}-1\big| < \alpha_j^{2n}$) and $ \alpha_j^{2n} < 2$ (in which case $|\alpha_j^{2n}-1\big| < 1$).}
\end{proof}

We will use Lemma~\ref{lem:Cn real roots} and Lemma~\ref{lem:disc Mahler bound} in steps \eqref{alg-v} and \eqref{alg-vii} of the algorithm. 

\subsection{Third tool: the conductor and class field theory}
\label{subsec:Ramification theory}
By the Kronecker--Weber theorem, every abelian number field is contained in some cyclotomic field. The \emph{conductor} $\mf{f}$ of an abelian extension $K/\Q$ is the smallest integer such that $K \subset \Q(\zeta_{\mf{f}})$. The next lemma collects (probably well-known) results concerning the set of values of the conductor and its relation to the field discriminant in our set-up of a cyclic extension of prime order. 

\begin{lemma}
\label{lem:class field theory}
Let $q$ be an odd prime and $K/\Q$ a cyclic field extension of degree $q$ and conductor $\mf{f}$. Then 
\begin{enumerate} \item \label{cft1} $\Delta_K = \mf{f}^{q-1}$;
\item \label{cft2} there is a $j \in \{0,2\}$ such that $\mf{f} = q^j \mf{f}_0$, where $\mf{f}_0$ is a product of distinct primes, each of which is congruent to $1 \bmod{q}$ and divides $\mf{f}_0$ exactly once.
\end{enumerate}
\end{lemma}
\begin{proof} We prove \eqref{cft1} locally for every prime $p$. The $p$-valuation of the discriminant is given by the \emph{F\"uhrerdiskriminantenproduktformel} (see, e.g., \cite[VI.4.3 \& VI.4.4]{CF})
\begin{equation}
\label{eq:p-valuation of field discriminant}
v_p(\Delta_K) = \sum_{\chi \in \widehat{\Z/q}} \mf{f}(\chi, p),
\end{equation}
where the local Artin conductor $\mf{f}(\chi, p)$ is given by 
\begin{equation}
\label{eq:local Artin conductor}
\mf{f}(\chi, p) \coloneqq  \sum_{j \geq 0} \frac{\#{I_j}}{\#{I_0}} (1 - \chi(I_j));
\end{equation}
here, $I_j$ is the $j$-th ramification group of a prime $\mf{p}$ over $p$, and $\chi(I_j)$ is the average value of $\chi$ on $I_j$. In our setup where $\Z/q$ is cyclic of prime order, $\chi(I_j)$ is $1$ if $\chi$ is trivial on $I_j$ and zero otherwise; and $I_j \leq \Z/q$ is either the trivial group or the full group $\Z/q$. Thus, $\mf{f}(\chi, p)$ is the smallest $j$ for which $I_j$ is trivial, and \eqref{eq:p-valuation of field discriminant} simplifies to $v_p(\Delta_K) = (q-1)\mf{f}(\chi, p)$. Since the conductor $\mf{f}$ (as we have defined it) is equal to the global Artin conductor $\mf{f}(\chi) \coloneqq  \prod_p p^{\mf{f}(\chi, p)}$ \cite[Chapter~VI.6]{Neukirch}, we find that $v_p(\Delta_K) = (q-1)v_p(\mf{f})$ holds for all $p$, proving the first part of the statement.

To prove \eqref{cft2}, notice that the values of $\mf{f}$ such that $\Q(\zeta_\mf{f})/\Q$ has a subextension $K$ with Galois group $\Z/q$ are precisely the $\mf{f}$ with $q \mid \phi(\mf{f})$---in other words, the $\mf{f}$ divisible by $q^2$ or by a prime $p \equiv 1 \bmod{q}$. Since the wild inertia subgroup $I_1 = I_1(p)$ is a $p$-group \cite[I.8, Theorem~1(ii)]{CF}, $I_1$ is trivial if $p \neq q$, in which case $v_p(\mf{f}) = 1$ (if $p$ ramifies) or $0$ (if $p$ does not ramify).

Finally, we consider the prime $q$. If $q$ does not ramify, then $v_q(\mf{f}) = 0$, whereas the above shows that if $q$ ramifies, then $v_q(\mf{f}) \geq 2$. That $v_q(\mf{f})$ is in fact equal to $2$ can be derived from the following more precise local computation. By class field theory, the valuation $v_q(\mf{f})$ equals the smallest $j \geq 0$ such that the $j$-th unit group $U_{\Q_q}^j = 1 + q^j\Z_q$ is contained in $N_{K_{\mf{q}}/\Q_q}(K_{\mf{q}}^*)$, i.e., any element congruent to $1 \bmod{q^j}$ is a norm of some invertible element of the localisation of $K$ at a prime $\mf{q}$ extending $q$ \cite[Chapter~V, Definition~1.6]{Neukirch}. Take $a \in \Z_q$; then $$N_{K_{\mf{q}}/\Q_q}(1+aq) = (1+aq)^q = 1 + bq^2$$ for some $b \in \Z_q$. 
To show $v_q(\mf{f}) = 2$, it suffices to show that any $b$ occurs as we vary $a$. The proof relies on Hensel's lemma. Note that $1+aq^2 \equiv 1+bq^2 \bmod{q^3}$, so $a \equiv b \bmod{q}$. Write $\FF_q \ni b_0 \coloneqq  b \bmod{q}$. Set $A = (a-b_0)/q, B = (b-b_0)/q \in \Z_q$. Then solving $(1+aq)^q = 1+bq^2$ for $a \in \Z_p$ is equivalent to finding a zero $A \in \Z_p$ of the polynomial
\begin{equation*}
F(x) \coloneqq  x-B + \sum_{j=2}^q \binom{q}{j} (qx+b_0)^jq^{j-3}.
\end{equation*}
Modulo $q$, the polynomial $F$ vanishes at $B + b_0^2/2 \bmod{q}$ (recall that $q$ is an odd prime); this is a simple zero since $F'(x) \equiv 1 \bmod{q}$. By Hensel's lemma, this zero lifts to a zero $A \in \Z_p$ of $F$.
\end{proof}

Collecting equalities and inequalities from the previous three lemmas yields the following conductor-Mahler measure inequality.
\begin{lemma}
\label{lem:conductor Mahler bound}
Let $f \in \Z[x]$ be a monic polynomial of odd prime degree $q$ with Galois group $\Z/q$, splitting field $K$ of conductor $\mf{f}$, and Mahler measure bounded above by some real number $B$. Then $\mf{f} < B^{2q/(q-1)}$.  
\end{lemma}
\begin{proof} Since the polynomial discriminant of a primitive element of $K/\Q$ is always lower bounded by the field discriminant of $K$, by combining 
Lemma~\ref{lem:class field theory}, Lemma~\ref{lem:Cn real roots}, and Lemma~\ref{lem:disc Mahler bound} we find  
\[
\mf{f}^{q-1} = \Delta_K \leq \Delta(f) \leq \abs{f(1)f(-1)}\Delta(f) < \MM(f)^{2q} \leq B^{2q}, 
\] 
which implies the desired inequality. 
\end{proof}

We will use the bound from Lemma \ref{lem:conductor Mahler bound} as well as the restrictions on the factorisation of the conductor from Lemma \ref{lem:class field theory} in step \eqref{alg-i} and \eqref{alg-ii} of the algorithm. 

\subsection{Fourth tool: splitting types modulo primes} 

The last tool is to use the possible splitting types of $f$ over finite fields, which are severely restricted since the Galois group is cyclic. Most of this follows from Dedekind's theorem on the relation between splitting of primes and the factorisation of $f$ modulo primes; only for primes dividing the polynomial discriminant but not the field discriminant, some extra work is necessary. 

\begin{lemma}
\label{lem:merging}
Let $f \in \Z[x]$ be a monic polynomial of odd prime degree $q$ with splitting field $K$ and Galois group $\Z/q$. For a prime $p$, denote $f_p \coloneqq (f \bmod{p}) \in \FF_p[x]$. Then one of the following holds:
\begin{enumerate}
    \item (irreducible) $p \nmid \Delta(f)$ and $f_p$ is irreducible, or
    \item (totally split) $p \nmid \Delta(f)$ and $f_p$ is the product of $q$ distinct linear factors, or
    \item \label{totsplit} (totally ramified) $p \mid \Delta_K$ and $f_p = (x+a)^q$ for some $a \in \FF_p$, or
    \item (\emph{au{\ss}erwesentliche Diskriminantenteiler}) $p \mid \Delta(f)$ but $p \nmid \Delta_K$ and $f_p$ is the product of $q$ linear factors of which at least two are the same.
\end{enumerate}
\end{lemma}
\begin{proof}
The Galois group of $f$ only consists of $q$-cycles and the identity, which is of cycle type $(1,1,...,1)$. By Dedekind's theorem, if $p \nmid \Delta(f)$ but $f_p$ is reducible, it follows that $p$ must be totally split. This covers cases (i) and (ii). For the case that $p \mid \Delta(f)$, \cite[Lemma~11.3]{BKK} implies, given the aforementioned cycle types, that $f_p$ is the product of $q$ linear factors of which at least $2$ are the same.
If $p \mid \Delta_K$, this can be strengthened further. Indeed, then $p$ ramifies in $K$, so the inertia subgroup $I_0$ at any prime ideal $\mf{p}$ over $p$ is nontrivial, so it must be $\Z/q$ since $q$ is prime. Since the decomposition group $G_p$ satisfies $I_0 \leq G_p \leq \Gal(f)$, also $G_p = \Z/q$; this is a transitive group, so $f$ is also irreducible over $\Q_p$, with splitting field $K_{\mf{p}}$. Let $\alpha \in \mc{O}_{K_{\mf{p}}}$ be a root of $f$ over $\Q_p$, and let $a \in k$ be the image of $-\alpha$ under the canonical projection of $\mc{O}_{K_{\mf{p}}}$ to the residue class field $k$ of $K_{\mf{p}}$. Since $I_0 = G_p$, we have in fact $k = \FF_p$. By \cite[Chapter~7, Exercise~7, p.~139]{Cassels}, it follows that $f_p(x) = (x+a)^q$.
\end{proof}


Lemma \ref{lem:merging} will be used in step \eqref{alg-ii} of the algorithm. 

\subsection{Fifth tool: higher congruences in the totally, wildly ramified case}

In the totally ramified case, Lemma~\ref{lem:merging} severely restricts the factorisation, which is very helpful in the algorithm. However, when the conductor of $K$ is simply $q^2$, in which case the extension is totally and wildly ramified, it leaves many potential candidates; without additional restrictions, the computation of, e.g., $m_7$ is infeasible. In the next lemma, we show that there are further congruences modulo $q^2$.

\begin{lemma}
	\label{lem:totally wild}
	Let $f \in \Z[x]$ be a monic polynomial of odd prime degree $q$ with Galois group $\Z/q$ and splitting field $K$ of conductor divisible by $q^2$. Fix any $a \in \Z$ for which $f(x) \equiv (x-a)^q \bmod{q}$ and set $h(x) \coloneqq f(x+a) \in \Z[x]$. Write $h(x) = x^q + \sum_{j=0}^{q-1} b_j x^j$. Then $q^2 \mid b_j$ for $j=1, \ldots, q-2$. 
	
	Let $v_q$ denote the standard $q$-adic valuation on $\Q$ and set $c \coloneqq v_q(b_0) \in \Z_{\geq 1}$. 
	More precisely, for any $j = 1, \ldots, q-1$, we have
	\begin{equation*}
		v_q(b_j) \geq 
		\begin{cases}
			1 + \left\lceil \frac{c(q-j)-1}{q} \right\rceil \quad &\text{if } q \nmid c, \\
			q-j \quad &\text{if } q \mid c.
		\end{cases}
	\end{equation*}
\end{lemma}
\begin{proof}
	We start with some preparatory work. By Lemma~\ref{lem:class field theory}(i), the assumption that the splitting field $K$ has conductor divisible by $q^2$ implies that $q^{2q-2} \mid \mid \Delta_K$. Lemma~\ref{lem:merging}(iii) then implies the existence of an integer $a$ as in the lemma statement. Furthermore, the proof of Lemma~\ref{lem:merging}(iii) implies that $f$ remains irreducible over $\Q_q$.
	Thus the polynomial $h$ of the lemma statement is irreducible over $\Q_q$ as well. It also satisfies $h(x) \equiv x^q \bmod{q}$. Define $d_j \coloneqq b_j/q \in \Z$ for $j = 1, \ldots, q-1$, and
	\begin{equation*}
		g(x) \coloneqq \frac{h'(x)}{q} = x^{q-1} + \sum_{j=1}^{q-1} jd_j x^{j-1} \in \Z[x].
	\end{equation*}
	Let $\alpha$ be a root of $h$ and let $\mfq$ be the unique prime ideal above $q$. Then $K_{\mfq} = \Q_q(\alpha)$ is the splitting field of $h$ over $\Q_q$. Denote by
	\begin{equation*}
		N(x) \coloneqq N_{K_{\mfq}/\Q_q}(x) \coloneqq \prod_{\sigma \in \Gal(K_{\mfq}/\Q_q)} \sigma(x)
	\end{equation*}
	the norm from $K_{\mfq}$ to $\Q_q$, and by $v$ the extension to $K_{\mfq}$ of the standard $q$-adic valuation $v_q$ on $\Q$, which for $x \in K_{\mfq}^{\times}$ is given by
	\begin{equation*}
		v(x) \coloneqq \frac{1}{q} v_q(N(x)).
	\end{equation*}
	In particular, we have
	\begin{equation}
		\label{eq:v(a)}
		v(\alpha) = \frac{1}{q} v_q((-1)^{q}b_0) = \frac{c}{q}.
	\end{equation}
	Denote by $\alpha_1, \ldots, \alpha_q$ the roots of $h$, and by
	\begin{equation*}
		\Delta(h) = \prod_{1 \leq i < j \leq q} (\alpha_i - \alpha_j)^2 = (-1)^{q(q-1)/2} N_{K/\Q}(h'(\alpha))
	\end{equation*}
	the discriminant of $h$, where $N_{K/\Q}$ denotes the norm from $K$ to $\Q$.

	We begin with the case that $q \mid c = v_q(b_0)$. Consider the monic, irreducible polynomial $H(x) \coloneqq h(qx)/q^q \in \Q[x]$. Then $H$ is the minimal polynomial of $\alpha/q$ over $\Q$. Since $v(\alpha)$ is a positive integer, we have $v(\alpha/q) \geq 0$, so $\alpha/q \in \mc{O}_{K_{\mfq}}$, implying that $H \in \Z[x]$. Therefore $q^q$ divides the coefficients of $h(qx)$, so $q^q$ divides $b_j q^{j}$ for all $j = 0, \ldots, q-1$. Hence $v_q(b_j) \geq q-j$ for all $j$.
	
	From here on, we assume that $q \nmid c = v_q(b_0)$. With the goal to obtain an expression for $v(g(\alpha))$, we compute $v_q(\Delta(h))$ in two ways. Firstly, since $h'(\alpha) = qg(\alpha)$, we have
	\begin{equation}
		\label{eq:v_q-1}
		v_q(\Delta(h)) = v_q(N_{K/\Q} (h'(\alpha))) = v_q(N_{K_{\mfq}/\Q_q} (h'(\alpha))) = q + v_q(N(g(\alpha))).
	\end{equation}
	For the second computation of $v_q(\Delta(h))$, we start by claiming that $v(\alpha_i-\alpha_j) = (c+1)/q$ for any distinct $i$, $j$. Assuming the claim, we find
	\begin{equation}
		\label{eq:v_q-2}
		v_q(\Delta(h)) = 2 \sum_{1 \leq i < j \leq q} v(\alpha_i - \alpha_j) = 2 \cdot \frac{q(q-1)}{2} \cdot \frac{c+1}{q} = (q-1)(c+1).
	\end{equation}
	Combining \eqref{eq:v_q-1} and \eqref{eq:v_q-2} yields $q + q v(g(\alpha)) = v_q(\Delta(h)) = cq + q - c -1$,
	which implies
	\begin{equation}
		\label{eq:v(g(a))-explicit}
		v(g(\alpha)) = \frac{cq-c-1}{q}.
	\end{equation}
	We now prove the claim that $v(\alpha_i - \alpha_j) = (c+1)/q$ for all distinct $i$, $j$. In the proof of Lemma~\ref{lem:class field theory}(ii), we saw that the smallest $\ell$ such that the $\ell$-th unit group $U_{\Q_q}^\ell$ is contained in $N_{K_{\mfq}/\Q_q}(K_{\mfq}^{\times})$ is $\ell=2$. By local class field theory, this means that $v(\sigma(x)-x) \geq 2/q$ for all $x \in \mathcal{O}_{K,q}$, and that $v(\sigma(\pi) - \pi) = 2/q$ for any uniformizer $\pi$ of $K_{\mfq}$. Fix $\alpha_j$. Since $v(\alpha_j) = c/q$, there exists a unit $u \in \mathcal{O}_{K,q}^{\times}$ such that $\alpha_j = u \pi^c$. By transitivity, there exists $\sigma \in \Gal(K_{\mfq}/\Q_q)$ such that $\sigma(\alpha_j) = \alpha_i$. Hence 
	\begin{equation}
		\label{eq:aiaj}
		\alpha_i - \alpha_j =  \sigma(u) \sigma(\pi)^c - u \pi^c = \sigma(u) (\sigma(\pi)^c - \pi^c) + \pi^c(\sigma(u) - u).
	\end{equation}
	Recall that the ultrametric inequality implies, for all $x, y \in K_{\mfq}^{\times}$, that
	\begin{equation}
		\label{eq:v(x+y)}
		v(x+y) \geq \min\{v(x), v(y)\} 
	\end{equation}
	with equality if $v(x) \neq v(y)$.
	We now compute the valuations of the terms on the right-hand side of \eqref{eq:aiaj}. Since $v(\pi^c) = c/q$ and $v(\sigma(u)-u) \geq 2/q$, we have
	\begin{equation}
		\label{eq:secondterm}
		v(\pi^c(\sigma(u) - u)) \geq (c+2)/q.
	\end{equation}
	Write the first term on the right-hand side of \eqref{eq:aiaj} as
	\begin{align*}
		\sigma(u)(\sigma(\pi)^c - \pi^c) = \sigma(u) \pi^{c-1} ( \sigma(\pi)-\pi) \sum_{\ell=0}^{c-1} \Big( \frac{\sigma(\pi)}{\pi} \Big)^{\ell}. 
	\end{align*}
	Then $v(\pi^{c-1}) = (c-1)/q$ and $v(\sigma(\pi)-\pi) = 2/q$, and $v(\sigma(u)) = 0$ since $u$ is a unit. Lastly, since $\sigma(\pi)/\pi \equiv 1 \bmod{\mathfrak{q}}$, we have
	\begin{equation*}
		\sum_{\ell=0}^{c-1} \Big( \frac{\sigma(\pi)}{\pi} \Big)^{\ell} \equiv c \not \equiv 0 \bmod{\mathfrak{q}},
	\end{equation*}
	where we used the assumption that $q \nmid v_q(b_0)$ in the last step.
	Therefore
	\begin{equation*}
		v(\sigma(u)(\sigma(\pi)^c - \pi^c)) = \frac{c+1}{q}.
	\end{equation*}
	Combining this with \eqref{eq:aiaj}, \eqref{eq:v(x+y)} and \eqref{eq:secondterm} yields the claimed $v(\alpha_i-\alpha_j) = (c+1)/q$ and thus concludes the proof of \eqref{eq:v(g(a))-explicit}.
	
	To finish the proof, we will again use \eqref{eq:v(x+y)} to analyze $v(g(\alpha))$. Note that $v(jd_j \alpha^{j-1}) = v(d_j) + c(j-1)/q$ for $j = 1, \ldots, q-1$, and $v(\alpha^{q-1}) = c(q-1)/q$; furthermore, $v(d_j) \in \Z$. Since $v(\alpha) \not \in \Z$, the valuations of the monomials of $g$ are distinct, since they are already distinct modulo $1$. Combining this with \eqref{eq:v(g(a))-explicit} and  \eqref{eq:v(x+y)}, we conclude that
	\begin{equation*}
		v(d_j) \geq \left\lceil \frac{cq-c-1-c(j-1)}{q} \right\rceil = \left\lceil \frac{c(q-j)-1}{q} \right\rceil
	\end{equation*}
	for $j = 1, \ldots, q-1$. 
\end{proof}


\subsection{The algorithm}
\label{subsec:The algorithm}
We describe a general algorithm to find the smallest Mahler measure for an integral polynomial with fixed cyclic Galois group, possibly with specific ramification constraints. In the algorithm, the verification that the Galois group is cyclic is postponed to one of the final steps, and the first steps use the bounds and factorisation restrictions from the previous subsections. 

\bigskip 

\hrule

\medskip

\noindent \underline{\texttt{Input}} $q$, an odd prime; $B>0$ a positive real number. \\[2mm]
\underline{\texttt{Output}} A set of monic polynomials $f\in \Z[x]$ of degree $q$ with Galois group $\Z/q$ representing all possible values of $\MM(f) < B$. \\[2mm]
\underline{\texttt{Procedure}} 
\begin{enumerate}
\item \label{alg-i} Let $K$ be the splitting field of such an $f$, and let $\mf{f}$ be its conductor. By Lemma~\ref{lem:conductor Mahler bound} we have $\mf{f} < B^{2q/(q-1)}$. Lemma~\ref{lem:class field theory} further restricts what values of $\mf{f}$ can occur. Collect such $\mf{f}$ in a set $\mf{F}$. 
\item\label{alg-ii}  Select $\mf{f} \in \mf{F}$. 
Let $p$ be a prime divisor of $\mf{f}$. Let $f$ be a candidate polynomial. Then $p \mid \Delta_K$ by Lemma~\ref{lem:class field theory}, so by Lemma~\ref{lem:merging} we have $f_p = (x+a)^q$ for some $a \in \FF_p$; in addition, by the coefficient bounds of Lemma~\ref{lem:Mahler coefficient bound}, we need $f_p(0) = a^q$ to have a lift in $\Z_{>0}$ smaller than $B$. (Indeed, after replacing $f$ by $-f(-x)$, which is still monic and has the same splitting field and Mahler measure as $f$, we may always assume that the constant coefficient of $f$ is positive.) Denote by $\mf{f}_1$ the product of the distinct prime divisors of $\mf{f}$. Varying over the prime divisors of $\mf{f}$ and using the Chinese remainder theorem, we end up with the set $S \subset \Z/\mf{f}_1\Z[x]$ of at most $\mf{f}_1$ polynomials characterized by the property that $h \in S$ if and only if for all $p \mid \mf{f}$ there exists an $a \in \FF_p$ such that $h_p = (x+a)^q$ and $a^q$ lifts to a positive integer smaller than $B$.

\item \label{alg-iii} For every $h \in S$, we check if the reversal $x^q h(1/x)$ of $h$ is also in $S$. If this is the case, then remove one of the two from $S$. (This simplification is possible because the coefficient bounds of Lemma~\ref{lem:Mahler coefficient bound} are symmetric, and the reversal of $h$ has the same Mahler measure and splitting field as $h$.)

\item \label{alg-iv} Select $h \in S$. Determine all monic lifts of $h$ from $\Z/\mf{f}_1\Z[x]$ to $\Z[x]$ with positive constant coefficient satisfying the coefficient bounds from Lemma~\ref{lem:Mahler coefficient bound}. When $\mf{f} = q^2$, these lifts should also satisfy the congruences from Lemma~\ref{lem:totally wild}. This gives a set $S'$ of candidate polynomials over $\Z$.

\item \label{alg-v} For each $f \in S'$, verify if $\Delta(f)$ is a square, satisfies $0 < \abs{f(1)f(-1)}\Delta(f) < B^{2q}$, and if $\mf{f}^{q-1}$ divides $\Delta(f)$ (which is a multiple of the discriminant of the number field defined by $f$).

\item \label{alg-vi} For each surviving $f$, check if it is irreducible.

\item \label{alg-vii} For each surviving $f$, check if it has only real roots and, if so, whether $1 < \MM(f) < B$.

\item \label{alg-viii} Calculate the Galois group of $f$ if $f$ survives up to this point and add $f$ to the output set if the Galois group is $\Z/q$.

\item \label{alg-ix} Repeat from step \eqref{alg-iv} onwards until $S$ is exhausted.

\item \label{alg-x} Repeat from step \eqref{alg-ii} onwards until $\mf{F}$ is exhausted.
\end{enumerate}
\hrule 

\bigskip

In particular, if the output of the algorithm is empty and the input $B$ is the Mahler measure of such a polynomial itself, then $B=m_q$. 

\subsection{Implementation} The above algorithm was implemented in \texttt{SageMath} (see Section~\ref{sec:code}) in degree $3$, $5$ and $7$ and used to prove Theorem~\ref{mainMM}, i.e., that $$m_q  = \m(\alpha_q) \qquad \text{for } q=3,5,7.$$  
To illustrate the procedure for the most challenging case $n=7$, the list of remaining conductors (given by their prime factorisations) is 
\begin{align*} \mf{F} =
\{ &29, 43, 7^2, 71, 113, 127, 197, 211, 239, 281, 337, 343, 379, 421, 449, 463, 491, 547, 617, \\ 
& 631, 659, 673, 701, 743, 757, 827, 883, 911, 953, 967, 1009, 1051, 1093, 1163, 29 \cdot 43, \\
& 1289, 1303, 1373, 7^2 \cdot 29, 1429, 1471, 1499, 1583, 1597, 1667 \}.  
\end{align*} 
In general, the list of polynomials to check is larger for small or powerful conductor (due to the congruences excluding fewer polynomials).  
To give an indication of the running time, it took around a week to test all approximately $10^{10}$ polynomials for conductor $29$ on a standard user-grade laptop (Intel i5, 8GB RAM).   

\begin{remark}
	When $q=7$ and we denote by $\zeta$ a primitive $49$-th root of unity, the field $K$ of conductor $49$ and degree $7$ over $\Q$ (which is unique by Kronecker--Weber) is $K = \Q(\eta)$ where
	\begin{equation*}
		\eta = \zeta + \zeta^{18} + \zeta^{19} + \zeta^{30} + \zeta^{31} + \zeta^{48}
	\end{equation*}
	is the corresponding Gaussian period.
	The minimal polynomial of $\eta$ is $f_{\eta}(x) \coloneqq x^7 - 21 x^5 - 21 x^4 + 91 x^3 + 112 x^2 - 84 x - 97$, with Mahler measure $M(f_\eta(x)) \approx 109.7392$. The polynomial $f_\eta(x-1) = x^7 - 7x^6 + 49x^4 - 98x^2 - 49x+7$ has smaller Mahler measure, namely $\approx 68.7528$. The minimal Mahler measure among primitive elements of $K$, however, is $\approx 60.1076$ and is attained by the polynomial 
	\begin{equation*}
		f_{49}(x) = x^7 + 7x^6 - 21x^5 - 49x^4 + 91x^3 + 7x^2 - 56x + 19.
	\end{equation*}
	This is the minimal polynomial of
	 $-\zeta^{36} - \zeta^{34} + \zeta^{33} + \zeta^{31} - \zeta^{27} - \zeta^{25} - \zeta^{23} - \zeta^{22} + \zeta^{19} + \zeta^{17} - \zeta^{15} - \zeta^{13} + \zeta^{12} - \zeta^{11} + \zeta^{10} - \zeta^9 + \zeta^5 - \zeta^4 + \zeta^3 - \zeta^2 - 1$.
\end{remark}

\subsection{Computation of other small Mahler measures} 

In Table \ref{tab:Mn}, we have listed some more values of $m_n$ for (very) small $n$. For general $n$, it is not true that $m_n=m(\alpha_n)$; this fails for $n=4$, for example. We do not have a conjectured value for $m_n$ in case $n$ is even. 

{\small \begin{table}[ht!]
\centering
\begin{tabular}{cccc}
\toprule
$n$ 	& $M_n$ 	& $\beta_n$ & $f_n$     \\ \midrule
1 		& $2$ 		& $2$   & $x-2$	  	 \\
2 		& $1.618$	& $2\cos(2\pi/5)$      & $x^2+x-1$	      \\
3 		& $2.247$	& $2\cos(2\pi/7)$        & $x^3 + x^2 - 2x - 1$  	 \\
4 		& $2.618$	& $\zeta_5+\zeta_5^3 \notin \R$ &	$x^4 + 2x^3 + 4x^2 + 3x + 1$ \\
5 		& $4.229$	& $2 \cos(2 \pi/11)$	& $x^5 + x^4 - 4x^3 - 3x^2 + 3x + 1$	\\
7 		& {$24.217$} 		& {$2\cos(\frac{2 \pi}{29})+2\cos(\frac{34 \pi}{29})$} 	& {$x^7 + x^6 - 12x^5 - 7x^4 + 28x^3 + 14x^2 - 9x + 1$}	 \\ \bottomrule
\end{tabular}
\caption{The value of $M_n = \exp(m_n)$ for $n \leq 5$ and $n=7$, rounded to three decimal places, an algebraic integer $\beta_n$ and its cyclic degree-$n$ minimal polynomial $f_n \in \Z[x]$ such that $\MM(f_n) = M_n$.}\label{tab:Mn}
\end{table}
}

Through direct search, we have found some small upper bounds on $m_n$ for a few other values of $n$, which we list in Table \ref{tab:uppMn}. There can be several polynomials (that are not each others reversal) whose Mahler measure is the candidate minimum; this happens, for example, for $n=6$ and $n=10$. We list only one of the possibilities. 

{\footnotesize \begin{table}[ht!]
\centering
\begin{tabular}{cccc}
\toprule
$n$ 	& $M_n \leq$ 	& $\beta_n$ & $f_n$     \\ \midrule
6 & {$5.049$} 	&  $\zeta_7+\zeta_7^3 \notin \R$ 	&	{$x^6 + 2x^5 + 4x^4 + x^3 + 2x^2 - 3x + 1$}		 \\
8 & {$11.107$} & {$2 \cos(2\pi/17)$}	&	{$x^8 + x^7 - 7x^6 - 6x^5 + 15x^4 + 10x^3 - 10x^2 - 4x + 1$} \\
9 		& {$15.350$} 		& {$2 \cos(2 \pi/19)$}  	&	{$x^9 + x^8 - 8 x^7 - 7 x^6 + 21 x^5 + 15 x^4 - 20 x^3 - 10 x^2 + 5 x + 1$}		 \\
10 		& {$17.882$}		&$\zeta_{11}+\zeta_{11}^3 \notin \R$	& {$x^{10} + 2x^9 + 4x^8 + 8x^7 + 16x^6 + 10x^5 + 20x^4 + 7x^3 + 3x^2 - 5x + 1$} \\ \bottomrule
\end{tabular}
\caption{The smallest upper bound we found for $M_n = \exp(m_n)$ for $n=6$ and $8 \leq n \leq 10$, rounded to three decimal places, an algebraic integer $\beta_n$ and its cyclic degree-$n$ minimal polynomial $f_n \in \Z[x]$ realising the upper bound. }\label{tab:uppMn}
\end{table}
}

\section{Conjectural asymptotics of $m_n$; Conjecture~\ref{mainconj}} 
\label{sectconj}

\subsection{A conjecture on the growth of $\m(\alpha_n)$} 

In Section \ref{sectalg}, we have shown some results that indicate that $m_n = \m(\alpha_n)$ for odd $n$. Here we discuss Conjecture~\ref{mainconj}, which concerns the growth rate of $\m(\alpha_n)$ as a function of $n$. This conjecture is part of a larger framework of height growth in families, as explained in the following two remarks. 

\begin{remark} \label{RemGalGrow} The current best lower bound on $\m(\alpha)$, where $\Q(\alpha)/\Q$ is any Galois extension of degree $n$, tends to zero in $n$ (see, e.g., Amoroso--Masser \cite{AmorosoMasser}). 
Consider the following general question: given a structured collection $\{G_n\}_{n=1}^\infty$ of finite groups of increasing size, what is the `average asymptotics' of 
$$ \min \{ \m(\alpha) \colon \alpha \in \overline \Z - \mu_\infty, \Gal(\Q(\alpha)/\Q) \cong G_n\} $$ as $n$ increases, depending on the properties of the groups $G_n$? Amoroso and Zannier \cite[Corollary~1.3]{Amoroso-Zannier} show a lower bound linear in $n$ when $G_n = D_n$ is the dihedral group of order $2n$. 
The problem seems interesting also for $G_n=S_n$ and $n>2$. The minimal polynomial is then not reciprocal, and Smyth has proven that the minimum over all varying $n$ is attained by the roots of $x^3-x-1$ \cite{SmythBLMS}, but the above question on the asymptotic dependence on $n$ appears to be open also in this case. Amoroso \cite{Amoroso} has proven that $S_n$-extensions admit generators for which the Mahler measure grows faster than linear in $n$, and conjectured that this holds for any generator of such extensions. For related work on $A_n$, see the recent \cite{Jenvrin}.  Our Conjecture~\ref{mainconj} is a quantitative version of this result for the collection of cyclic Galois groups $G_n = \Z/n\Z$. 
\end{remark} 

We start with the following simple, unconditional result.
\begin{lemma}
	\label{lem:m(an) unconditional bounds}
	For all $n$, we have $n \ll m(\alpha_n) \ll n \log{n}$.
\end{lemma}
\begin{proof}
The lower bound follows from \eqref{eq:Schinzel}. For the upper bound, we use the observation --- already made in the proof of Theorem~\ref{mainas} --- that $\alpha_n$ and all its Galois conjugates have absolute value at most $\kappa(n)$. Hence $\m(\alpha_n) \leq n \log{\kappa(n)}$. Lastly, the bound $\log{\kappa(n)} \ll \log{n}$ follows from a result by Linnik, who proved that $\kappa(n) \ll n^{L-1}$ for some constant $L>1$ (see, e.g., \cite{Xylouris}).
\end{proof}

To improve conditionally on this result for odd $n$, in this section, we will use the following convenient terminology. 
\begin{definition}
We say that a property of odd integers $n$ holds \emph{almost always} if it holds up to a subset of the odd integers of natural density zero, i.e., it holds for all but $o(X)$ odd values $n \leq X$ as $X \rightarrow + \infty$. 
\end{definition} 

\begin{remark} \label{RemLauGrow} In \cite[p.~461]{Boyd-range}, Boyd suggests that the minimal Mahler measure of all irreducible Laurent polynomials with integer coefficients and for which the Newton polytope has dimension $m$ tends to infinity with $m$. In our situation, this suggestion implies $\m(C_k) \rightarrow + \infty$ with $k$ (indeed, $x_0 + F_k(\mathbf x)$ is irreducible, since linear in $x_0$, and the dimension of the corresponding Newton polytope---including the variable $x_0$---is $\phi(k)+1$, which tends to infinity with $k$). We will see below that $k = \kappa(n)$ almost always grows with $n$, so that Boyd's suggestion would imply a weaker version of Conjecture~\ref{mainconj}, namely, $\m(\alpha_n)/n \rightarrow + \infty$ for a suitable sequence of values $n \rightarrow +\infty$. 
\end{remark}

The purpose of the remainder of this section is to formulate two other conjectures, related to random walks and to numerical integration, that together imply Conjecture~\ref{mainconj}. 
Up to this point, we have considered (asymptotic) properties of $\m(\alpha_n)$ and $\m(C_k)$ in a regime of values of $k$ that were either fixed (Section \ref{sec:asymp_m_alpha_n_fixed_k}) or at least controlled (Section \ref{sec:random_walks_and_asymp_mCk}). By studying the asymptotics in $n$ instead, such control on $k$ is no longer a given.

\subsection{How Conjecture~\ref{mainconj} follows from other conjectures}

The first of two conjectures that we introduce concerns the growth behaviour of the Mahler measure of the cyclovariety $C_k$, but as a function of $n$ (i.e., with $k=\kappa(n)$). 

\begin{conjecture}
\label{con:mCk logk}
For odd $n$, $\m(C_{\kappa(n)}) \gg \log{\kappa(n)}$ almost always.
\end{conjecture}

The second conjecture concerns the convergence of the Mahler measure $\m(\alpha_n)$ to $\m(C_{\kappa(n)}) \cdot n$ with $n$ tending to infinity. 

\begin{conjecture}
\label{con:malphan mCk}
For odd $n$, $\m(\alpha_n)/n \asymp \m(C_{\kappa(n)})$ almost always.
\end{conjecture}

We will also need the following  special case of the Bateman--Horn conjecture for linear polynomials in precise quantitative form (see \cite[p.~310]{Granville}). 

\begin{conjecture}[Specific case of prime tuplets conjecture]
\label{con:PT_with_error_term}
Fix a constant $b$ and integers $r_1, r_2, \ldots, r_\ell \leq b \log{x}$. For a prime $p$, define $w(p)$ to be the number of distinct solutions $q \bmod{p}$ of the congruence $\prod_{j=1}^\ell (qr_j+1) \equiv 0 \bmod{p}$. 
Then
\begin{equation*}
\# \{q : x \leq q < 2x, \text{ each } q r_j + 1 \text{ prime} \} = (1+o_{b,\ell}(1)) \frac{x}{(\log{x})^\ell} \prod_p \frac{1-w(p)/p}{(1-1/p)^\ell}
\end{equation*}
where the product extends over all primes.
\end{conjecture}

In the rest of this section, we will often leave out the `for odd $n$', assuming it tacitly. The following proposition explains the relation between the conjectures. 

\begin{proposition}
\label{prop:conjs imply conj}
Conjectures~\ref{con:mCk logk}, \ref{con:malphan mCk} and \ref{con:PT_with_error_term} imply Conjecture \ref{mainconj}. 
\end{proposition}
\begin{proof}
Together with Proposition~\ref{prop:upper-bound}, Conjecture~\ref{con:mCk logk} and Conjecture~\ref{con:malphan mCk} imply that $\m(\alpha_n) \asymp n \log{\kappa(n)}$ almost always. It thus remains to prove $\log{\kappa(n)} \asymp \log{\log{n}}$ almost always. 

For a lower bound, the following lemma is more than sufficient. 

\begin{lemma} 
\label{lem:k(n) lower bound}
For all but at most $O(X/\log{\log{X}})$ values of $n \leq X$, we have $$\kappa(n) > (\log{n})/\log{\log{n}}.$$
\end{lemma}
\begin{proof}
Let $S(X)$ be the cardinality of the set $\{1 \leq n \leq X : \kappa(n) \leq \log{X}/\log{\log{X}}\}$. By \eqref{eq:primes_in_AP} (i.e., Siegel--Walfisz), for fixed $1 \leq k \leq \log{X}/\log{\log{X}}$, the number of positive integers $1 \leq n \leq X$ such that $kn+1$ is prime is 
\begin{equation*}
\frac{kX}{\phi(k) \log{kX}} + O \Big( \frac{kX}{(\log{kX})^2} \Big) = O \Big( \frac{X}{\log{X}} \Big( \frac{k}{\phi(k)} + \frac{k}{\log{X}} \Big)  \Big) = O\Big( \frac{kX}{\phi(k)\log{X}} \Big),
\end{equation*}
where the implied constants are independent of $k$. Therefore 
\begin{equation} \label{sx}
S(X) \ll \frac{X}{\log{X}} \sum_{1 \leq k \leq \log{X}/\log{\log{X}}} \frac{k}{\phi(k)} \ll \frac{X}{\log{\log{X}}},
\end{equation} 
where we used that $\sum_{1\leq k \leq Y} k/\phi(k) = O(Y)$, see, e.g., \cite[Equation (2.12)]{Sitaramachandra}.
\end{proof}

The required almost always upper bound $\log{\kappa(n)} \ll \log{\log{n}}$ follows from a result of Granville \cite[Theorem 4]{Granville}, stating that under Conjecture~\ref{con:PT_with_error_term}, we have $\kappa(n) \leq \log(n) f(n)$ almost always, for any function $f$ on the integers for which $\lim\limits_{n \rightarrow + \infty} f(n) = + \infty.$ Taking, e.g., $f(n) = \log{n}$ suffices to conclude the proof. 
\end{proof}

\begin{remark}[About the size of $\kappa(n)$]
As already mentioned in the proof of Lemma~\ref{lem:m(an) unconditional bounds}, unconditional results on the quantity $\kappa(n)$ go back as far as Linnik, who proved that $\kappa(n) \ll n^{L-1}$ for some constant $L>1$, which can be taken to be $5$ by results of Xylouris \cite{Xylouris}, and almost always equal to $2$ by \cite{Bombieri-et-al}. This is far from what is believed to be optimal. Under the extended Riemann hypothesis, $\kappa(n) \leq n \log(n)^2$ \cite[Cor.~1.2]{LLS}. Wagstaff \cite{Wagstaff} has proposed a heuristics and numerics that support the idea that `typically', $\kappa(n) \approx \log n$. Granville's argument also shows that $\kappa(n) \ll \log n$ does not hold almost always.
\end{remark} 

\begin{remark} Replacing $\asymp$ with $\gg$ in Conjecture~\ref{con:malphan mCk} and Conjecture~\ref{mainconj}, the weakening of the former implies together with Conjecture~\ref{con:mCk logk} the weakening of the latter, without assuming Conjecture~\ref{con:PT_with_error_term}. Since Conjecture~\ref{con:PT_with_error_term} appears to be widely believed, we have worked under this assumption. 
\end{remark}

\subsection{Discussion of Conjecture~\ref{con:mCk logk} and \ref{con:malphan mCk}}

\subsubsection*{Concerning Conjecture~\ref{con:mCk logk}} 
Theorem~\ref{mainRW} implies Conjecture~\ref{con:mCk logk} if $n$ has the property that $\kappa(n) \in Q \coloneqq  \{ q, 2 q^r \colon q \text{ prime} \}$. For general odd $n$, the size of $\m(C_{\kappa(n)})$ is related to the expected distance of a random walk in the plane by $k=\kappa(n)$ unit steps with certain dependencies between the angles dictated by the cyclopolytope; there are $d+1=\phi(k)+1$ independent random steps, and then $k-\phi(k)$ deterministic steps. The number of deterministic steps can be large compared to the number of independent steps (e.g., for primorials, as in Remark \ref{remprimorial}) and it is not clear how to analyze the expected distance from the origin of the random walk in such a situation. 

\subsubsection*{Concerning Conjecture~\ref{con:malphan mCk}} 
Conjecture~\ref{con:malphan mCk} is about the quantitative convergence rate of $\m(\alpha_n)/n$ to $\m(C_k)$, but now with $k=\kappa(n)$. The upper bound for this rate in \eqref{was}, deduced using Wasserstein distance, is divergent in $n$ if we use the typical value $\kappa(n) \approx \log n$ and $d =\varphi(k) \approx k$.
Using the Hlawka--Koksma upper bound, we notice that the variation of $f$ seems to be of order exponential in $k  \approx \log n$ (see the crude estimate in  \eqref{upperVbad}), whilst numerics suggest the discrepancy to be of order $O(\log(p)^a/p)$ for some integer $a$ depending on $k$. In conclusion, to prove the conjecture unconditionally, either one of the Wasserstein or Hlawka--Koksma bound would need to be improved.

\section{Some open problems} \label{open} 

For the convenience of the reader, we collect some open problems.
\begin{enumerate}
\item Improve the integration error bound on $\big| \m(\alpha_n)/n - \m(C_k)\big|$.  
\item Is the polytope $N_k$ reflexive for all $k$? See Remark \ref{guessreflexive}. 
\item Is $C_k(\Cc) \cap \T^{d+1} \neq \emptyset$ for all $k$? See Subsection \ref{toricpointsDS}. 
\item Explore possible relations between $\m(C_k(\lambda))$ and Alexander polynomials of links. See Remark \ref{rem:links}. 
\item Compute $\m(C_9)$ up to ten correct digits, cf.\ Table \ref{tab:mCk}.
\item Find the first order asymptotics for $\m(C_k)$ in $k$ for all $k$, cf.\ Theorem \ref{mainRW}. 
\item Are the singularities of the Picard--Fuchs operators $\mathcal L_k$ always real for all $k$? See Remark \ref{realsing}.
\item Find a better algorithm that computes $m_n$ for $n \geq 9$, cf.\ Section \ref{sectalg}. 
\item Study Conjecture~\ref{con:mCk logk} and Conjecture~\ref{con:malphan mCk}. 
\item What happens in the general framework of Remark \ref{RemGalGrow} for other interesting families of groups, e.g., if $G_n$ is a collection of elementary abelian groups $(\Z/m)^n$ for fixed $m$; or Chevalley groups $G(\mathbf{F}_q^n)$ for some suitable reductive algebraic group $G$? The absolute minimum in Lehmer's original question is attained by a reciprocal polynomial of degree $10$ with Galois group $((\Z/2)^5 \rtimes S_5) \cap A_{10}$, and this suggests to consider the family $G_n = ((\Z/2)^n \rtimes S_n) \cap A_{2n}$ of Galois groups of generic reciprocal polynomials with square discriminant.  
\end{enumerate}

\bigskip

\section{Code}
\label{sec:code}

The Github repository 
\begin{quote} 
\url{https://github.com/davidhokken/cyclic-mahler} 
\end{quote} 
contains the following \texttt{SageMath}/Python routines: 
\begin{itemize} 
\item An algorithm to compute the cyclopolytopes and check reflexivity as described in Remark \ref{algocycsage}, as well as a check that the transportation polytope $T(3,5)$ is combinatorially, but not polar dual to the cyclopolytope $N_{3 \cdot 5}$. 
\item A workbook containing the calculations of the graphs of the density functions $\rho_k$ for $k=6,8,10$, and of the corresponding Mahler measures $\m(C_k)$. 
\item A workbook containing the routines to compute $m_n$ for small $n$, as described in Subsection \ref{subsec:The algorithm}.
\end{itemize} 

\bibliographystyle{amsplain}
\providecommand{\bysame}{\leavevmode\hbox to3em{\hrulefill}\thinspace}
\providecommand{\MR}{\relax\ifhmode\unskip\space\fi MR }
\providecommand{\MRhref}[2]{%
  \href{http://www.ams.org/mathscinet-getitem?mr=#1}{#2}
}
\providecommand{\href}[2]{#2}

\end{document}